\documentclass[11pt]{article}
\usepackage[margin=1in]{geometry}
\usepackage{graphics}
\usepackage{xcolor}
\usepackage{enumerate}
\usepackage{commath}
\usepackage{algorithm}
\usepackage{float}
\usepackage{algorithmic}
\usepackage{listings}
\usepackage{cancel} 
\usepackage{bm}
\usepackage{lmodern}
\usepackage[font=scriptsize, margin=1cm]{caption}

\def\cS{{\cal S}}

\newcommand{\EE}{\mathbb{E}}

\newcommand{\cX}{\mathcal{X}}
\newcommand{\cD}{\mathcal{D}}

\newcommand{\cV}{\mathcal{V}}

\newcommand{\RR}{\mathbb{R}}

\newcommand{\cF}{\mathcal{F}}

\newcommand{\cO}{\mathcal{O}}

\newcommand{\QED}{\ \hfill\rule[-2pt]{6pt}{12pt} \medskip}

\newcommand{\inner}[2]{\left\langle #1,#2 \right\rangle}
\newcommand{\cL}{\mathcal{L}}

\usepackage[shortlabels]{enumitem}

\usepackage{amsmath}
\usepackage{algorithm}
\usepackage{algorithmic}
\usepackage{appendix}

\def \Var {\mathrm{Var}}

\DeclareMathOperator*{\argmax}{argmax}
\DeclareMathOperator*{\argmin}{argmin}

\usepackage{natbib}

\bibpunct[, ]{(}{)}{,}{a}{}{,}%

\usepackage{tcolorbox}

\usepackage[hidelinks]{hyperref}       
\hypersetup{
    colorlinks=false
} 
\usepackage{url}            
\usepackage{booktabs}       
\usepackage{amsfonts}       
\usepackage{nicefrac}       
\usepackage{microtype}      

\usepackage{authblk}

\newtheorem{theorem}{Theorem}
\newtheorem{proposition}{Proposition}
\newtheorem{lemma}{Lemma}

\newtheorem{assumption}{Assumption}

\title{Stochastic Compositional Optimization with Compositional Constraints}

\author[1]{Shuoguang Yang}
\author[1]{Wei You}
\author[2]{Zhe Zhang}
\author[3]{Ethan X. Fang}

\affil[1]{Department of Industrial Engineering and Decision Analytics, The Hong Kong University of Science and Technology}
\affil[2]{School of Industrial Engineering, Purdue University}
\affil[3]{Department of Biostatistics and Bioinformatics, Duke University}

\begin{document}

\maketitle

\begin{abstract}
	Stochastic compositional optimization (SCO) has attracted considerable attention due to its broad applicability to important real-world problems. 
	However, existing work on SCO typically assumes that the projection within a solution update is straightforward, which is not the case for problem instances where constraints are in the form of expectations, such as empirical conditional value-at-risk constraints. 
	In this paper, we introduce a novel model that integrates single-level expected value and two-level compositional constraints into the existing SCO framework. 
	Our model has wide applicability to data-driven optimization, fairness optimization, and risk management, including risk-averse optimization and high-moment portfolio selection, and is capable of handling multiple constraints. 
	Additionally, we propose a class of primal-dual algorithms that generate sequences converging to the optimal solution at a rate of $\cO(\frac{1}{\sqrt{N}})$ under both single-level and two-level compositional expected-value constraints, where $N$ is the iteration counter, thus establishing benchmarks in expected value constrained SCO. Numerical experiments show the efficiency of our algorithm over real-world applications.
\end{abstract}


	\section{Introduction}\label{sec:intro}
	Over the past few years, stochastic compositional optimization (SCO) problems have attracted increasing attention. These problems are characterized by their compositional structure, which involves the optimization of functions that are composed of multiple layers of expectations:
	\begin{equation*}
	\min_{x \in \cX} \ F(x) = \EE_{\xi_1}\bigl[ f_1 (  \EE_{\xi_2}[ f_2 (x , \xi_2) ], \xi_1 )\bigr],
	\end{equation*}
	where $\cX \subset \RR^{d_x}$ is  convex, and $f_2(x, \xi_2)$ and $f_1(y,\xi_1)$ are two  functions depending on the random vectors $\xi_2$ and $\xi_1$, respectively. The SCO problems have widespread applications across various fields, including reinforcement learning \citep{wang2017stochastic}, risk management \citep{ahmed2007coherent,liu2022solving,ruszczynski2006optimization}, engineering \citep{ge2015escaping,ghadimi2020single}, and machine learning \citep{chen2021solving}. 
	Below, we provide a few examples to illustrate these applications.
	
	For example, in reinforcement learning, the policy evaluation problem aims to determine the value function using empirically observed samples. The corresponding Bellman minimization problem can be recast as a two-level SCO problem \citep{wang2017accelerating}. 
	Another instance is the low-rank matrix estimation problem, which is widely used in engineering applications such as image analysis and recommendation systems. When data arrives sequentially, the online low-rank estimation problem can also be reformulated as a two-level SCO problem \citep{ghadimi2020single}. Additionally, SCO has been applied to Model-Agnostic Meta-Learning (MAML) \citep{finn2017model,Finn2019oninemetalearning}, a powerful approach for training new models from prior related tasks. It was demonstrated by \cite{chen2021solving} that MAML falls within the two-level SCO category.
	
	Given the success of SCO in addressing real-world applications, various algorithms have been developed to solve these problems by iteratively updating the solution using projected stochastic gradient descent (SGD) with an approximated gradient \citep{wang2016stochastic,wang2017accelerating}. Despite the rapid advancements in both the theoretical and computational aspects of SCO, most existing studies on constrained SCO assume that the feasible region $\cX$ is simple, allowing for efficient projections onto $\cX$. Unfortunately, this assumption is not valid in many applications where the constraints are complex, such as nonlinear or high-order polynomials, making projection computationally expensive. The challenge is further exacerbated when the constraints involve expected value functions---whether single-level or compositional---whose explicit forms are not available. These constraints are particularly crucial in practice when incorporating risk-type constraints, such as conditional value-at-risk (CVaR) or high-moment constraints. However, due to the complexity of these constraints, existing SCO approaches are not directly applicable.
	
	In this paper, we aim to solve SCO problems with single-level and two-level compositional expected-value constraints. Specifically, we focus on the following expected value constrained SCO problem that
	\begin{equation}\label{prob:sco_expectation}
	\begin{split}
	\min_{x \in \cX} \  & F(x)   = \EE_{\xi_1}\bigl[ f_1 (  \EE_{\xi_2}[ f_2 (x , \xi_2) ], \xi_1 )\bigr] ,
	\\
	\mbox{ s.t.  } &  G^{(j)}(x) \leq 0, \text{ for } j = 1,2,\cdots, m.
	\end{split}
	\end{equation}
	Here we assume that  set $\cX$ is convex, and that projection onto it is straightforward. Additionally, we assume that $F(x)$ is convex, and that each $G^{(j)}(x)$ is a convex single-level or two-level compositional function expressed as an expectation, whose explicit form is unknown and must be estimated using observed samples. 
	In Section~\ref{sec:single}, we consider the  case where $G(x)$ is a single-level expected value function and write $G(x) = g(x) =\EE_{\zeta}[ g(x, \zeta)]$.  In Section~\ref{sec:composition}, we consider the case where $G(x)$ is a two-level compositional function and write $G(x) = \EE_{\zeta_1} [g_1( \EE_{\zeta_2}[g_2(x, \zeta_2)], \zeta_1)]$.
	
	Our model serves as a powerful tool for incorporating data-driven constraints that arise from large-scale datasets or online scenarios \citep{yang2024data}. 
	Consider a scenario where a decision-maker has access to a large-scale dataset consisting of $n$ data points $\{ w_i \}_{i=1}^n$ and aims to satisfy the data-driven constraint $\frac{1}{n} \sum_{i=1}^n g(x,w_i) \leq 0$ in a decision-making process. 
	When $n$ is large, evaluating the empirical average $\frac{1}{n} \sum_{i=1}^n g(x,w_i)$ and performing projections becomes computationally intractable, even if each function $g(x,w_i)$ is simple. 
	Alternatively, we can reformulate this as a single-level expected value constraint $\EE_{w} [g (x,w)] \leq 0 $ where $w$ is a random sample drawn uniformly from the dataset, which is captured by our model \eqref{prob:sco_expectation}. 
	Similarly, let $\{ w_{1,i}\}_{i=1}^{n_1}$ and $\{ w_{2,j}\}_{j=1}^{n_2}$ be two datasets, and let $g_1(z,w_{1,i})$ and $g_2(x,w_{2,j})$ be two functions associated with $w_{1,i}$ and $w_{2,j}$, respectively. The empirical two-level constraint $\frac{1}{n_1}\sum_{i=1}^{n_1} g_1(   \frac{1}{n_2}\sum_{j=1}^{n_2} g_2(x,w_{2,j})   , w_{1,i}) \leq 0 $ can be reformulated as a two-level compositional expected value constraint $ \EE_{w_1}[ g_1( \EE_{w_2} [ g_2(x,w_2)] , w_1)] \leq 0$, where $w_1$ and $w_2$ are random samples from $\{ w_{1,i}\}_{i=1}^{n_1}$ and $\{ w_{2,j}\}_{j=1}^{n_2}$, respectively. 
	Furthermore, in online scenarios where data arrive sequentially and each data point $w$ induces a loss $g(x,w)$ under decision variable $x$, we can formulate an expected value constraint $\EE_w[g(x,w)] \leq 0$.  
	
	\subsection{Motivating Applications} \label{sec:applications}
	In optimization, most existing approaches treat constraints as penalties on the objective function, with weights assigned to each penalty based on its relative importance. However, in certain cases, hard constraints are necessary to ensure they are strictly satisfied. For example, in finance, regulatory requirements often mandate that the risks of a portfolio do not exceed specific thresholds. In such cases, the commonly used penalty approach may not be suitable, as it could result in solutions that violate these critical hard constraints. In this paper, we propose a model that directly incorporates hard constraints within the big-data regime and introduce an algorithm to solve these problems. Our algorithm can handle multiple hard constraints and guarantees that all are satisfied. Furthermore, our modeling approach is not in conflict with the penalty approach and can be applied when both soft constraints (which can be treated as penalties) and hard constraints (which must be directly modeled) are present. 
	While the penalty approach remains powerful and widely adopted, we believe that our modeling approach offers distinct advantages in addressing the aforementioned scenarios.
	
	We present four motivating applications that requires hard compositional expected-value constraints.
	
	\medskip
	\noindent  (a) \emph{Conditional value-at-risk.} We first consider the conditional value-at-risk (CVaR) problem \citep{rockafellar2000optimization}. Let $x\in \RR^{d_x}$ be the decision variable and $\ell(x,\zeta)$ be its associated random loss, with $\zeta$ being random variable. 
	We denote by $\text{CVaR}_{\alpha}(x) $ the $\alpha$-quantile of CVaR. 
	To control the risk of loss such that the $\alpha$-quantile conditional mean of $\ell(x,\zeta)$ is no more than $\gamma$, we can impose a constraint that
	\begin{equation*}\label{eq:cvar}
	\begin{split}
	\Big \{ x: \text{CVaR}_{\alpha}(x) \leq \gamma  \Big \} & = \Big \{x : \inf_{u} \big \{ u  + \frac{1}{1-\alpha} \EE_\zeta \big [ ( \ell(x,\zeta)  - u)_+ \big ]  \big \}  \leq \gamma \Big \}, 
	\end{split}
	\end{equation*}
	where $u$ is an auxiliary variable. When the distribution of $\zeta$ the above constraint cannot be expressed explicitly and must instead be formulated as an expectation. Additionally, when the loss function $\ell(x,\xi)$ is convex in $x$, the set  $\big \{ x: \text{CVaR}_{\alpha}(x) \leq \gamma  \big \}$  is also convex in $x$. 
	
	\medskip
	\noindent  (b) \emph{Risk-averse mean-deviation constraint.} Risk-averse optimization has attracted significant attention across various communities due to its wide-ranging applications in finance and other industries \citep{cole2017does,bruno2016risk,ruszczynski2006optimization,so2009stochastic}. Traditionally, studies in risk-averse optimization have focused on controlling the risk associated with a single utility function.
	\cite{haneveld2006integrated} was one of the first that considered multiple risk levels, and proposed an optimization 
	model that incorporates these risk concerns as constraints. They derived the closed-form of these constraints so that the model can be solved by classical method. However, such an approach is computationally intractable in the big-data regime when the risk-constraints is evaluated on a large dataset.

	We consider the big-data regime and propose a model that incorporates multiple risk concerns into the decision-making process using a constrained optimization approach. Specifically, we consider the case where the closed-form of risk constraints are unknown and they have to be written in the form of expectation. Our approach is applicable to instances where constraints can be empirically evaluated, and we don't require the knowledge of constraints' closed form. Letting $x \in \RR^{d_x}$ represent the decision variable and $\ell_j(x,\xi)$ for $j=1,\cdots, K$ be $K$ heterogeneous utility functions, we consider the following set of first-order risk-averse mean-deviation constraints 
	\begin{equation*}
	\Big \{ x :   \EE_\xi [ \ell_j(x,\xi)] -   \EE_\xi \big [ ( \ell_j(x,\xi)  - \EE_\xi [ \ell_j(x,\xi)]  )_+ \big ]    \geq \gamma_j, \text{ for } j = 1,\cdots, K \Big  \}.
	\end{equation*}
	The above set can be reformulated as a region defined by two-level compositional expected-value constraints:
	\begin{equation*}
	\begin{split}
	& \Big  \{ x : \EE_\xi [ g_1^{(j)} (\EE_\xi \bigl[g_2^{(j)}(x,\xi)], \xi)\bigr]  \leq  - \gamma_j , \text{ for } j = 1,2,\cdots,K \Big  \} , 
	\\
	& \ \text{ where } 
	g_1^{(j)} \big ((z,x),\xi \big )  = z +  ( \ell_j(x,\xi) + z)_+ \text{ and } g_2^{(j)}(x, \xi ) = (-\ell_j(x,\xi), x). 
	\end{split}
	\end{equation*}

	\medskip
	\noindent  (c) \emph{Portfolio optimization with high-moment constraints.} Since the early 1960s, portfolio managers have been incorporating high-moment considerations into portfolio selection. Let $w \in \RR^{d_x}$ represent the random returns of $d_x$ assets, and let $x \in \RR^{d_x}$ denote the corresponding portfolio weights. For any even number $p$, the $p$-th order moment of the portfolio's return is given by
	\begin{equation*}
	M_p(x) = \EE_w \bigl[ | w^\top x - \EE_w[ w^\top x ] |^p \bigr]. 
	\end{equation*}
	High-moment optimization has aroused considerable interest due to its ability to efficiently characterize the distributional properties of portfolios \citep{harvey2010portfolio}. For example, the standardized 4th order moment of a portfolio, known as Kurtosis, is commonly used to capture the tail distribution of portfolios. Portfolios with high Kurtosis often exhibit long tails.
	
	We consider the following portfolio optimization problem with $p$-th order moment constraints 
	\begin{equation*}
	\begin{split}
	\max_x \EE_w [ w^\top x], \  \ 
	\mbox{s.t. }  M_p(x) \leq c_p, \ \text{ for }p = 2,4,6,8,\cdots. 
	\end{split}
	\end{equation*}
	Each of the $p$-th order moment constraint can be recast as a two-level compositional expected value constraint 
	\begin{equation*}
	\begin{split}
	&  \EE_w \bigl[  g_1 ( \EE[  g_2(x, w) ] , w ) \bigr] \leq c_p,
	\\
	\text{ where }
	& g_1 \big ( (z, x ) , w \big )  =   | w^\top x - z |^p  \ \text{ and }  g_2(x,w) =  (w^\top x , x). 
	\end{split}
	\end{equation*}

	\medskip
	\noindent (d) \emph{Optimization with group fairness constraints.} Recently, fairness learning has attracted significant attention. 
	In practice, improper use of sensitive attributes, such as gender and race, can lead to discrimination against certain groups, resulting in unfair outcomes. 
	A common approach to address this issue is to introduce fairness constraints. 
	For example, let $w$ be a random feature vector and $z \in {0,1}$ be its associated sensitive attribute. 
	In linear logistic regression, we may enforce the following constraint:
	$$
	\text{Prob}( w^\top x \geq 0 | z = 0) = \text{Prob}(w^\top x \geq 0  | z  = 1).
	$$
	This is known as \textit{statistical parity}, which means that the classification decision is independent of the sensitive attribute $z$.
	\cite{zafar2019fairness} reformulates the statistical parity constraint as 
	$$
	\EE [ ( z  - \EE[z])w^\top x] = 0, 
	$$
	where $\EE[\cdot] $ represents the expectation taken over the distribution of $w$ and $z$. 
	In the online setting with streaming data, we typically have no access to the expected value $\EE[z]$ of the sensitive variable $z$. 
	To address this challenge, we can introduce a tolerance parameter $c \geq 0$ and impose the following compositional expected-value constraint for fair classification:
	\begin{equation*}
		\begin{split}
	  \ \  \EE[ ( z  - \EE[ z ] )w^\top x] \leq c \text{ and } \EE[ ( z  - \EE[ z] )w^\top x]  \geq -c. 
		\end{split}
	\end{equation*}
	
	\subsection{Related Work}

	\paragraph{Stochastic compositional optimization.} SCO was first studied by \cite{Ermoliev88}, who analyzed the asymptotic convergence behavior of a two-timescale algorithm employing two sequences of step sizes that preserve different time scales, and established an almost sure convergence guarantee. \cite{wang2016stochastic} introduced the SCGD algorithm to address two-level SCO, providing the first non-asymptotic convergence result. These convergence rate guarantees were further improved through acceleration techniques \citep{wang2017accelerating}. Subsequently, \cite{ghadimi2020single} proposed a single-timescale algorithm that achieves an $\cO(\frac{1}{\sqrt{N}})$ convergence rate for nonconvex two-level SCO, with various methods being applied to achieve faster convergence rates \citep{chen2021solving,zhang2021multilevel}. The extension to multi-level SCO was first explored by \cite{yang2019multilevel}, and $\cO(\frac{1}{\sqrt{N}})$ convergence rates for multi-level nonconvex problems were established by \cite{balasubramanian2022stochastic} and \cite{ruszczynski2021stochastic}. \cite{liu2022solving} investigated the compositional difference of convex (DC) regime and extensively studied applications in risk management. The convex regime was examined by \cite{zhang2024optimal}, who established $\cO(\frac{1}{\sqrt{N}})$ convergence rates for both two-level and multi-level SCO using Fenchel conjugates. In all of the above works, the uncertainty only comes from the objective but not the constraints.
	
	\paragraph{Expected value constrained optimization.}
	Studies of constrained optimization can be broadly divided into two streams. 
	On one hand, constraints are often treated as penalty terms within the targeted optimization problem to avoid costly projections. 
	Various penalization approaches have been extensively studied \citep{bertsekas1997nonlinear}. 
	However, these approaches often require solving subproblems, which can be computationally expensive and thus limit their applicability. 
	On the other hand, significant efforts have been made to handle constrained optimization problems directly. 
	For deterministic constraints, \cite{nemirovski2009robust} proposed a class of gradient-type methods for solving problems with functional constraints. 
	The level-set method, also known as the bundle method, has been explored by \cite{lemarechal1995new} and \cite{lin2020data}. 
	In scenarios where first-order information is unavailable, \cite{lan2021conditional} developed a class of Frank-Wolfe-based algorithms for solving deterministic optimization problems with complex constraints.
	Recently, there has been growing interest in addressing problems with stochastic constraints. \cite{lan2020algorithms} developed the cooperative stochastic approximation (CSA) algorithm, which updates the solution based on the estimated feasibility of stochastic constraints. Following this, \cite{boob2023stochastic} focused on the smooth scenario and proposed the constraint extrapolation (ConEx) method, a primal-dual algorithm that iteratively updates the primal and dual solutions using accelerated SGD. ConEx achieves a convergence rate of $\cO(\frac{1}{\sqrt{N}})$ and sets the benchmark for expected value constrained optimization. \cite{yang2024data} studied expected value constrained problems with minimax objectives, while \cite{yu2017online} investigated the online scenario where an arbitrary objective function is introduced at each online round. For more detailed reviews on constrained optimization, see \cite{boob2023stochastic}.
	It is important to note that all of the aforementioned works focus on single-level stochastic objectives and constraints, and do not address the compositional setting. 
	
	\subsection{Contributions}
	Our work seeks to incorporate various complex expected-value constraints into the current SCO framework to better address real-world applications. The contributions of this paper are four-fold:
	\begin{itemize}
	\item We introduce a new model that integrates single-level and two-level compositional expected-value constraints into the SCO framework. This model is designed to handle complex expectation constraints that arise in real-world applications involving online or large-scale datasets. We illustrate the applicability of our model with three examples in risk management.
	
	\item  We address SCO with single-level expected-value constraints by reformulating it as a saddle point optimization problem and proposing a primal-dual algorithm that iteratively updates the primal and dual variables. We demonstrate that the solution path converges to the optimal solution at the optimal rate of $\cO(\frac{1}{\sqrt{N}})$.  Numerical experiments validate the theoretical convergence rate and showcase the practical effectiveness of our approach.
	
	\item Extending our work to SCO with two-level compositional expected-value constraints, we propose a modified primal-dual algorithm to address these challenges. We prove that the solution path achieves the same $\cO(\frac{1}{\sqrt{N}})$ convergence rate to the optimal solution. 
	
	\item We conducted numerical experiments on a stochastic portfolio optimization problem with CVaR constraints. The results suggest that our algorithm is both efficient and stable across problem instances of varying dimensions.
	\end{itemize}

	\paragraph{Notation.} For any vectors $a,b \in  \RR^n$, we denote by $\inner{ a }{ b } = \sum_{i=1}^n a_i b_i$ their inner product and denote by $\max \{ a,b\} = (\max\{a_1,b_1\},\cdots, \max \{a_n,b_n\})^\top \in \RR^{n}$ their component-wise maximum. We denote by $\| a\| = \| a\|_2$ the Euclidean norm for any vector $a \in \RR^n$ and denote by $\| A\| = \| A\|_2$ the Euclidean norm for any matrix $A \in \RR^{n \times n'}$. For any non-smooth stochastic function $f(x,\omega)$, with some slight abuse of notations, we denote by $f(x) = \EE[f(x,\omega)]$, denote by $\partial f(x)$ the collection of its deterministic subgradients, and denote by $\partial f(x,\omega)$ its stochastic subgradient  such that $\EE_\omega[\partial f(x,\omega)] \in \partial f(x)$.

	\section{Stochastic Compositional Optimization with Single-level Expected-Value Constraints}\label{sec:single}
	In this section, we consider the single-level expected value constrained stochastic compositional optimization (EC-SCO) problem
	\begin{equation}\label{prob:1}
	\begin{split}
	\min_{x\in \cX} &  \,\,\,\, F(x) = f_1 \circ f_2(x) = \EE_{\xi_1}\bigl[ f_1 (  \EE_{\xi_2}[ f_2 (x , \xi_2) ], \xi_1 )\bigr] ,\\
	\mbox{s.t. }  & \,\,\,\,G^{(j)}(x) = g^{(j)}(x)  = \EE_{\zeta_j} [g^{j}(x,\zeta_j)] \leq 0, \mbox{ for }j=1,2,\cdots,m, 
	\end{split}
	\end{equation}
	where $F$ is \emph{convex}, $\cX \subset \RR^{d_x}$ is a convex  set. Here, $\xi_1$, $\xi_2$, and $\{ \zeta_j \}_{j=1}^m$ are random vectors. 
	We assume $\xi_1$, $\xi_2$, and $\zeta = [\zeta_1,\cdots, \zeta_m]$ are independent of each other, and we write $f_1(y,\xi_1), f_2(x,\xi_2)$, $\{ g^{j}(x,\zeta_j)\}_{j=1}^m$ as the sampled function values.
	With a bit abuse of notation, we write $f_1(y) = \EE_{\xi_1}[f_1(y,\xi_1)] : \RR^{d_y} \mapsto  \RR$ as a convex and differentiable function, write $f_2(x)  =  \EE_{\xi_2}[f_2(x,\xi_2)]  : \RR^{d_x} \mapsto   \RR^{d_y} $  as a convex function, and  write $g^{(j)}(x) = \EE_{\zeta_j}[g^{(j)}(x,\zeta_j)] : \RR^{d_x} \mapsto  \RR $ for $j=1,2,\cdots, m$ as convex functions. To be specific, we write $f_2(x) = \bigl( f_2^{(1)}(x), \cdots , f_2^{(d_y)}(x) \bigr) ^\top $ and assume each  $f_{2}^{(j)}$ is a convex function. Note that the functions $f_1(y), f_2(x)$, and $g(x)$ are taken in the forms of expectation, whose explicit forms are unknown, and can only be estimated from observed samples. For notational convenience, the functions $f_1, f_2, g$ always refer to their sample-value form when they have two arguments, and always refer to their expectation form when they have one argument.
	
	To facilitate our discussion, we briefly review the Fenchel conjugate here. In particular, for any convex function $h(z): \RR^{d_z} \to \RR$,  its Fenchel conjugate is $h^*(z) =  \max_{u \in \RR^{d_z}} \langle  u, z \rangle - h(u)$. With a slight abuse of notation, for any function $\tilde h(z):\RR^{d_z} \mapsto  \RR^{k}$, where each component, $\tilde{h}_j(z)$ for $j=1,...,k$, of the function is convex, we let $\tilde{h}^*_j(z) = \max_{w\in\RR^{d_z}} \langle w,z \rangle - \tilde{h}_{j}(w)$, and let $\tilde{h}^*(z) = \bigl(\tilde{h}^*_1(z),...,\tilde{h}^*_k(z)\bigr)^\top$. That is, we let the Fenchel conjugate of $\tilde{h}(z)$ be $\tilde{h}^*(z): \RR^{d_z} \to \RR^{k}$, where the $j$-th component of $\tilde{h}^*(z)$, $\tilde{h}^*_j(z)$, is the Fenchel conjugate of $\tilde{h}_j(z)$, the $j$-th component of $\tilde{h}(z)$.
	 For any convex smooth function $f$, we denote by $D_f(x,y) = f(y) - f(x) - \inner{ \nabla f(x) }{  y - x } $ its associated Bregman's distance for ease of presentation. 
	Given $x \in \RR^{n_2}$, we denote by $g(x) = \bigl( g^{(1)}(x), \cdots , g^{(m)}(x) \bigr)^\top  \in \RR^m$ the vector of constraint values and denote by   $g(x, \zeta) = \bigl( g^{(1)}(x,\zeta_1), \cdots , g^{(m)}(x,\zeta_m) \bigr)^\top  \in \RR^m$ the vector of the sample constraint values. 
	
	We begin by specifying the stochastic sampling environment and assume access to a blackbox sampling oracle that provides stochastic zeroth- and first-order information for the objective functions $f_1, f_2$ and constraint function $g$ upon each query. Specifically, we assume the availability of the following \emph{Sampling Oracle $(\cS\cO)$} such that  
	\begin{itemize}
	\item Given $x \in \cX$, the $\cS\cO$ returns a sample value $ f_2(x , \xi_2)\in \RR^{d_y}$ and a sample subgradient  $$\partial { f_2}(x , \xi_2) \in \RR^{d_x \times d_y  }.$$
	\item Given $x \in \cX $, the $\cS\cO$ returns a sample value $ g(x , \zeta) \in \RR^m $ and a sample subgradient $$ \partial g(x , \zeta) \in \RR^{d_x \times m}.$$
	\item Given $y \in \RR^{n_1}$, the $\cS\cO$ returns a sample value $ f_1(y , \xi_1) \in \RR $ and a sample gradient $$\nabla { f_1}(y , \xi_1) \in \RR^{ d_y }.$$
	\end{itemize}
	Note that in the above $\cS\cO$, we do not impose any smoothness assumption for the inner level function~$f_2$ and the  constraint function $g$, and we let $\partial f_2(x,\xi_2)$ and $\partial  g(x , \zeta)$ be their sample subgradients, respectively.  
	We then introduce the Lagrangian  dual function of problem \eqref{prob:1} and briefly illustrate the challenges of solving the problem. 
	Specifically,   the Lagrangian dual function of problem \eqref{prob:1} is 
	\begin{equation}\label{def:lagrangian}
	\cL(x,\lambda)  = F(x) + \sum_{j=1}^{m} \lambda_j g^{(j)}(x ), \text{ where }\lambda \ge 0.
	\end{equation}
	By the convexity of $F$ and $g_j$'s, the above Lagrangian $\cL(x,y)$ is convex in $x$ and concave in $\lambda$. 
	We assume that there  exists a \emph{saddle point} $(x^*,\lambda^*)$ such that
	\begin{equation*}
	(x^*,\lambda^* )  = \argmin_{x \in \cX} \max_{\lambda \in \RR_+^m} \cL(x, \lambda ),
	\end{equation*}
	and we have that for any $(x,\lambda) \in \cX \times  \RR_+^m$, it holds that
	\begin{equation}\label{eq:lag_sad_exist}
	\cL(x^*,\lambda) \leq \cL(x^*,\lambda^*) \leq \cL(x,\lambda^*).
	\end{equation}
	Note that the existence of such a saddle point is guaranteed under the mild Slater condition \citep{slater2014lagrange}. 
	With a slight abuse of notation, the first-order partial subgradients of $\cL$ satisfy
	\begin{equation*}
	\partial_x \cL(x,\lambda )  = \partial F(x) + \sum_{j=1}^m  \lambda_j  \partial g^{(j)}(x)  \text{ and } \nabla_\lambda \cL(x,\lambda)  = g(x), \ \ \forall x \in \cX, \lambda \in \RR_+^m.
	\end{equation*}

	\paragraph{Challenges.} The key difference between expected value constrained SCO and classical minimax optimization \citep{liu2022adaptive,tan2018stochastic,yan2019stochastic}  is the unavailablity of Lipschitz continuity within the Lagrangian dual function $\cL(x , \lambda)$~\eqref{def:lagrangian}. When the classical minimax algorithms  are employed to compute the saddle point $(x^*,\lambda^*) \in \argmin_{x\in \cX} \max_{\lambda \in \RR_+^m} \cL(x , \lambda)$, existing works often assume $\cL(x , \lambda)$ is Lipschitz continuous in $(x,\lambda)$ for some constants $C_x, C_\lambda >0$ that
	\begin{equation*}
	\| \cL(x_1, \lambda_1) -  \cL(x_2, \lambda_2) \|\leq  C_x \| x_1 - x_2\| + C_\lambda \| \lambda_1 - \lambda_2\|. 
	\end{equation*}
	However, such a Lipschitz continuous property does not naturally hold for expected value constrained optimization because the dual variable is defined on the region $\lambda \in \RR_+^m$ of an infinite diameter. As a result, when $F(x)$ and $g(x)$ are  Lipschitz continuous such that $\| F(x_1) - F(x_2)\| \leq C_f \| x_1 - x_2\|$ and $\| g(x_1) - g(x_2)\| \leq C_g \| x_1 - x_2\|$, we only have 
	\begin{equation}\label{eq:Lip_continuity}
	\| \cL(x_1, \lambda) -  \cL(x_2, \lambda) \|\leq  ( C_f  + C_g \| \lambda \| )\| x_1 - x_2\|, \ \ \text{ for all } \lambda \in \RR_+^m,
	\end{equation}
	which loses the Lipschitz continuity when $\| \lambda \| \to \infty$.  Consequently, existing results on stochastic minimax optimization do not apply to the expected value constrained setting. It requires the development of new algorithms and theories to ensure the boundedness of the generated dual solution $\{ \lambda_t\}$ so that $\cL(x_t,\lambda_t)$ preserves Lipschitz continuity over the solution path $\{ (x_t,\lambda_t )\}$.
	
	In this paper, we aim to develop a class of efficient primal-dual algorithms to solve the expected value constrained optimization, which is challenging from several aspects.
	
	\begin{itemize}
	\item \emph{Primal update:} 
	At iteration $t$, we  first update the primal variable $x_t$ by obtaining an estimator of $\partial_x \cL(x_t,\lambda_t)$.  
	The key challenge to solving the EC-SCO problem lies in the unavailability of an unbiased estimator of the subgradient of $F$. Specifically, by employing the chain rule, an unbiased estimator of  $\partial F(x) $ is  
	\begin{equation}\label{eq:true_gradient_F}
	\partial  f_2 (x, \xi_2)  \nabla f_1  \big (  \EE_{\xi_2} [ f_2(x, \xi_2)], \xi_1 \big )  .
	 \end{equation}
	Unfortunately, the absence of knowledge  about $\EE_{\xi_2} [f_2(x_t, \xi_2)]$  induces a bias when we use the plug-in sample estimator $\nabla f_1( f_2 (x_t, \xi_2 ) , \xi_1)$  \citep{wang2017stochastic}. Without an unbiased estimator of $\partial F(x)$, the primal update is challenging.
	\item \emph{Dual update:}
	Letting $\{ \lambda_t\}$ be a sequence of generated dual variables, based on \eqref{eq:Lip_continuity}, we note that $\cL(x,\lambda_t)$ is  $( C_f  + C_g \| \lambda_t \| )$-Lipschitz continuous in $x$ if $\lambda_t$ is bounded. Unfortunately, the boundedness of dual variables $\{ \lambda_t\}$ only holds under certain \emph{restrictive} settings, such as vanilla expected value constrained stochastic optimization \citep{boob2023stochastic,madavan2021stochastic}. In comparison with the classical single-level stochastic optimization problems, the compositional objective $F(x)$ is much more complicated, and existing results cannot be applied to prove the boundedness of the dual variables. 
	\end{itemize}
	
	To address the above-mentioned obstacles, in the next subsection, we develop an efficient algorithm to tackle problem \eqref{prob:1}. Before presenting this algorithm, we impose some  assumptions on the   objective functions $f_1$ and $f_2$.
	\begin{assumption}\label{assumption:01}
	Let $C_{f_1}, C_{f_2}, \sigma_{f_1}, \sigma_{f_2}, L_{f_1}$, and $D_X$ be positive constants. We assume that:
	\begin{enumerate}[(a)]
	\item The set $\cX$ is convex and bounded such that $\sup_{x\in \cX} \| x - x^*\|\leq D_X$ where $x^*$ satisfies \eqref{eq:lag_sad_exist}.
	\item $f_1, f_2$ are convex, and  have  equivalent representations by their Fenchel conjugates $f_1^*(\pi_1)$ and $f_2^*(\pi_2)$ that for all $x\in \cX, y \in \RR^{d_y}$,
	$$
	f_1(y) = \max_{ \pi_1 \in \Pi_1} \langle  \pi_1 ,  y  \rangle - f_1^*(\pi_1),\text{ and  }f_2(x) = \max_{ \pi_2 \in \Pi_2} \langle  \pi_2, x  \rangle - f_2^*(\pi_2),
	$$
	where $\Pi_1 \subset \RR^{ d_y }$ and $\Pi_2 \subset \RR^{ d_x \times d_y}$ are the domains of $f_1^*$ and $f_2^*$, respectively. 
	\item $f_1$ is an $L_{f_1}$-smooth function with a Lipschitz continuous gradient
	$$\| \nabla f_1(y_1) - \nabla f_1(y_2) \| \leq L_{f_1} \|y_1 - y_2\|, \ \forall y_1, y_2 \in \RR^{d_y}. $$
	\item The sample value $f_2(x, \xi_2)$ returned by the $\cS\cO$ is unbiased and has a bounded second moment
	$$
	\EE_{\xi_2}  [ f_2(x, \xi_2)] = f_2(x),\mbox{ and } \EE_{\xi_2}   [ \| f_2(x, \xi_2) - f_2(x)\|^2] = \sigma_{f_2}^2,  \ \forall x  \in \cX. 
	$$
	\item The sample value $f_1(y,\xi_1)$ returned by the $\cS\cO$ is unbiased and has a bounded second moment 
	$$
	\EE_{\xi_1}  [ f_1(y, \xi_1)] = f_1(y),\mbox{ and } \EE_{\xi_1} [ ( f_1(y, \xi_1) - f_1(y))^2] = \sigma_{f_1}^2  , \ \forall y \in \RR^{d_y}. 
	$$
	\item The sample subgradient $\partial f_2(x, \xi_2)$ returned by the $\cS\cO$ is unbiased and has a bounded second moment  
	$$
	\EE_{\xi_2} [ \partial  f_2(x, \xi_2)]  \in \partial   f_2(x), \mbox{ and } \EE_{\xi_2}  [ \| \partial  f_2(x, \xi_2) \|^2 ] \leq C_{f_2}^2, \forall x \in \cX.
	$$
	\item  The sample gradient $ \nabla f_1(y,\xi_1)$ returned by the $\cS\cO$ is unbiased and has a bounded second moment 
	$$
	\EE_{\xi_1} [ \nabla f_1(y, \xi_1)]  = \nabla  f_1(y), \mbox{ and } \EE_{\xi_1}  [ \| \nabla f_1(y, \xi_1) \|^2 ] \leq C_{f_1}^2, \ \forall y \in \RR^{d_y}.
	$$
	\end{enumerate}
	\end{assumption}
	As stated in (b) above, when $f_1$ and $f_2$ are convex, we have their equivalent representations through their Fenchel conjugates. 
	Also, we recall that for any convex but possibly nonsmooth function $f$, 
	 for a dual variable $\pi_y$ associated with a primal variable $y$ such that $\pi_y \in \argmax_{\pi \in \Pi_y } \langle  \pi,  y  \rangle - f^*(\pi )$, we have that $\pi_y $ is  a subgradient of {$f(y)$}.  In particular, we have 
	\begin{equation*}
	 \pi_y \in  \argmax_{ \pi \in \Pi_y } \langle  \pi ,  y  \rangle - f^*(\pi) \iff \pi_y \in \partial f(y) \iff f(y) = \langle \pi_y ,  y  \rangle - f^*(\pi_y).
	\end{equation*}
	
	Further, we note that for any $L_f$-smooth convex function $f$, the Bregman distance function generated by its Fenchel conjugate function is $\tfrac{1}{L_f}$-strongly convex (see, for example,  Theorem 5.26 in~\cite{Beck2017First}) that
	\begin{equation*}
	\cD_{f^*}(\bar \pi;  \pi) = f^*( \bar  \pi) - f^*(  \pi ) - \inner{\nabla  f^*(  \pi)}{ \bar \pi -  \pi  } \geq \tfrac{1}{2 L_f} \norm{\pi - \bar \pi}^2, \quad  \forall \pi , \bar\pi \in {\rm dom}(f^*).
	\end{equation*}
	Thus the $L_{f_1}$-smoothness assumption of $f_{1}$ implies that $\cD_{f_1^*}$ is $\frac{1}{L_{f_1}}$-strongly convex.
	
	Then, we make the following assumptions for the constraint functions $g^{(j)}$'s. 
	\begin{assumption}\label{assumption:02}
	Let $C_g$ and $\sigma_g$ be positive constants. The constraint functions $g^{(j)}$'s satisfy:
	\begin{enumerate}[(a)]
	\item For $j=1,\cdots, m$, each component $g^{(j)}$ is convex such that it admits an equivalent representation using its Fenchel conjugate that
	$$
	g^{(j)}(x) = \max_{\bar v^{(j)} \in V^{(j)}} \langle  x , \bar v^{(j)}  \rangle  - [g^{(j)}(v)]^{\ast}, \ \ \forall x \in \cX, 
	$$ 
	where $V^{(j)} \subset \RR^{n_2}$ is the domain of $[g^{(j)}]^*$. 
	\item The sample value $g(x,\zeta)$ returned by the $\cS\cO$ is unbiased and has  a bounded second moment such that 
	$$ \EE_{\zeta}  [ g(x,\zeta)] = g(x), \text{ and } \EE_{\zeta}   [ \| g(x,\zeta) -  g(x) \|^2] \leq \sigma_g^2, \ \ \forall x \in \cX.$$
	\item The sample  subgradient $\partial  g(x,\zeta)$ returned by the $\cS\cO$ is unbiased and has a bounded second moment 
	$$
	 \EE_{\zeta}   [\partial  g(x,\zeta)]  \in \partial  g(x), \mbox{ and } \EE_{\zeta}   [  \| \partial  g(x,\zeta) \|^2] \leq C_{g}^2, \ \ \forall x \in \cX .
	 $$ 
	\end{enumerate}
	\end{assumption}
	
	In this paper, we focus on two scenarios: (i) the inner function $f_2$ is affine; (ii) the inner function $f_2$ is non-affine and the outer function is monotone non-decreasing, which is specified in the next assumption. 
	\begin{assumption}\label{assumption:03}
	For a non-affine function $f_2$, $\nabla f_1(x)$ is always component-wise nonnegative, i.e., $ \nabla f_1(y) \geq 0$ for any $ y \in \RR^{d_y}$. That is, the outer function $f_1$ is monotone non-decreasing if $f_2$ is non-affine.
	\end{assumption}
	
	For ease of notation, in what follows, we  omit the subscripts $\xi_1,\xi_2$, and $\zeta$ within $\EE_{\xi_1}[\cdot ], \EE_{\xi_2}[\cdot ]$, and $\EE_{\zeta}[\cdot]$.

	\subsection{A Primal-Dual Algorithm}\label{sec:alg_single}
	We now propose our  algorithm to tackle problem \eqref{prob:1}. 
	Our algorithm operates iteratively, alternately updating the primal variable $x$ and dual variable $\lambda$.
	At the $t$-th iteration, we begin by querying the $\cS\cO$  \emph{twice} at the previous solution $x_t$ to obtain a sample of $ f_2(x_t, \xi_{2,t}^1)$ and an independent sample of subgradient $\partial f_2(x_t, \xi_{2,t}^0)$. Then, given the previous estimator $y_t$ for $\EE[f_2(x_t, \xi_2 )]$ from the previous iteration, we update our estimator. Specifically, we take a weighted average of $y_t$ and $ f_2(x_t, \xi_{2,t}^1)$ such that for some $\tau_t\in(0,1)$
	\begin{equation}\label{eq:update_y} 
	 y_{t+1} =  \frac{f_2(x_{t},\xi_{2,t}^1) + \tau_{t} y_{t}}{1+\tau_{t}}.
	\end{equation}
	Next, given $y_{t+1}$, we query the $\cS\cO$ at $y_{t+1}$ to obtain a sample gradient $\nabla f_1(y_{t+1}, \xi_{1,t}^0)$. In addition, we query the  $\cS\cO$ at $x_t$, where $x_t$ is the solution from the previous iteration, for a sample subgradient  $\partial  g(x_t,\zeta_t^0)$. Then, we update the primal solution $x_t$ by a projected stochastic gradient step such that 
	\begin{equation}\label{eq:update_x}
	\begin{split}
	x_{t+1} & = \argmin_{x\in \cX}  \Big \{ \bigl\langle  \partial  f_2(x_t ,\xi_{2,t}^0 )   \nabla f_1(y_{t+1} , \xi_{1,t}^0 )     +   \partial  g(x_{t},\zeta_t^0 ) \lambda_t ,  x \bigr\rangle       + \frac{\eta_t}{2} \| x - x_t\|^2 \Big \} 
	\\
	& = \text{Proj}_{\cX} \left  [ x_t - \bigl\{{ \partial  f_2(x_t ,\xi_{2,t}^0 )   \nabla f_1(y_{t+1} , \xi_{1,t}^0 )     +   \partial  g(x_{t},\zeta_t^0 ) \lambda_t }\bigr\}/{\eta_t} \right ]
	,
	\end{split}
	\end{equation}
	where $\eta_t>0$ is prespecified. 
	
	For the dual update part, we query an independent sample of  $g(x_{t},\zeta_t^1)$, and update the dual variable $\lambda_{t+1}$ by a projected stochastic {gradient ascent} step that 
	\begin{equation}\label{eq:update_update}
	\lambda_{t+1} \in   \argmax_{ \lambda \in \RR_+^m } \Big \{ \inner{ g(x_{t},\zeta_t^1)  }{ \lambda } - \frac{\alpha_{t}}{2} \| \lambda_t - \lambda \|^2 \Big \}
	=  \Big [ \lambda_t +  g(x_{t},\zeta_t^1) /\alpha_t  \Big ]_+,
	\end{equation}
	where $\alpha_t>0$ is prespecified.
	We summarize the proposed algorithm in Algorithm~\ref{alg:01}. 
	
	\begin{algorithm}[t]  
	 \caption{Expected Value Constrained Stochastic Compositional Gradient Descent (EC-SCGD)} \label{alg:01}
	\begin{algorithmic}
	 \STATE{\bfseries Input : } Step-sizes $\{ \alpha_t \}$, $\{ \eta_t\}$, $\{ \tau_t\}$, initial points $x_0 \in \cX$, $y_0= 0$, $\lambda_0 \in \RR_+^m $, sampling oracle $\cS\cO$\\
	 \FOR{$t = 0, 1, 2, ..., N-1$}
	\STATE Query the $\cS\cO$ at $x_t$ twice for the sample values $f_2(x_t, \xi_{2,t}^i)$ and sample subgradients $ \partial f_2(x_t,\xi_{2,t}^i )$,  for $i = 0,1$.
	 \\
	 \STATE Update $$y_{t+1} =  \frac{f_2(x_{t},\xi_{2,t}^1) + \tau_{t} y_{t}}{1+\tau_{t}}.$$ 
	 Query $\cS\cO$ once at $y_{t+1} $ to obtain a sampled gradient  $ \nabla f_1(y_{t+1} ,\xi_1^0 )$. 
	 \\
	 \STATE Query the $\cS\cO$ once at $x_t$ to obtain a sampled subgradient $\partial g(x_{t},\zeta_t^0)$. 
	 \\
	 \STATE Update the main solution $x_{t+1}$ by 
	\begin{equation*}
	x_{t+1} =  \text{Proj}_{\cX} \left  [ x_t - \frac{ \partial  f_2(x_t ,\xi_{2,t}^0 )   \nabla f_1(y_{t+1} , \xi_{1,t}^0 )     +   \partial  g(x_{t},\zeta_t^0 ) \lambda_t }{\eta_t} \right ]. 
	\end{equation*}
	\STATE  Query the $\cS\cO$ once at $x_t$ to obtain $g(x_{t},\zeta_t^1 ) $, update the dual variable $\lambda_{t}$ by 
	\begin{equation*}
	\begin{split}
	\lambda_{t+1} =  \Big [ \lambda_t +  g(x_{t},\zeta_t^1) /\alpha_t  \Big ]_+.
	\end{split}
	\end{equation*}
	\ENDFOR
	 \STATE{\bfseries Output : }  $\bar x_N = \tfrac{1}{N}\sum_{t=1}^N x_t$
	  \end{algorithmic}
	\end{algorithm}

	\subsection{Stochastic Sequential Dual Interpretation} \label{sec:ssd_single}
	We  provide a   stochastic sequential dual (SSD)  interpretation \citep{zhang2019efficient,zhang2024optimal} of Algorithm \ref{alg:01}, which is the key for our  later analysis. The SSD interpretation relies on the Fenchel conjugate reformulations for the convex objective and constraints. 
	To build up the SSD framework, we first define the following composite Lagrangian functions for $f_1, f_2$ and constraint functions $g^{(j)}$:  
	\begin{equation}\label{def:Lagrangian_single}
	\begin{split}
	\cL_{f_2} (x,\pi_2) & =   \pi_2^\top  x  - f_2^*(\pi_2), \\
	\cL_F(x,\pi_{1:2}) & = \cL_{f_1} (x,\pi_1, \pi_2) =   \pi_1^\top  \cL_{f_2} (x;\pi_2)    - f_1^*(\pi_1),\\
	\cL_{g^{(j)}} (x,v^{(j)} ) & =   x^\top   v^{(j)}     - [g^{(j)}]^\ast(v^{(j)}), \mbox{ for } j=1,\cdots, m.
	\end{split}
	\end{equation}
	In what follows, we denote  by 
	$v  = [ v^{(1)} , \cdots,  v^{(m)} ] \in \RR^{ d_x  \times m}$ the collection of $v^{(j)}$'s, denote by $\cL_{g} (x, v) = [\cL_{g^{(1)}} (x,v^{(1)} ),\cdots,\cL_{g^{(m)}} (x,v^{(m)} )   ]^\top  \in \RR^m$, and denote by $\cV = [V^{(1)},\cdots, V^{(m)}]  \subset \RR^{ d_x  \times m}$ the domain of $g^*$ for ease of presentation. 
	Let $ \lambda \in \RR_{+}^m$ be the Lagrangian dual variable associated with the constraints $g^{(j)}$'s. We have the following overall composite Lagrangian function
	\begin{equation*}
	 \cL(x,\lambda, \pi_{1:2},  v ) = \cL_F(x,\pi_{1:2}) + \sum_{j=1}^m \lambda_j  \cL_{g^{(j)}}(x,v^{(j)}),  \ \text{ where } \lambda \geq 0.
	\end{equation*}
	Letting $\Pi = \Pi_1 \times \Pi_2$, 
	problem~\eqref{prob:1} has an equivalent formulation as the following saddle point problem that 
	\begin{equation}\label{prob:Lag_01}
	 \min_{x\in \cX}\Big ( \max_{\pi_{1:2} \in \Pi}\cL_F(x,\pi_{1:2}) +  \max_{\lambda \in \RR_{+}^m} \max_{ v  \in \cV}\lambda^\top  \cL_g(x, v ) \Big ).
	\end{equation}
	The existence of a saddle point  to \eqref{prob:1} and {the first-order optimality condition of $(x^*,\lambda^*)$}  ensure that there exist $\pi_2^* \in \partial f_2(x^*)$, $\pi_1^* = \nabla f_1(f_2(x^*))$, and $v^{*}  \in \partial g(x^*)$ such that $(x^*,\lambda^*,\pi_{1:2}^*,v^*)$ constitutes a saddle point to problem~\eqref{prob:Lag_01} (see, for example, Proposition 1 of \cite{zhang2024optimal}), i.e., for any $(x,\lambda,\pi_{1:2},v) \in \cX \times \RR_+^m \times \Pi \times \cV$, 
	\begin{equation}\label{eq:saddle_point}
	\cL(x,\lambda^*,\pi_{1:2}^*, v^*) \geq  \cL(x^*,\lambda^*,\pi_{1:2}^*, v^*) \geq   \cL(x^*,\lambda,\pi_{1:2}, v ). 
	\end{equation}
	
	We then provide interpretations for both primal and dual update steps of Algorithm~\ref{alg:01} under the stochastic sequential dual framework. 
	Specifically, 
	for the primal update, given a solution $x_t$, we consider the  inner layer $f_2$ and let
	\begin{equation}\label{eq:pi_2t}
	\pi_{2, t+1} \in \argmax_{ \pi_2 \in \Pi_2} \{  \pi_2^\top  x_t   - f_2^*(  \pi_2)  \}. 
	\end{equation}
	Here $\pi_{2,t+1}$ is a subgradient of the \emph{deterministic} function $f_2$ at $x_t$ that $\pi_{2,t+1} \in \partial f_2(x_t)$. Correspondingly, we have that 
	$$\cL_{f_2}(x_t, \pi_{2,t+1}) =  \pi_{2,  t+1}^\top  x_t   -  f_2^*(\pi_{2, t+1} )  = f_2(x_t).$$
	Under the stochastic setting, we let $\pi_{2,t+1}^0  \in  \partial f_2(x_t, \xi_{2,t}^0)$ be the \emph{stochastic} dual variable  such that  
	$$
	\EE[  \pi_{2,t+1}^0  \mid x_t] \in   \partial f_2(x_t ). 
	$$
	Then, recall \eqref{eq:update_y} that  $y_{t+1} =  ({f_2(x_{t},\xi_{2,t}^1) + \tau_{t} y_{t}}) / ({1+\tau_{t}})$. Using $y_{t+1}$ as an estimator for $f_2(x_t)$, we let 
	\begin{equation*}
	\pi_{1, t+1}  \in \argmax_{\pi_1 \in \Pi_1} \bigl\{   \pi_1^\top   y_{t+1}    - f_1^*(\pi_1)  \bigr\},
	\end{equation*}
	which equals the gradient of $f_1$ at $y_{t+1}$ that 
	$
	\pi_{1, t+1}  = \nabla f_1(y_{t+1}).
	$
	We also define a stochastic dual variable $ \pi_{1, t+1}^0  = \nabla f_1(y_{t+1},\xi_{1,t}^0 )
	$ such that 
	$$
	\EE[ \pi_{1, t+1}^0 \mid y_{t+1}] = \nabla f_1(y_{t+1}).
	$$
	Meanwhile, for the constraints $g(x)$, we consider its Fenchel conjugate and let
	\begin{equation}\label{eq:dual_LG_single}
	\begin{split}
	v_{t+1} & \in  \partial  g(x_t), \ \ v_{t+1}^i   \in  \partial    g(x_t,\xi_t^i), \text{ for }i =0,1, \\
	g(x_t)  & = \cL_{g}(x_t, v_{t+1}), \text{ and } g(x_t,\zeta_t^1)  = \cL_{g}(x_t, v_{t+1}^1 ).
	\end{split}
	\end{equation}
	Using the above dual interpretation, the primal update step  \eqref{eq:update_x} is equivalent to 
	\begin{equation*}
	\begin{split}
	x_{t+1} & \in \argmin_{x\in \cX}  \Big \{ \big  \langle  \pi_{2, t+1}^0  \pi_{1, t+1}^0   +  v_{t+1}^0  \lambda_t,  x  \big \rangle   + \frac{\eta_t}{2} \| x - x_t\|^2 \Big \} 
	\\
	& = \text{Proj}_{\cX} \left [  x_t - \frac{\pi_{2, t+1}^0  \pi_{1, t+1}^0   +  v_{t+1}^0  \lambda_t}{\eta_t}\right ].
	\end{split}
	\end{equation*}
	For the dual update step \eqref{eq:update_update},  using the dual interpretation \eqref{eq:dual_LG_single}, it is equivalent to
	\begin{equation}\label{eq:lambda_ssd_update}
	\lambda_{t+1} \in  \argmax_{ \lambda \in \RR_+^m } \Big \{ \inner{  \cL_{g}\big (x_{t},v_{t+1}^1  \big ) }{  \lambda } - \frac{\alpha_{t}}{2} \| \lambda_t - \lambda \|^2 \Big \}
	=  \Big [ \lambda_t +  \frac{ \cL_{g} \big (x_{t},v_{t+1}^1  \big )}{ \alpha_t }  \Big ]_+. 
	\end{equation}
	
	The above dual interpretation helps us analyze the structures within the compositional objectives and expectation constraints. In particular, under the above SSD framework, we provide the convergence results  of Algorithm~\ref{alg:01} in the next subsection.
	
	\subsection{Convergence Analysis}
	After presenting Algorithm \ref{alg:01} and its dual interpretation, we now analyze its convergence properties. LetAfter presenting Algorithm \ref{alg:01} and its dual interpretation, we now analyze its convergence properties. Let $\{ ( x_t, \lambda_t ) \}$ be the sequence generated by Algorithm \ref{alg:01}, and let $\{  ( \pi_{1:2,t},v_t)\}$ represent its corresponding dual sequence. We denote by $\bar x_N = \tfrac{1}{N}\sum_{t=1}^N x_t$. To evaluate the performance of our algorithm, we consider the following metrics: the objective optimality gap $F(\bar x_N) - F(x^*)$ and the feasibility residual $\| g(\bar x_N)_{+}\| = \|\max\{ g(\bar x_N), 0 \} \|$. 
	It is clear that for any $x \in \cX$ is an optimal solution to problem \eqref{prob:1} if and only if $F(x) - F(x^*) = 0$ and $\| g(x)_+\| = 0$. 
	In the remainder of this section, we demonstrate that both of these metrics converge to zero and derive their respective convergence rates.

	\subsubsection{Convergence of gap function}\label{sec:gap_rate}
	Letting $z_t = ( x_t, \lambda_t,  \pi_{1:2,t}, v_t ) $ be the solution generated by Algorithm \ref{alg:01} and $z =  ( x^*, \lambda,  \pi_{1:2}, v ) \in \cX \times \RR_+^m \times \Pi \times \cV $ be a general feasible point, 
	 to facilitate our analysis, we consider the \emph{gap function}
	\begin{equation}\label{eq:gap_function}
	Q(z_t, z) =   \cL(x_t,\lambda,\pi_{1:2}, v) -  \cL(x^*, \lambda_t,  \pi_{1:2,:t} , v_t)  .
	\end{equation}
	In what follows, we provide our first result to characterize the contraction of the above gap function. 
	\begin{theorem}\label{thm:01}
	Suppose Assumptions \ref{assumption:01}, \ref{assumption:02}, and \ref{assumption:03} hold and Algorithm \ref{alg:01} generates  $\{ (x_t,\lambda_{t}, \pi_{1:2,t},v_t)  \}_{t=1}^N $ by setting $\tau_t = t/2$, $\alpha_t = \alpha $, and $\eta_t = \eta $ for $t \leq N$. Let $\lambda \in \RR_+^m$ be a nonnegative {bounded} (random)  variable \ such that $\| \lambda\| \leq  M_{\lambda}$ {uniformly} and let $z = (x^*,\lambda,\pi_{1:2},v)$ be a feasible point, then for any integer $K \leq N$,  we have
	\begin{equation}\label{eq:thm_single}
	\begin{split}
	& \sum_{t=1}^K \EE [Q(z_t, z) ] 
	 +  \frac{\alpha }{2}\EE[ \| \lambda_K - \lambda\|^2] 
	\\
	 & \leq     \frac{95 K C_{f_1}^2 C_{f_2}^2 }{2\eta} 
	 + \frac{\eta }{2} \| x_{0}- x^*\|^2 
	+     \frac{K}{\alpha } D_X^2 C_g^2 + \sqrt{K} \Big (  M_{\lambda}  \sigma_g  +   L_{f_1} \sigma_{f_2}^2  +   3  C_{f_1}\sigma_{f_2}  \Big )  \\ 
	& \quad 
	+   \frac{\alpha }{2}  \EE [ \| \lambda_{0}- \lambda  \|^2   ]
	+ \frac{K \sigma_{g}^2}{ \alpha }
	 +   \sum_{t=1}^K 5  C_g^2 \Big ( \frac{   \EE[  \|  \lambda_{t-1}\|^2]  }{2\eta   }  +  \frac{2 M_{\lambda}^2 }{ \eta  } + \frac{ \EE[  \|  \lambda  - \lambda_{t} \|^2 ] }{2\eta }  \Big ) .
	 \end{split}
	\end{equation}
	\end{theorem}
	
	The above result serves as a crucial building block when analyzing the convergence behavior of Algorithm~\ref{alg:01}. A key observation is that  the dual sequence $\{ \lambda_t\}$ appears on both sides of \eqref{eq:thm_single}, so that it is not guaranteed to be bounded. In this case, we cannot readily establish the overall convergent property of the gap function $Q(z_t,z)$. 
	In the next two subsections, without imposing any additional assumptions, we show that the dual sequence $\{\lambda_t\}$ is \emph{bounded} by using \eqref{eq:thm_single}, and further derive the convergence rates of both objective optimality gap  and feasibility residual.
	 
	\subsubsection{Boundedness of \texorpdfstring{$\{ \lambda_t\}$}{lambda}}
	To show the boundedness of $\{\lambda_t\}$, we observe that \eqref{eq:thm_single} in Theorem \ref{thm:01} preserves a recursive structure. Specifically, in iteration $K$, the dual solution $\lambda_K$ can be bounded by the sum of dual solutions $\{\lambda_t\}_{t=0}^{K-1}$ generated in previous iterations plus an extra constant. 
	Here, we restate a technical result from Lemma 2.8 of \cite{boob2023stochastic} to study the boundedness of such sequences. 
	\begin{lemma}\label{lemma:recursive_bound}
	Let $\{a_t\}$ be a nonnegative sequence and $m_1,m_2 \geq 0 $ be constants such that $a_0 \leq m_1$. Suppose the following relationship holds for all $K\geq 1$.
	$$
	a_K \leq m_1 + m_2 \sum_{t=0}^{K-1}a_t.
	$$
	Then we have $a_K \leq m_1 (1+m_2)^{K}$.
	\end{lemma}
	Using the above result, we then show the boundedness of  $\{ \lambda_t \}$ in the following proposition.
	\begin{proposition}\label{prop:bounded_lambda_01}
	Suppose Assumptions \ref{assumption:01}, \ref{assumption:02}, and \ref{assumption:03} hold, and  Algorithm \ref{alg:01} generates $\{ (x_t,\lambda_{t}, \pi_{1:2,t},v_t)  \}_{t=1}^N $ by setting $\tau_t = \frac{t}{2}$,  $\alpha_{t} = \alpha = 2 \sqrt{N}$, and $\eta_t = \eta =  \frac{15C_g^2 \sqrt{N}}{2}$ for all $t\leq N$. Then for any $N \geq 2$ and any integer $K \leq N$, we have 
	$$
	\EE[ \| \lambda_K - \lambda^*\|^2 ] \leq 2R e^{2 }, 
	$$
	where 
	\begin{equation}\label{eq:R_1}
	\begin{split}
	R & = \frac{  7C_{f_1}^2 C_{f_2}^2 }{C_g^2} +   \frac{D_X^2 C_g^2}{2}  + 2  \|  \lambda^* \|^2  +  \frac{ \sigma_{g}^2}{2}  +  
	 L_{f_1} \sigma_{f_2}^2  +   3 C_{f_1}\sigma_{f_2}  +   \| \lambda^*\|\sigma_g 
	 \\
	 & \quad + \frac{15C_g^2}{4} \| x_{0}- x^*\|^2 +  \| \lambda_{0} - \lambda^*  \|^2  .
	\end{split}
	\end{equation}
	\end{proposition}
	
	The above result implies that, with properly chosen step-sizes, the dual sequence $\| \lambda_t\|^2$ generated by Algorithm~\ref{alg:01} remains bounded in expectation, despite the increased randomness and bias introduced by the compositional setting. Although it is discussed in Section~\ref{sec:single} that $\cL(x,\lambda)$ generally loses Lipschitz continuity, the above result allows us to conclude that, in expectation, $\cL(x,\lambda)$ is $( C_f  + C_g (\sqrt{2R} e + \| \lambda^*\|) )$-Lipschitz continuous in $x$ over the trajectory $\{ (x_t,\lambda_t) \}_{t=1}^N$. This finding serves as a crucial building block for our subsequent convergence analysis.

	\subsubsection{Convergence of value and feasibility metrics}\label{sec:feasibility_rate}
	After proving the boundedness of $\{ \lambda_t\}$, we derive the convergence rates of the objective optimality gap  and the feasibility residual in an ergodic sense.
	
	\begin{theorem}\label{thm:rate}
	Suppose that Assumptions \ref{assumption:01}, \ref{assumption:02}, and \ref{assumption:03} hold,  and Algorithm \ref{alg:01} generates $\{ (x_t, \lambda_t,\pi_{1:2,t},v_t) \}_{t=1}^N$  by setting $\lambda_0 = 0$, $\tau_t = t/2$, $\alpha_{t} = \alpha = 2 \sqrt{N}$, and $\eta_t =  \eta = \frac{15 C_g^2 \sqrt{N}}{2}$ for all $0 \leq t \leq N$. Let $R$ be the constant defined in \eqref{eq:R_1}. Letting $\bar x_N = \sum_{t=1}^Nx_t/N$, for any $N \geq 2$, we have 
	\begin{equation*}
	\begin{split}
	 &\EE [ F(\bar x_N) ] - F(x^*) 
	 \leq    \frac{1}{ \sqrt{N}} \Big ( \frac{ 7 C_{f_1}^2 C_{f_2}^2  }{C_g^2}
	+ L_{f_1} \sigma_{f_2}^2  +   3C_{f_1}\sigma_{f_2}  + \frac{15 C_g^2 }{4} \| x_{0}- x^*\|^2 
	+  \frac{ 4\| \lambda^*\|^2  + 8Re^{2}  }{3}   \Big )
	,
	\end{split}
	\end{equation*}
	and 
	\begin{equation*}
	\begin{split}
	  \EE[ \| g(\bar x_N)_{+} \| ] 
		& \leq    \frac{1}{\sqrt{N} }  \Big ( \frac{ 7 C_{f_1}^2 C_{f_2}^2  }{C_g^2}
	+  L_{f_1} \sigma_{f_2}^2  +   3  C_{f_1}\sigma_{f_2}  + \frac{15C_g^2}{4} \| x_{0}- x^*\|^2 
	+     \frac{D_X^2 C_g^2}{2}  
	+ \frac{ \sigma_{g}^2}{ 2}  \Big ) \\ 
		& \qquad 
	 + \frac{ (2\| \lambda^*\|^2  + 4Re^{2 } )   +   3 ( \| \lambda^*\| + 1 )^2   +  ( \| \lambda^*\| + 1 )\sigma_g }{\sqrt{N} }  .
	\end{split}
	\end{equation*}
	\end{theorem}
	
	It is important to note that, given any $\lambda_0 \in \RR_+^m $, our algorithm still enjoys the $\cO(\frac{1}{\sqrt{N}})$ convergence rate for both the objective optimality gap and the feasibility residual. Notably, our approach is sufficiently general to accommodate multiple expected-value constraints.
	It is also worth noting that the above result matches the optimal $\cO(\frac{1}{\sqrt{N} })$ convergence rate for standard convex stochastic optimization, despite the complexity of the compositional objective and expected-value constraints involved. This establishes a new benchmark for stochastic compositional optimization under expected-value constraints. Finally, we emphasize that the primitive $C_{g}^2$ scales with the number of constraints $m$. Consequently, the hidden constant within $\cO(\frac{1}{\sqrt{N}})$ grows \emph{polynomially} in the number of constraints $m$. This highlights the efficiency of our algorithm in handling real-world applications with a large number of expected-value constraints.

	\section{Stochastic Compositional Optimization under Compositional Expected-Value Constraints}\label{sec:composition}
	In Section \ref{sec:single}, we have presented a primal-dual algorithm for expected value constrained stochastic compositional optimization that achieves the optimal rate of convergence. It is natural to ask whether our analytical framework can also handle the compositional expected value constrained applications discussed in Section~\ref{sec:applications}. We extend the framework by considering the following two-level compositional expected value constrained stochastic compositional  optimization (CoC-SCO) problem:
	\begin{equation}\label{prob:2}
	\begin{split}
	\min_{x\in \cX} \ & F(x) = f_1 \circ f_2 (x)  ,\\
	\mbox{s.t. }  \ & G(x) = g_1 \circ g_2(x)  \leq 0, \, 
	\end{split}
	\end{equation}
	where  $\cX \subset \RR^{d_x}$ is convex and closed,  $f_1( y  ) =  \EE_{\xi_1}[f_1( y,\xi_1)] :\RR^{d_y} \mapsto \RR $ is convex and differentiable,  $f_2(x )  = \EE_{\xi_2}[f_2(x,\xi_2)]   :\RR^{d_x} \mapsto \RR^{d_y} $ is convex, $F(x):\RR^{d_x} \mapsto \RR$ is convex, $G(x):\RR^{d_x} \mapsto \RR^m$ is convex, $g_1(z) =  \EE_{\zeta_1} [g_1(z, \zeta_1)] : \RR^{d_z} \mapsto  \RR^m$ is convex and differentiable, and $g_2(x ) = \EE_{\zeta_2} [g_2( x, \zeta_2)]: \RR^{ d_x } \mapsto  \RR^{d_z}$ is convex. In more details, we write $f_2(x) = (f_2^{(1)}(x),
	\cdots, f_2^{(d_y )}(x) )^\top \in \RR^{ d_y }$, $G(x) = (G^{(1)}(x),\cdots, G^{(m)}(x))^\top \in \RR^m$, $g_2(x) = (g_2^{(1)}(x),
	\cdots, g_2^{(d_z)}(x) )^\top \in \RR^{d_z }$, and $g_1(z) = (g_1^{(1)}(z),
	\cdots, g_1^{(m)}(z) )^\top \in \RR^{m}$, and assume that $f_2^{(j)}, G^{(j)}, g_2^{(j)}$, and $g_1^{(j)}$ are convex. As in the previous section, we consider that the functions $f_1,f_2,g_1$, and $g_2$ are  in the form of expectations and can only be estimated from observed samples.  Note that each $G^{(j)}$ represents one compositional constraint and problem \eqref{prob:2} allows us to handle $m$ compositional constraints.  
	
	We first  specify the sampling environment for CoC-SCO. In particular,  we assume  access to the following \emph{Compositional Sampling Oracle} $(\cS\cO_{\text{co}})$ such that
	\begin{itemize}
	\item Given $x \in \cX$, the $\cS\cO_{\text{co}}$ returns a sample value $ f_2(x , \xi_2) \in \RR^{d_y}$ and a sample subgradient 
	$$\partial  f_2(x , \xi_2) \in \RR^{d_x\times d_y} .$$
	\item Given $y \in \RR^{d_y}$, the $\cS\cO_{\text{co}}$ returns a sample value $ f_1(y , \xi_1) \in \RR$ and a sample gradient
	$$ \nabla f_1(y , \xi_1) \in \RR^{d_y}.$$
	\item Given $x \in \cX$, the $\cS\cO_{\text{co}}$ returns a sample value $ g_2(z ,  \zeta_2) \in \RR^{d_z}$ and a sample subgradient
	$$ \partial  g_2(x, \zeta_2) \in \RR^{d_x\times d_z }.$$
	\item Given $z \in \RR^{d_z}$, the $\cS\cO_{\text{co}}$ returns a sample value $ g_1(z ,  \zeta_1) \in \RR^m$ and a sample gradient 
	$$ \nabla g_1(z, \zeta_1) \in \RR^{d_z \times m}.$$
	\end{itemize}
	Here we do not impose any smoothness assumptions on the inner-level  functions $f_2$ and $g_2$. For the compositional expected value constrained problem \eqref{prob:2}, with a slight abuse of notation,   
	we write its Lagrangian dual function as 
	\begin{equation*}
	\cL(x,\lambda)  = F(x) + \langle G (x )  , \lambda \rangle, \text{ where }\lambda \in  \RR_+^m.
	\end{equation*}
	By the convexity of $F$ and $G^{(j)}$'s, the above Lagrangian $\cL(x,y)$ is convex in $x$ and concave in $\lambda$. We assume the Slater condition holds such that  there  exists a \emph{saddle point} $(x^*,\lambda^*)$ satisfying 
	\begin{equation*}
	\cL(x^*,\lambda) \leq \cL(x^*,\lambda^*) \leq \cL(x,\lambda^*),  \,  \, \forall (x,\lambda) \in \cX \times \RR_+^m.
	\end{equation*}
	Letting $\partial F(x) \in \RR^{d_x} $ and $\partial G(x) \in \RR^{d_x \times m}$ be the subgradients of $F(x)$ and $G(x)$, respectively, the first-order partial subgradients of $\cL$ satisfy
	\begin{equation*}
	\partial_x \cL(x,\lambda )  = \partial F(x) +  \partial G(x)\lambda  \text{ and } \nabla_\lambda \cL(x,\lambda)  = G(x), \ \ \forall x \in \cX, \lambda \in \RR_+^m.
	\end{equation*}
	Compared with  EC-SCO discussed in Section~\ref{sec:single}, the compositional expected value constrained problem is more challenging, because of the lack of  unbiased estimators of $\partial G(x)$ and $G(x)$. Specifically,  the unbiased estimators of $\partial G(x)$ and $ G(x)$ are 
	$$
	\partial g_2(x, \zeta_2) \nabla g_1(\EE_{\zeta_2} [g_2(x,\zeta_2)],\zeta_1) \text{ and } g_1(\EE_{\zeta_2} [g_2(x,\zeta_2)],\zeta_1), 
	$$
	respectively. Lacking the knowledge of the expected value $ \EE_{\zeta_2} [g_2(x,\zeta_2)]$, the plug-in estimators, $\partial g_2(x, \zeta_2) \nabla g_1( g_2(x,\zeta_2) ,\zeta_1)$ and $g_1( g_2(x,\zeta_2),\zeta_1)$,  induce biases in both primal and dual update steps. Such biases would corrupt both our primal and dual updates. 
	For the primal update step, in addition to the bias in estimating $\partial F(x)$ discussed earlier, the biased estimator of $\partial G(x)$ induces further biases to our estimation of $\partial_x \cL(x,\lambda)$.
	Furthermore, for the dual update step, the bias induced by estimating $G(x)$ would compound with the randomness within updating $\lambda_t$. Consequently, it is more challenging to analyze the boundedness of the dual sequence $\{ \lambda_t\}$ and the convergence behavior of the solution sequence $\{x_t\}$ for problem~\eqref{prob:2}.
	
	To facilitate our discussion, 
	in addition to Assumptions \ref{assumption:01} and  \ref{assumption:03} for the objective functions $f_1$ and $f_2$, we impose the following convexity, smoothness, and boundedness assumptions for the constraint functions $g_1$ and $g_2$. 
	\begin{assumption}\label{assumption:compositional_g}
	Let $C_{g_1}, C_{g_2}, \sigma_{g_1}, \sigma_{g_2}$ and $L_{g_1}$ be positive scalars. The constraint functions $g_1$ and  $g_2 $ satisfy:
	\begin{enumerate}[(a)]
	\item $g_1, g_2$ are convex such that they could be reformulated using the following Fenchel conjugates
	$$
	g_1(z) = \max_{v_1 \in \cV_1 }v_1^\top  z - g_1^*  (v_1),\text{ and } g_2(x) = \max_{ v_2 \in \cV_2 } v_2^\top  x - g_2^* (v_2),  \ \  \forall \  x \in \cX, \, z \in \RR^{d_z},
	$$
	where $\cV_1 \subset \RR^{d_z \times m }$ and $\cV_2 \subset \RR^{d_x \times d_z}$ are the domains of $g_1^*$ and $g_2^*$, respectively. 
	\item The sample value $g_1(z,\zeta_1)$ returned by the $\cS\cO_{\text{co}}$ is unbiased and  has a bounded second moment 
	$$
	\EE_{\zeta_1} [ g_1(z, \zeta_1)] = g_1(z) \mbox{ and }  \EE_{\zeta_1} [ ( g_1(z, \zeta_1) - g_1(z))^2] \leq \sigma_{g_1}^2,   \  \forall \ z \in \RR^{d_z}. 
	$$
	\item The sample value $g_2(x, \zeta_2)$ returned by the $\cS\cO_{\text{co}}$ is unbiased and has a bounded second moment 
	$$
	\EE_{\zeta_2} [ g_2(x, \zeta_2)] = g_2(x) \mbox{ and } \EE_{\zeta_2} [ \| g_2(x, \zeta_1) - g_2(x)\|^2] \leq \sigma_{g_2}^2, \ \forall \ x \in \cX. 
	$$
	\item The sample subgradient $\partial  g_2(x, \zeta_2)$ returned by the $\cS\cO_{\text{co}}$ is unbiased and has a bounded second moment
	$$
	\EE_{\zeta_2} [ \partial  g_2(x, \zeta_2)] \in \partial   g_2(x) \mbox{ and } \EE_{\zeta_2} [ \| \partial    g_2(x, \zeta_2) \|^2 ] \leq C_{g_2}^2, \ \forall \ x \in \cX.
	$$ 
	\item The sample gradient $ \nabla g_1(y,\zeta_1)$ returned by the $\cS\cO_{\text{co}}$ is  unbiased and has a bounded second moment
	$$
	\EE_{\zeta_1} [  \nabla g_1(z, \zeta_1)]  = \nabla  g_1(z) \mbox{ and } \EE_{\zeta_1} [ \| \nabla g_1(z, \zeta_1) \|^2 ] \leq C_{g_1}^2, \ \forall \ z \in \RR^{d_z}.
	$$
	\item $g_1$ is an $L_{g_1}$-smooth function with a Lipschitz continuous gradient such that
	$$\norm{\nabla g_1 (z_1) - \nabla g_1(z_2)} \leq L_{g_1} \norm{z_1 - z_2}, \ \ \forall z_1, z_2 \in \RR^{d_z}. $$
	
	\item For each constraint $G^{(j)}$, its outer level function  is monotone non-decreasing if its inner-level function  is  non-affine. 
	\end{enumerate}
	\end{assumption}
	Note that (g) above is essentially a generalization of Assumption~\ref{assumption:03} to handle the mixture of affine and non-affine inner-level compositional constraints.  
	Specifically, we allow some constraints $G^{(j)}$ to have affine inner-level functions, while allowing others to have non-affine inner-level and monotone non-decreasing outer-level functions. Also, we only assume smoothness for the outer-level constraint $g_1$ while allowing the inner-level constraint $g_2$ to be non-smooth. We  omit the subscripts $\xi_1,\xi_2,\zeta_1$, and $\zeta_2$ within $\EE_{\xi_1}[\cdot ], \EE_{\xi_2}[\cdot ], \EE_{\zeta_1}[\cdot]$, and $\EE_{\zeta_2}[\cdot]$ for notational convenience. 
	
	{The above assumptions are mild and our framework can handle the compositional expected value constrained risk management applications provided in Section~\ref{sec:applications} (b) and (c). Specifically, in the risk-averse mean-deviation application (b),  when each utility function $\ell_j(x, \xi)$ is concave  and monotone nondecreasing in $x$, we observe that $g_2^{(j)}(x,\xi)$ is convex in $x$, and $g_1^{(j)} \big ((z,x),\xi \big )$ is convex and monotone nondecreasing in $(z,x)$. {Thus, this application can be handled after smoothing the $(\cdot)_+$ operator within the outer-layer function.}
	For the high-moment portfolio selection application~(c),  given any random return $w$,  $g_2(x,w)$ is an affine function while $g_1\big ( (z, x ),w \big )$ is convex in $(z,x)$. Consequently, by taking expectation over the random return $w$, $\EE_w [g_2(x,w)]$ is affine and $\EE_w[ g_1\big ( (z, x ),w \big )]$ is convex in $(z,x)$.}

	\subsection{Algorithm}
	Now we propose our algorithm to solve problem~\eqref{prob:2}. In iteration $t$, we first construct an estimator $y_{t+1} $ for $\EE_{\xi_2} [ f_2(x_t,\xi_2)]$, or the limit point of $\{\EE_{\xi_2} [ f_2(x_t,\xi_2)]\}_t$ to be precise, by letting $y_{t+1} = \frac{f_2(x_{t},\xi_{2,t}^1) + \tau_{t} y_{t}}{1+\tau_{t}} $ for some $\tau_t >0$.
	We query  $\cS\cO_{\text{co}}$ for a sample gradient $\nabla f_1(y_{t+1}, \xi_{1,t}^0)$  and a sample subgradient $\partial f_2(x_t,\xi_{2,t}^0)$. Letting $z_{t}$ be an estimator of $g_2(x_{t-1})$, we query $\cS\cO_{\text{co}}$ to obtain $ g_2(x_t, \zeta_{2,t}^1)$, and update 
	$$
	z_{t+1} =  \frac{g_2(x_{t},\zeta_{2,t}^1) + \rho_{t} z_{t}}{1+\rho_{t}},
	 $$
	where $\rho_t >0$ is prespecified. Next, we query  $\cS\cO_{\text{co}}$ at $x_t$ to obtain a sample subgradient $\partial g_2(x_t, \zeta_{2,t}^0)$ and query  $\cS\cO_{\text{co}}$ at $z_{t+1}$ to obtain another sample gradient $\nabla g_1(z_{t+1}, \zeta_{1,t}^{0})$. We then employ $ \partial    f_2(x_t ,\xi_{2,t}^0 )  \nabla f_1(y_{t+1} , \xi_{1,t}^0 ) $ as an estimator of $\partial F(x_t)$, use $\partial g_2(x_t, \zeta_{2,t}^0) \nabla g_1(z_{t+1}, \zeta_{1,t}^{0})$ as an estimator of $\partial G(x_t)$, and  update the primal solution $x_{t}$ by a projected stochastic gradient step that 
	\begin{equation}\label{eq:update_x_comp}
	x_{t+1}
	 = \text{Proj}_{\cX} \left [ x_t - \frac{ \partial    f_2(x_t ,\xi_{2,t}^0 )  \nabla f_1(y_{t+1} , \xi_{1,t}^0 )  +   \partial    g_2(x_t,\zeta_{2,t}^0 )   \nabla g_1(z_{t+1}, \zeta_{1,t}^0) \lambda_t }{ \eta_t }\right ],
	\end{equation}
	where $\eta_t >0$ is prespecified. 
	Finally, we consider the dual variable $\lambda \in \RR_+^m$, query $\cS\cO_{\text{co}}$ to obtain $g_1(z_{t+1},  \zeta_{1,t}^1)$, $\nabla g_1(z_{t+1},\zeta_{1,t}^1) $, and $g_2(x_{t}, \zeta_{2, t}^2)$, and set 
	\begin{equation}\label{def:H}
	H_{t+1} = g_1(z_{t+1},  \zeta_{1,t}^1) -  \nabla g_1(z_{t+1},\zeta_{1,t}^1) z_{t+1}  +   \nabla g_1(z_{t+1},\zeta_{1,t}^1)g_2(x_{t}, \zeta_{2, t}^2)
	\end{equation}
	as an estimator of $G(x)$. 
	We update $\lambda_{t+1}$ by a projected stochastic gradient step that
	\begin{equation}\label{eq:update_update_comp}
	\lambda_{t+1} \in   \argmax_{ \lambda \in \RR_+^m } \Big \{ \lambda^\top H_{t+1}   - \frac{\alpha_{t}}{2} \| \lambda_t - \lambda \|^2 \Big \}
	=  \Big [ \lambda_t +  H_{t+1}/\alpha_t  \Big ]_+,
	\end{equation}
	for some $\alpha_t >0$. 
	We summarize the details of the above process  in Algorithm~\ref{alg:02}. We point out that given $x_t, \lambda_t,z_t$, and $y_t$,  the subsequent updates $x_{t+1}, \lambda_{t+1}$ are independent, because we employ independent stochastic samples within our update scheme. 
	
	\begin{algorithm}[t] 
	 \caption{Compositional-Constrained Stochastic Compositional Gradient Descent (CC-SCGD)}
	 \begin{algorithmic} \label{alg:02}
	 \STATE{\bfseries Input: } Step-sizes $\{ \alpha_k \}$, $\{ \eta_k \}$, $\{ \tau_k \}$, $\{ \rho_k \}$,  initial points $x_0 \in \cX$, $y_0= 0$, $z_0 = 0$, $\lambda_0 \in \RR_+^m $, sampling oracle $\cS\cO_{\text{co}}$.
	 \FOR{$t = 1, 2, ..., N$}
	 \STATE Query the $\cS\cO_{\text{co}}$ at $x_t$ for the sample values $f_2(x_t, \xi_{2,t}^1), g_2(x_t, \zeta_{2,t}^1),g_2(x_t, \zeta_{2,t}^2)$ and sample subgradients $ \partial f_2(x_t,\xi_{2,t}^0 ), \partial g_2(x_t,\zeta_{2,t}^0 )$.
	 \\
	 \STATE  Update 
	 $$y_{t+1} =  \frac{f_2(x_{t},\xi_{2,t}^1) + \tau_{t} y_{t}}{1+\tau_{t}} \text{ and } z_{t+1} =  \frac{g_2(x_{t},\zeta_{2,t}^1) + \rho_{t} z_{k}}{1+\rho_{t}}.$$
	  \STATE Query $\cS\cO_{\text{co}}$ at $y_{t+1}  $ to obtain $ \nabla f_1(y_{t+1} ,\xi_1^0 )$ and query $\cS\cO_{\text{co}}$ at $z_{t+1}$ to obtain $\nabla g_1(z_{t+1},  \zeta_1^0)$. 
	 \\
	 \STATE  Update the main solution $x_{t+1}$ by 
	 $$
	x_{t+1} = \text{Proj}_{ \cX}  \left [ x_t - \frac{  \partial  f_2(x_t,\xi_{2, t}^0)  \nabla  f_1(y_{t+1},  \xi_{1,t}^0)  +  \partial  g_2(x_t,\zeta_{2,t}^0)  \nabla g_1(z_{t+1},  \zeta_{1, t}^0)  \lambda_t  }{\eta_t } \right ].
	 $$
	 Query $\cS\cO_{\text{co}}$ at $y_{t+1}  $ to obtain $g_1(z_{t+1},  \zeta_{1,t}^1)$   and $\nabla g_1(z_{t+1},\zeta_{1,t}^1)  $. 
	Set $H_{t+1}$ by 
	 $$
	 H_{t+1}  = g_1(z_{t+1},  \zeta_{1,t}^1) -  \nabla g_1(z_{t+1},\zeta_{1,t}^1) z_{t+1}  +   \nabla g_1(z_{t+1},\zeta_{1,t}^1)g_2(x_{t}, \zeta_{2, t}^2) .
	 $$
	Update the dual variable $\lambda_{t+1}$ by 
	 $$
	 \lambda_{t+1} = \bigl[ \lambda_t +  H_{t+1}    /\alpha_t  \bigr]_+.
	 $$
	\ENDFOR
	 \STATE{\bfseries Output : } $\bar x_N = \tfrac{1}{N} \sum_{t=1}^N x_t$.
	 \end{algorithmic}
	\end{algorithm}

	\subsection{Stochastic Sequential Dual Interpretation}
	After introducing Algorithm~\ref{alg:02}, we now present its interpretation under the SSD framework to facilitate our analysis. 
	We first define the following composite Lagrangian functions for $f_1, f_2$:
	\begin{equation*}
	\begin{split}
	\cL_{f_2} (x,\pi_2)  =  \pi_2^\top  x - f_2^*(\pi_2),\ 
	\cL_F(x,\pi_{1:2}) =  \pi_1^\top \cL_{f_2} (x, \pi_2) - f_1^*(\pi_1),
	\end{split}
	\end{equation*}
	and define the composite Lagrangians for constraint functions $g_1, g_2$ that
	\begin{equation*}
	\begin{split}
	\cL_{g_2} (x,v_2) & =  v_2^\top  x -  g_2^*(v_2),\\
	\cL_{G}(x,v_{1:2}) & =   \cL_{g_1} (x,v_1, v_2) =  v_1^\top  \cL_{g_2} (x, v_2) - g_1^{*}( v_1).
	\end{split}
	\end{equation*}
	We then have the overall composite Lagrangian function 
	\begin{equation*}
	 \cL(x,\lambda, \pi_{1:2}, v_{1:2}) = \cL_F(x,\pi_{1:2}) + \lambda^\top \cL_G(x,v_{1:2}).
	\end{equation*}
	Letting $\Pi = \Pi_1 \times \Pi_2$ and $\cV = \cV_1 \times \cV_2$, 
	problem \eqref{prob:2} can be reformulated as a saddle point problem 
	\begin{equation}\label{def:saddle_2}
	\min_{x\in \cX}\Big ( \max_{\pi_{1:2} \in \Pi}\cL_F(x,\pi_{1:2}) +  \max_{\lambda \in \RR_{+}^m } \max_{v_{1:2} \in \cV}\lambda^\top  \cL_G(x,v_{1:2})\Big ). 
	\end{equation}
	There exist 
	$\pi_2^* \in \partial f_2(x^*), \pi_1^* = \nabla f_1( f_2(x^*)), v_2^* \in \partial g_2(x^*)$, and $v_1^* = \nabla g_1 (g_2(x^*))$ such that  $(x^*,\lambda^*,\pi_{1:2}^*, v_{1:2}^*)$ serves as a saddle point to \eqref{def:saddle_2}. Specifically,  for any feasible $(x,\lambda, \pi_{1:2}, v_{1:2}) \in \cX \times \RR_+^m \times \Pi \times \cV$, 
	\begin{equation}\label{def:saddle_2_ssd}
	\cL(x,\lambda^*,\pi_{1:2}^*, v_{1:2}^*) \geq \cL(x^*,\lambda^*,\pi_{1:2}^*, v_{1:2}^* ) \geq \cL(x^*,\lambda,\pi_{1:2}, v_{1:2}).
	\end{equation}
	For any solution $x_t$, we define $\pi_{1,t}, \pi_{1, t+1}^0, \pi_{2,t}, \pi_{2,t+1}^0$ as 
	\begin{equation*}
	\begin{split}
	\pi_{2, t+1}  \in \argmax_{ \pi_2 \in \Pi_2} \{  \pi_2^\top  x_t   - f_2^*(  \pi_2)  \}, &  \ \ \pi_{2,t+1}^0  \in  \partial f_2(x_t, \xi_{2,t}^0), 
	\\
	\pi_{1, t+1}   \in \argmax_{\pi_1 \in \Pi_1} \{   \pi_1^\top   y_{t+1}    - f_1^*(\pi_1)  \} , &   \text{ and }\pi_{1, t+1}^0  = \nabla f_1(y_{t+1},\xi_{1,t}^0 ). 
	\end{split}
	\end{equation*}
	Consequently, we have 
	\begin{equation*}
	\begin{split}
	\cL_{f_2}(x_t, \pi_{2,t+1})  = f_2(x_t), 
	\ \ \EE[  \pi_{2,t+1}^0  \mid x_t]  \in   \partial f_2(x_t ),   \text{ and }\EE[ \pi_{1, t+1}^0 \mid y_{t+1}] = \nabla f_1(y_{t+1}).
	\end{split}
	\end{equation*}
	Then for the constraints $g_1, g_2$,
	we define 
	\begin{equation}\label{def:v2}
	v_{2, t+1} \in \argmax_{ v_2 \in \cV_2} \{ v_2^\top  x_t - g_2^*( v_2)  \}, 
	\end{equation}
	and observe that $v_{2,t+1}$ is the subgradient of $g_2$ at $x_t$. That is, $v_{2,t+1} \in \partial g_2(x_t)$. Correspondingly, we have 
	$$\cL_{g_2}(x_t, v_{2,t+1}) =  v_{2,  t+1}^\top  x_t -  g_2^*(v_{2,   t+1} )  = g_2(x_t).$$
	Meanwhile, we let the stochastic dual variable be $v_{2,t+1}^0 \in \partial g_2(x_t, \zeta_{2,t}^0)$, and observe that
	$$
	\EE[  v_{2,  t+1}^0   \mid x_t] \in \partial g_2(x_t). 
	$$
	Then, recall that  $z_{t+1} =  \frac{g_2(x_{t},\zeta_{2,t}^1) + \rho_{t} z_{t}}{1+\rho_{t}}$, taking $z_{t+1}$ as an estimator for $g_2(x_t)$, we have that 
	$$
	v_{1, t+1}  \in \argmax_{v_1 \in \cV_1} \bigl\{   v_1^\top   z_{t+1}   - g_1^*( v_1)  \bigr\}
	$$
	is the gradient of $g_1$ at $z_{t+1}$. In particular, we have 
	$$
	v_{1,t+1} = \nabla g_1(z_{t+1}).
	$$
	We also define a stochastic dual variable $ v_{1, t +1}^0  = \nabla g_1(z_{t+1},\zeta_{1,t}^0 )
	$ such that 
	$$
	\EE[ v_{1, t+1}^0  \mid z_{t+1}] = v_{1,t+1}. 
	$$
	After understanding the stochastic gradients under the stochastic sequential dual framework, the primal update step \eqref{eq:update_x_comp} can be viewed as 
	\begin{equation}\label{eq:update_x_comp_ssd}
	x_{t +1} = \argmin_{x\in \cX}  \Big \{ \big  \langle  \pi_{2, t+1}^0 \pi_{1, t+1}^0    +   v_{2,  t+1}^0  v_{1, t+1}^0 \lambda_t ,  x \big 
	\rangle   + \frac{\eta_t}{2} \| x - x_t\|^2 \Big \}.
	\end{equation}
	For the dual update step \eqref{eq:update_update_comp}, by the following lemma, we observe that $H_{t+1}$ is an unbiased estimator of $\cL_G(x_{t}, v_{1:2,t+1} ) $ with finite variance. This implies that \eqref{eq:update_update_comp} is equivalent to
	\begin{equation*}
	\lambda_{t+1}  \in \argmax_{ \lambda \in \RR_+^m } \Big \{ \lambda^\top H_{t+1}  - \frac{\alpha_{t}}{2} \| \lambda_t - \lambda \|^2 \Big \}, \quad \text{where} \quad  \EE[ H_{t+1} | x_{t}, v_{1,t+1}, v_{2,t+1}] = \cL_G(x_{t}, v_{1:2,t+1} ).
	\end{equation*}
	
	\begin{lemma}\label{lemma:H}
	Suppose Assumptions \ref{assumption:01}, \ref{assumption:03}, and \ref{assumption:compositional_g} hold, and   Algorithm~\ref{alg:02} generates $\{ (x_t, \lambda_t,\pi_{1:2,t},v_{1:2,t}) \}_{t=1}^N$. We have
	\begin{enumerate}
	\item[(a)] $H_{t}$ is an unbiased estimator of $ \cL_G(x_{t-1}, v_{1:2,t} )$ such that 
	\begin{equation*}
	\EE \big [ H_{t}  \mid x_{t-1}, v_{1,t}, v_{2,t}  \big ] = \cL_G(x_{t-1}, v_{1:2,t} ).
	\end{equation*}
	\item[(b)] The variance of $H_t$ is upper bounded by a  constant $\sigma_H >0$ such that 
	\begin{equation}\label{eq:g_var}
	\begin{split}
	& \Var \big ( H_t  \mid x_{t-1}, v_{1,t},v_{2,t} \big )
	\leq \sigma_H^2.
	\end{split}
	\end{equation}
	\end{enumerate}
	\end{lemma}
	\textit{Proof:}
	(a) By recalling the definition of $H_t$ \eqref{def:H} and noting that $g_1(z_{t},  \zeta_{1,t-1}^1)$, $\nabla g_1(z_{t},\zeta_{1,t-1}^1) $, and $g_2(x_{t-1}, \zeta_{2, t-1}^2)$ are independently sampled, we have 
	\begin{equation*}
	\begin{split}
	\EE \big [ H_{t} \mid  x_{t-1}, v_{1,t}, v_{2,t}  \big ]  
	& =  \EE \bigl[ g_1(z_{t},  \zeta_{1,t-1}^1)   -  [v_{1, t }^1 ]^\top  z_{t}   + [v_{1,t}^1]^\top  g_2(x_{t-1}, \zeta_{2,t-1}^2) \ \big| \   x_{t-1}, v_{1,t},v_{2,t}  \bigr] \\
	& =  g_1(z_{t})  + v_{1,t}^\top g_2(x_{t-1}) -  v_{1, t }^\top  z_{t}  =  v_{1, t}^\top   \cL_{g_2}(x_{t-1}, v_{2,t})    -   g_1^*(v_{1,t} )   \\
	& = \cL_G(x_{t-1}, v_{1:2,t } ),
	\end{split}
	\end{equation*}
	where the third equality holds since $v_{1,t} = \nabla g_1(z_t)$ and $v_{2,t} \in \partial g_2(x_{t-1})$ such that $g_1(z_t) = v_{1,t}^\top z_t - g_1^*(v_{1,t}) $ and $g_2(x_{t-1}) = \cL_{g_2} (x_{t-1}, v_{2,t})$.
	This proves part (a). 
	
	(b)
	Next, given $x_{t-1}, v_{1,t},v_{2,t}$, by  the definition of $H_{t}$  \eqref{def:H} and  the fact that $g_1(z_{t},  \zeta_{1,t-1}^1), v_{1,t}^1$, and $ g_2(x_{t-1}, \zeta_{2,t-1}^2) $ are independently sampled, we have that each component of $H_t$  preserves a bounded second moment. Therefore, $H_t$ has finite variance, completing the proof. 
	\QED 
	
	\subsection{Convergence Analysis}
	We then analyze the convergence behavior of  Algorithm~\ref{alg:02}. 
	With a slight abuse notation, we denote by $ ( x_t, \lambda_t,  \pi_{1:2,t}, v_{1:2,t} ) $ the solution generated by Algorithm~\ref{alg:02} at the $t$-th iteration and denote by $ ( x^*, \lambda,  \pi_{1:2}, v_{1:2} ) \in \cX \times \RR_+^m \times \Pi \times \cV $ a general feasible point. 
	We start our  analysis by decomposing the composite Lagrangian difference:
	\begin{equation}\label{eq:decom_comp}
	\begin{split}
	& \cL(x_t,\lambda,\pi,v) -  \cL(x^*, \lambda_t,  \pi_t , v_t) 
	 \\
	&  =   \cL_F(x_t, \pi_{1:2}) - \cL_F( x^*,\pi_{1:2,t}) + \lambda^\top  \big (  \cL_G(x_t, v_{1:2} ) - \cL_G(x_t, v_{1:2,t} )   \big) 
	 \\
	& \quad  + \lambda_t^\top \big (  \cL_G(x_t, v_{1:2,t} )   - \cL_G(x^*, v_{1:2,t})  \big )  
	 + (\lambda  - \lambda_t)^\top \cL_G(x_t, v_{1:2,t} ) \\
	 &  =   \cL_F(x_t, \pi_{1:2}) - \cL_F( x_t,\pi_{1:2,t}) + \lambda^\top  \big (  \cL_G(x_t, v_{1:2} ) - \cL_G(x_t, v_{1:2,t} )   \big) 
	 \\
	 & \quad + (\lambda  - \lambda_t)^\top \cL_G(x_t, v_{1:2,t} )  + \cL_F( x_t,\pi_{1:2,t })  -  \cL_F( x^*,\pi_{1:2,t })  
	 \\
	 & \quad +  \lambda_t^\top \big (  \cL_G(x_t, v_{1:2,t } )   - \cL_G(x^*, v_{1:2,t})  \big )
	 . 
	\end{split}
	\end{equation}
	By  Lemma \ref{lemma:H} that $H_t $ is an unbiased estimator for $ \cL_G(x_{t-1}, v_{1:2,t} )$ with finite variance, 
	we  first bound the term $ (\lambda  - \lambda_t)^\top \cL_G(x_t, v_{1:2,t} )$ in expectation in the next lemma.
	\begin{lemma}\label{lemma:comp_lambda_G}
	Suppose Assumptions \ref{assumption:01},  \ref{assumption:03}, and \ref{assumption:compositional_g} hold. 
	Suppose Algorithm~\ref{alg:02} generates  $\{ (x_t, \lambda_t,\pi_{1:2 , t},v_{1:2 , t} )\}_{t=1}^N$ by setting $\tau_t = \rho_t = t/2$, $\alpha_{t}  = \alpha > 0$, and $\eta_t =  \eta >0 $ for $t=1,2,\cdots, N$. Let $\lambda \in \RR_+^{m}$ be a {bounded} random dual variable such that $\| \lambda\| \leq M_{\lambda}$ uniformly, then there exists a constant $\sigma_H >0 $ such that for any $ K \leq N$, 
	\begin{equation*}
	\begin{split}
	&\EE \Big [   \sum_{t=1}^{K} (\lambda  - \lambda_{t})^\top  \cL_G(x_t, v_{1:2,t} ) \Big ]  +  \frac{\alpha }{2} \EE [    \| \lambda_{K} - \lambda  \|^2 ]   \\
	 &\quad  \leq   \sum_{t=1}^K\Big (  \frac{\sigma_{H}^2}{\alpha }  - \frac{ \alpha \EE[  \| \lambda_{t-1}  - \lambda_{t}\|^2 ] }{4}   + \frac{5C_{g_1}^2 C_{g_2}^2 \EE[  \|  \lambda  - \lambda_{t} \|^2 ]  }{2\eta } + \frac{\eta \EE [ \| x_{t} - x_{t-1}\|^2 ] }{10} \Big )  
	 \\
	& \qquad +  \sqrt{K} M_{\lambda}  \sigma_H + \frac{\alpha }{2}\EE[  \| \lambda_{0}- \lambda  \|^2].
	 \end{split}
	\end{equation*} 
	\end{lemma}
	
	Next, we  show the contraction of the composite Lagrangian difference. Note that we first assume that the dual variable $\lambda$ is bounded, and we justify this assumption in  Proposition~\ref{prop:dual_comp}.
	\begin{theorem}\label{thm:composition_01}
	Suppose Assumptions \ref{assumption:01}, \ref{assumption:03}, and \ref{assumption:compositional_g} hold, and  Algorithm~\ref{alg:02} generates $\{ (x_t, \lambda_t,\pi_{1:2 , t},v_{1:2 , t} )\}_{t=1}^N$ by setting $\tau_t = \rho_t = t/2$, $\eta_t =  \eta  $ and $\alpha_t = \alpha$ for $1\leq t \leq N$. Letting $\lambda \in \RR_+^m$ be a bounded random variable such that $\| \lambda \|\leq M_{\lambda}$, then for any integer $K \leq N$ and any feasible point $(x^*,\lambda,\pi_{1:2},v_{1:2})$, we have 
	\begin{equation*}
	\begin{split}
	& \sum_{t=1}^K \EE \Big ( \cL(x_t,\lambda,  \pi_{1:2} , v_{1:2}) -  \cL(x^*, \lambda_t,  \pi_{1:2, t} , v_{1:2 , t} )   \Big )  +  \frac{\alpha }{2} \EE [    \| \lambda_{K} - \lambda  \|^2 ] 
	\\
	&\quad \leq 
	\frac{40 K C_{g_1}^2 C_{g_2}^2  M_{\lambda}^2 }{ \eta }  +      2\sqrt{ K } L_{g_1} \sigma_{g_2}^2 M_{\lambda} +   3 \sqrt{K} C_{g_1}\sigma_{g_2} M_{\lambda}  + \frac{\eta }{2} \| x_{0}- x^*\|^2   
	\\
	& 
	\qquad    +   \frac{95 K C_{f_1}^2 C_{f_2}^2 }{2\eta   }    +  \frac{K}{\alpha } ( D_X^2 C_{g_1}^2 C_{g_2}^2 + \sigma_H^2 )
	+ \frac{\alpha }{2}\EE[  \| \lambda_{0}- \lambda  \|^2]    +  \sqrt{K} M_{\lambda} \sigma_H  
	\\
	& 
	\qquad 
	 +      2\sqrt{ K } L_{f_1} \sigma_{f_2}^2  +   3 \sqrt{K} C_{f_1}\sigma_{f_2} 
	+  \sum_{t=1}^K  \frac{5  C_{g_1}^2 C_{g_2}^2 }{2\eta }  \Big ( \EE[  \| \lambda  - \lambda_{t}\|^2 ]+ \EE [ \| \lambda_{t-1}\|^2]  \Big ),
	\end{split}
	\end{equation*}
	where $\sigma_H >0$ is provided in \eqref{eq:g_var}.
	\end{theorem}
	
	In the rest of this section, Proposition 3.1 demonstrates the boundedness of the sequence $\{\lambda_t\}$, and while it references Theorem 3.1, the theorem itself does not assume the boundedness of $\{\lambda_t \}$. In Theorem 3.1, we construct a random variable $\lambda$ without assuming that either $\{ \lambda_t\}$ or $\{ \| \lambda_t - \lambda\| \}$ is bounded. 
	
	We now show that the dual variables $\{ \lambda_t \}$ are bounded in the next proposition.

	\begin{proposition}\label{prop:dual_comp} 
	Suppose Assumptions~\ref{assumption:01}, \ref{assumption:03}, and \ref{assumption:compositional_g} hold. 
	Let $\{ \lambda_t \}_{t=1}^N$ be the sequence of dual variables generated by Algorithm \ref{alg:02} with  $\tau_t = \rho_t = t/2$, $\alpha_{t} = 2  \sqrt{N}$, and $\eta_t = \frac{15 C_{g_1}^2 C_{g_2}^2 \sqrt{N}}{2}$ for all $t\leq N$. Then for any $N \geq 2$ and any integer $K \leq N$, we have 
	$$
	\EE[ \| \lambda_K - \lambda^*\|^2 ] \leq 2R e^{2}, 
	$$
	where $ \sigma_H >0$  is provided in \eqref{eq:g_var},
	\begin{equation}\label{eq:R_comp}
	\begin{split}
	Q  & = 2L_{g_1} \sigma_{g_2}^2 \|\lambda^* \| +   3C_{g_1}\sigma_{g_2} \|\lambda^* \|  +  \|\lambda^* \|  \sigma_H     +        2L_{f_1} \sigma_{f_2}^2  +   3 C_{f_1}\sigma_{f_2} ,
	 \\
	 \text{ and } R & = \frac{7C_{f_1}^2 C_{f_2}^2 }{C_{g_1}^2 C_{g_2}^2} +   6  \|\lambda^* \|      +  \frac{ D_X^2 C_{g_1}^2 C_{g_2}^2 + \sigma_H^2  }{2}   +  \frac{15 C_{g_1}^2 C_{g_2}^2 }{4}  \| x_{0}- x^*\|^2     + \| \lambda_{0}- \lambda^*  \|^2      +  Q. 
	 \end{split}
	 \end{equation}
	\end{proposition}

	In the above, we show that the dual sequence $\{\lambda_t \}$ generated by Algorithm~\ref{alg:02} is also bounded in expectation, serving as a building block for analyzing the convergence behavior of the solution path~$\{ x_t\}$.

	Letting $\bar x_N = \frac{1}{N} \sum_{t=1}^N x_t$, by  the boundedness property, we derive the convergence rates for the objective optimality gap $F(\bar x_N) - F(x^*)$ and feasibility residual $G(\bar x_N)_+ = \max\{ G(x),0\}$ in the next theorem.
	\begin{theorem}\label{thm:comp_rate}
	Suppose Assumptions \ref{assumption:01}, \ref{assumption:03}, and \ref{assumption:compositional_g} hold, and  Algorithm~\ref{alg:02} generates $ \{ (x_t, \lambda_t,\pi_{1:2,t},v_{1:2,t}) \}_{t=0}^N $ by setting $\lambda_0 = 0$, $\tau_t = \rho_t = t/2$,  $\alpha_{t} = 2 \sqrt{N}$, and $\eta_t = \frac{15C_{g_1}^2 C_{g_2}^2 \sqrt{N}}{2}$  for $t = 1,\cdots, N$, and $N \geq 2$. Letting $\bar x_N = \frac{1}{N} \sum_{t=1}^N x_t$, we have 
	\begin{equation*}
	\begin{split}
	& \EE[ F(\bar x_N) - F(x^*)] \\
	& \leq   \frac{1}{\sqrt{N}}\Big ( 
	 \frac{15 C_{g_1}^2 C_{g_2}^2 }{4} \| x_{0}- x^*\|^2        +  \frac{ 7  C_{f_1}^2 C_{f_2}^2  }{C_{g_1}^2 C_{g_2^2}}  +  \frac{ D_X^2 C_{g_1}^2 C_{g_2}^2 + \sigma_H^2}{2 }   +      2 L_{f_1} \sigma_{f_2}^2  +   3 C_{f_1}\sigma_{f_2}  \Big ) \\
	& 
	\qquad 
	+  \frac{4}{3\sqrt{N}}   \Big (    \| \lambda^*\|^2  + 2R e^{2} \Big ),
	 \end{split}
	\end{equation*}
	and 
	\begin{equation*}
	\begin{split}
	& \EE [ \| G(\bar x_N)_{+} \|_2 ]  \\
	& \leq \frac{1}{\sqrt{N}} \Big ( \frac{16}{3}    \| \tilde \lambda\|^2   +      2L_{g_1} \sigma_{g_2}^2   \| \tilde \lambda\| +   3  C_{g_1}\sigma_{g_2}   \| \tilde \lambda\| +  \frac{15 C_{g_1}^2 C_{g_2}^2 }{4}  \| x_{0}- x^*\|^2     +  \frac{7  C_{f_1}^2 C_{f_2}^2 }{C_{g_1}^2 C_{g_2}^2}\Big )
	\\
	& 
	\quad + \frac{1}{\sqrt{N}} \Big (      \frac{D_X^2 C_{g_1}^2 C_{g_2}^2 + \sigma_H^2}{2 } 
	+ \| \tilde \lambda\|^2    +  \sigma_g   \| \tilde \lambda\|    +      2L_{f_1} \sigma_{f_2}^2  +   3  C_{f_1}\sigma_{f_2}   \Big )
	\\
	& 
	\quad 
	+   \frac{1 }{\sqrt{N} }  \Big ( \| \tilde \lambda  \| + \| \lambda^* \|^2  +  2R e^{2}   \Big ),
	 \end{split}
	\end{equation*}
	where $\tilde \lambda =  (\| \lambda^*\| +1)\frac{G(\bar x_N)_+}{ \| G(\bar x_N)_+ \|_2}$ is bounded such that $\| \tilde \lambda\|  = \| \lambda^*\| +1$  and $\| \tilde \lambda\|^2  = (\| \lambda^*\| +1)^2$, and $\sigma_H, R>0$ are  constants defined in \eqref{eq:g_var} and \eqref{eq:R_comp}, respectively. 
	\end{theorem}
	
	The above result indicates that Algorithm~\ref{alg:02} achieves the optimal $\cO(\frac{1}{\sqrt{N}})$ convergence rate for more complex compositional expected value constrained problems. Similar to Theorem~\ref{thm:rate}, the hidden constant within the above $\cO(\frac{1}{\sqrt{N}})$ convergence rate still grows polynomially with the number of constraints $m$, underscoring the efficiency of our algorithm in managing real-world applications that involve a large number of compositional constraints.

	\section{Numerical experiments} \label{sec:numerical}
	
	In this section, we investigate the empirical performance of the proposed algorithm through application to a portfolio optimization problem involving CVaR constraints. 
	This scenario corresponds to the single-level expected value constrained problem introduced in Section \ref{sec:applications}.
	
	Let $w \in \mathbb{R}^{d_x}$ denote random returns for ${d_x}$ assets, and let $x$ represent corresponding portfolio weights. We focus exclusively on long positions, constraining $x$ to the $({d_x}-1)$-dimensional probability simplex $\Delta= \{ x \in \RR_+^d: \sum_{i=1}^{d_x} x_i = 1\}$.
	We formulate the portfolio optimization problem with multiple CVaR constraints as follows:
	\begin{equation}\label{eq:Opt_numerical}
	  \min_{x \in \Delta} F(x), \quad
		  \text{s.t. } \mathrm{CVaR}_{\delta_i}(x) \le \gamma_i, \quad \forall i=1,\dots,m,
	\end{equation}
	where $\delta_i$ are CVaR constraint levels, $\gamma_i$ are their corresponding upper bounds, the objective function is given by $F(x) = \mathbb{E}_w[- w^{\top}x] + c  \mathbb{E}_w\left[(w^{\top}x-\mathbb{E}[ w^{\top}x])^4\right]$, and the CVaR of the portfolio's negative returns is defined as $\mathrm{CVaR}_{\delta_i}(x) = \mathrm{CVaR}_{\delta_i}(-w^{\top}x)$.

	We now reformulate problem \eqref{eq:Opt_numerical} into the standard structure from Section \ref{sec:single}. Define the following:
	\begin{align*}
		f_1(y,w) & = f_1(x,z,w) = - z + \xi (w^{\top} x-z)^4\in \mathbb{R}^1, \\
		f_2(x,\bm{u},w) & = [x,w^{\top} x]\in \mathbb{R}^{d+1}, \\
		g_i(x,u_i,w) & = u_i + \frac{1}{1-\delta_i} [-w^{\top} x - u_i]_+ - \gamma_i \in \mathbb{R}^1,
	\end{align*}
	where $z \in \mathbb{R}$, $\bm{u} \in \mathbb{R}^m$, and $y = (x,z)$. Let the expected functions be:
	\[
	f_1(x,z) = \mathbb{E}_w[f_1(x,z,w)], \quad f_2(x,\bm{u}) = \mathbb{E}_w[f_2(x,\bm{u},w)], \quad g_i(x,u_i) = \mathbb{E}_w[g_i(x,u_i,w)].
	\]
	Then, problem \eqref{eq:Opt_numerical} takes the form:
	\[
	\min_{x \in \Delta, u \in \mathbb{R}} f_1\circ f_2(x,u), \quad \text{s.t. } g_i(x,u_i) \leq 0, \quad \forall i=1,\dots,m.
	\]

\begin{figure}[htbp]
		\centering
		\includegraphics[width=0.45\linewidth]{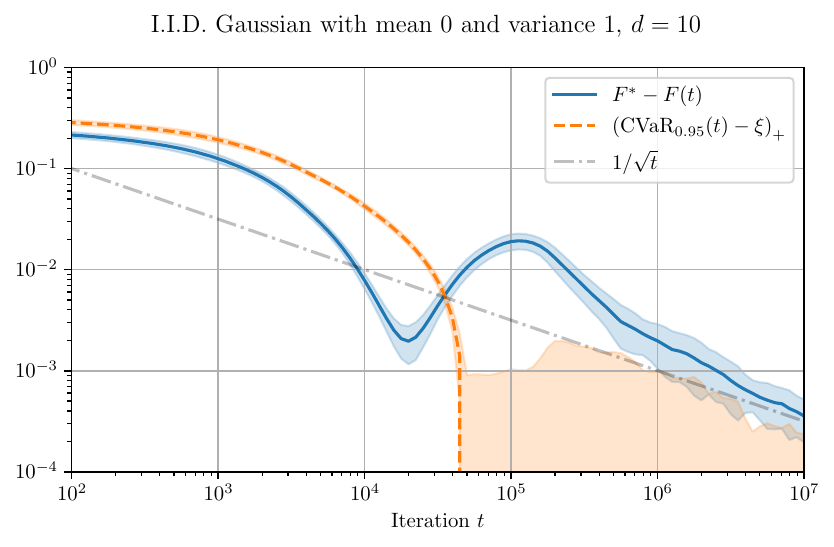}
		\includegraphics[width=0.45\linewidth]{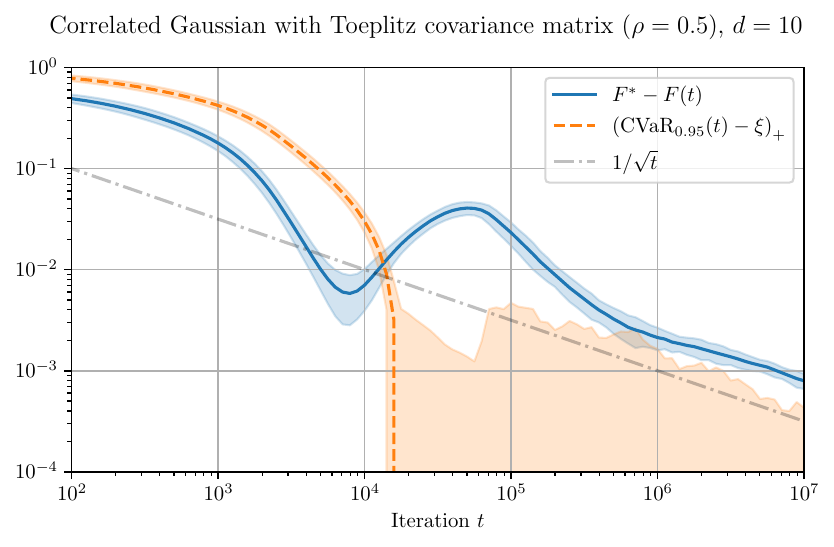}
		\includegraphics[width=0.45\linewidth]{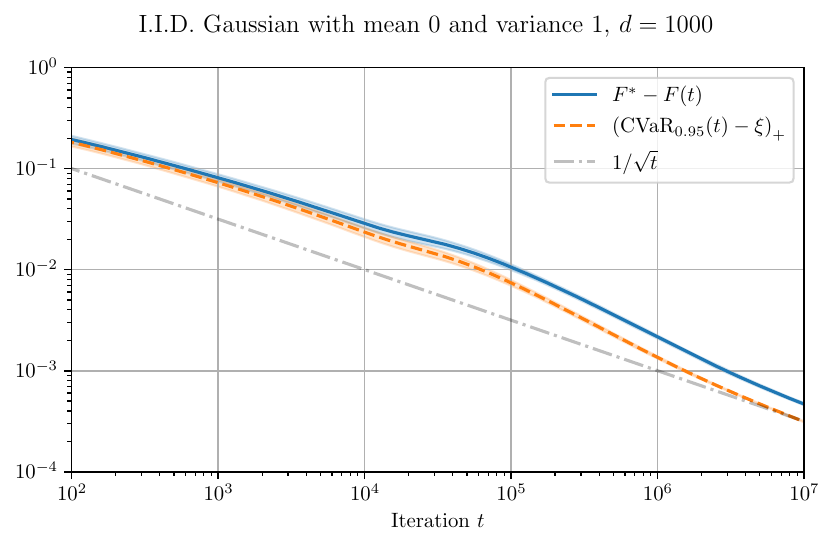}
		\includegraphics[width=0.45\linewidth]{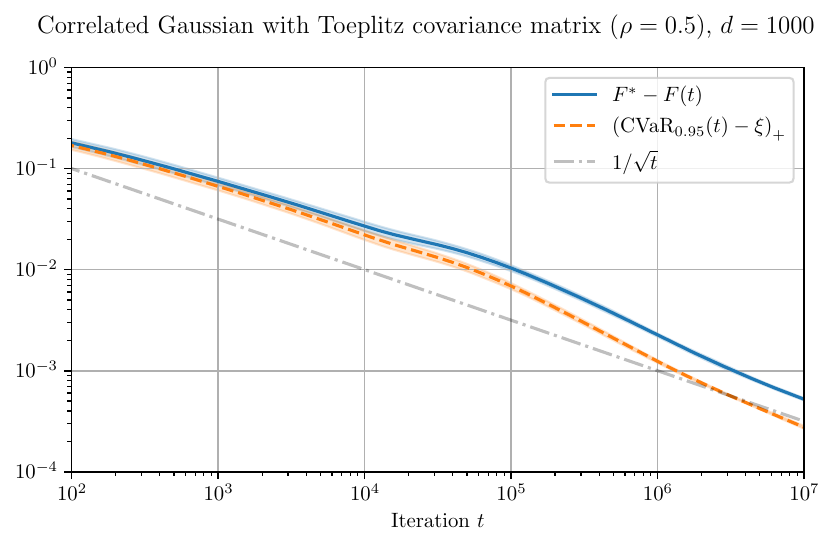}
		\caption{Empirical convergence of EC-SCGD algorithm for single CVaR constraint. Optimality gap and feasibility residual averaged over 10 runs with shaded 95\% confidence intervals.}
		\label{fig:num}
	\end{figure}

	\paragraph{Single CVaR constraint.}
	First, we analyze the single CVaR constraint scenario. Consider $w$ as a multivariate normal random vector with mean $\mu$ and covariance $\Sigma$. Random instances are generated with $\mu \sim \mathrm{Unif}[-1,1]^{d_x}$, using two covariance structures: independent case ($\Sigma = I_{d_x}$) and positively correlated Toeplitz matrix with entries $\Sigma_{ij} = 0.5^{|i-j|}$. Parameters are set as follows: risk-aversion parameter $c=0.5$, CVaR level $\delta=0.95$, and upper bound $\gamma = 0.6\mathrm{CVaR}_0 + 0.4\mathrm{CVaR}_{\min}, $
	where $\mathrm{CVaR}_0$ corresponds to the optimal portfolio without CVaR constraint and $\mathrm{CVaR}_{\min}$ is the minimum achievable CVaR.
	Algorithm \ref{alg:01} is implemented with step-sizes $\alpha_t = \max\{20d_x, 0.02 d_x\sqrt{t}\}, \quad \eta_t = 300\sqrt{t}, \quad \tau_t = 0.02t. $
	We consider dimensions $d_x = 10, 1000$ for $10^7$ iterations with 10 repetitions each.
	Results depicted in Figure \ref{fig:num} show the trajectories of optimality gaps $F(x^*) - F(x_t)$ and feasibility residuals $(\mathrm{CVaR}_{\delta}(x_t)-\gamma)_+$, averaged over 10 runs, with 95\% confidence intervals shaded. Observations confirm theoretical convergence rates from Theorem \ref{thm:rate}, showing slope close to $-1/2$, and feasibility residual quickly approaching zero.
	
\begin{figure}[htbp]
	  \centering
	  \includegraphics[width=0.45\linewidth]{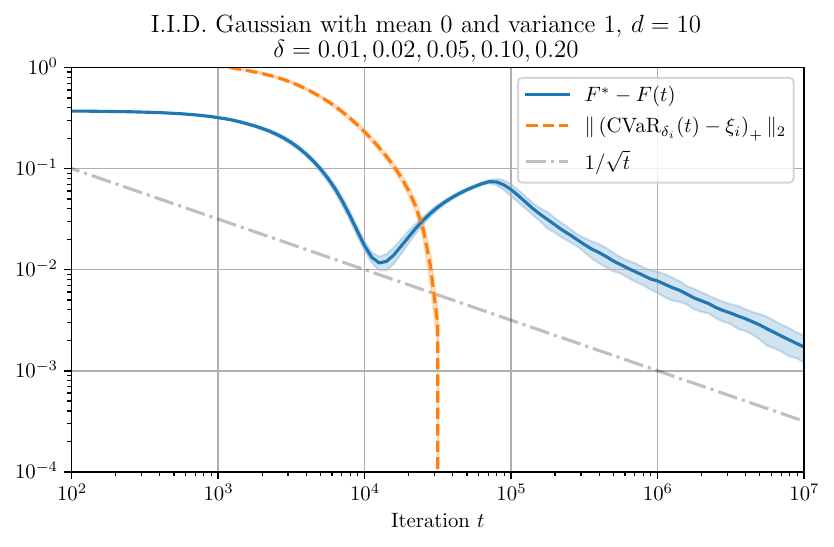}
	  \includegraphics[width=0.45\linewidth]{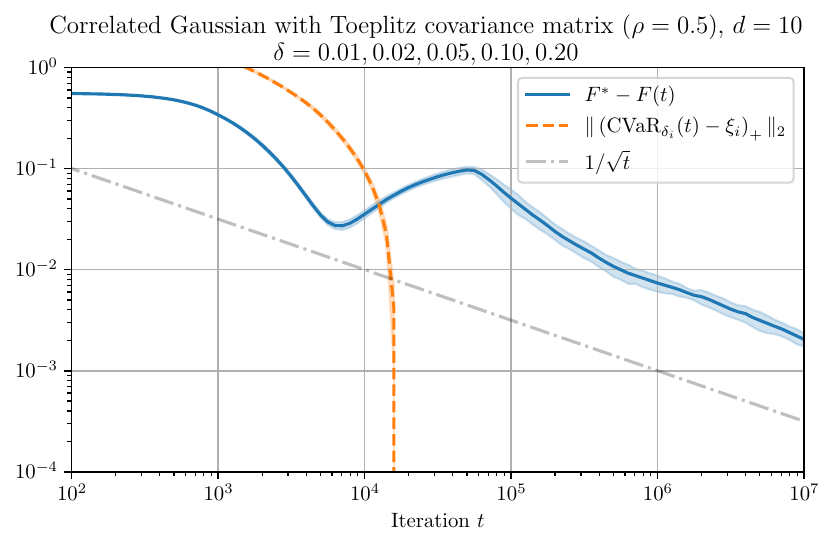}
	  \caption{Empirical convergence of EC-SCGD for multiple CVaR constraints. Optimality gap and feasibility residual averaged over 10 runs, with shaded 95\% confidence intervals.}
	  \label{fig:num_multi}
	\end{figure}

	\paragraph{Multiple CVaR constraints.} 
	We extend the analysis to multiple CVaR constraints, maintaining identical settings for noise distributions and objective function. The upper bounds are set as $\gamma_i = 0.6\mathrm{CVaR}_{0_i} + 0.4\mathrm{CVaR}_{\min_i},$ with CVaR levels $\delta_i = 0.01, 0.02, 0.05, 0.1, 0.2$, where $\mathrm{CVaR}_{0_i}$ is the CVaR at level $\delta_i$ of the optimal portfolio without the CVaR constraints, and $\mathrm{CVaR}_{\min_i}$ is the smallest achievable CVaR at level $\delta_i$.
	Algorithm \ref{alg:01} parameters remain unchanged, with dimension $d_x=10$ and $10^7$ iterations over 10 repetitions. 	
	Results presented in Figure \ref{fig:num_multi} show the trajectories of optimality gaps $F(x^*) - F(x_t)$ and feasibility residuals $\|(\mathrm{CVaR}_{\delta}(x_t) - \gamma)_+\|_2$, averaged over 10 runs, with 95\% confidence intervals shaded. We observe that the slope of the trajectory of the optimality gap closely matches $-1/2$, confirming the theoretical convergence rate established in Theorem \ref{thm:rate}.
	Additionally, the feasibility residual converges quickly to zero.
	
	\section{Conclusion}
	In this paper, we develop a novel model for stochastic compositional optimization (SCO) with single-level and two-level compositional expected-value constraints. Our formulation is handles multiple expected-value constraints and has broad applications in data-driven optimization and risk management. We propose a class of primal-dual algorithms, which we interpret within the stochastic sequential dual framework. For convex objectives and constraints, we demonstrate that our algorithms achieve the optimal convergence rate of $\cO(\frac{1}{\sqrt{N}})$ in both single-level and two-level compositional expected value constrained scenarios. These results match the optimal convergence rate for standard convex stochastic optimization and set new benchmarks. 
	
	\paragraph{Acknowledgement.}
	Shuoguang Yang received support from the Hong Kong Research Grants Council [ECS Grant 26209422]. 
	Wei You received support from the Hong Kong Research Grants Council [ECS Grant 26212320 and
	GRF Grant 16212823].

\clearpage

\bibliography{refs.bib}

\begin{thebibliography}{41}
\providecommand{\natexlab}[1]{#1}
\providecommand{\url}[1]{\texttt{#1}}
\expandafter\ifx\csname urlstyle\endcsname\relax
  \providecommand{\doi}[1]{doi: #1}\else
  \providecommand{\doi}{doi: \begingroup \urlstyle{rm}\Url}\fi

\bibitem[Ahmed et~al.(2007)Ahmed, {\c{C}}akmak, and Shapiro]{ahmed2007coherent}
Shabbir Ahmed, Ula{\c{s}} {\c{C}}akmak, and Alexander Shapiro.
\newblock Coherent risk measures in inventory problems.
\newblock \emph{European Journal of Operational Research}, 182\penalty0 (1):\penalty0 226--238, 2007.

\bibitem[Balasubramanian et~al.(2022)Balasubramanian, Ghadimi, and Nguyen]{balasubramanian2022stochastic}
Krishnakumar Balasubramanian, Saeed Ghadimi, and Anthony Nguyen.
\newblock Stochastic multilevel composition optimization algorithms with level-independent convergence rates.
\newblock \emph{SIAM Journal on Optimization}, 32\penalty0 (2):\penalty0 519--544, 2022.

\bibitem[Beck(2017)]{Beck2017First}
Amir Beck.
\newblock \emph{First-order methods in optimization}.
\newblock SIAM, 2017.

\bibitem[Bertsekas(1997)]{bertsekas1997nonlinear}
Dimitri~P Bertsekas.
\newblock Nonlinear programming.
\newblock \emph{Journal of the Operational Research Society}, 48\penalty0 (3):\penalty0 334--334, 1997.

\bibitem[Boob et~al.(2023)Boob, Deng, and Lan]{boob2023stochastic}
Digvijay Boob, Qi~Deng, and Guanghui Lan.
\newblock Stochastic first-order methods for convex and nonconvex functional constrained optimization.
\newblock \emph{Mathematical Programming}, 197\penalty0 (1):\penalty0 215--279, 2023.

\bibitem[Bruno et~al.(2016)Bruno, Ahmed, Shapiro, and Street]{bruno2016risk}
Sergio Bruno, Shabbir Ahmed, Alexander Shapiro, and Alexandre Street.
\newblock Risk neutral and risk averse approaches to multistage renewable investment planning under uncertainty.
\newblock \emph{European Journal of Operational Research}, 250\penalty0 (3):\penalty0 979--989, 2016.

\bibitem[Chen et~al.(2021)Chen, Sun, and Yin]{chen2021solving}
Tianyi Chen, Yuejiao Sun, and Wotao Yin.
\newblock Solving stochastic compositional optimization is nearly as easy as solving stochastic optimization.
\newblock \emph{IEEE Transactions on Signal Processing}, 69:\penalty0 4937--4948, 2021.

\bibitem[Cole et~al.(2017)Cole, Gin{\'e}, and Vickery]{cole2017does}
Shawn Cole, Xavier Gin{\'e}, and James Vickery.
\newblock How does risk management influence production decisions? evidence from a field experiment.
\newblock \emph{The Review of Financial Studies}, 30\penalty0 (6):\penalty0 1935--1970, 2017.

\bibitem[Ermoliev(1988)]{Ermoliev88}
Y.~Ermoliev.
\newblock Stochastic quasigradient methods.
\newblock In Y.~Ermoliev and R.~J.-B. Wets, editors, \emph{Numerical Techniques for Stochastic Optimization}, pages 141--186. Springer-Verlag, New York, 1988.

\bibitem[Finn et~al.(2017)Finn, Abbeel, and Levine]{finn2017model}
Chelsea Finn, Pieter Abbeel, and Sergey Levine.
\newblock Model-agnostic meta-learning for fast adaptation of deep networks.
\newblock In \emph{International conference on machine learning}, pages 1126--1135. PMLR, 2017.

\bibitem[Finn et~al.(2019)Finn, Rajeswaran, Kakade, and Levine]{Finn2019oninemetalearning}
Chelsea Finn, Aravind Rajeswaran, Sham Kakade, and Sergey Levine.
\newblock Online meta-learning.
\newblock In Kamalika Chaudhuri and Ruslan Salakhutdinov, editors, \emph{Proceedings of the 36th International Conference on Machine Learning}, volume~97 of \emph{Proceedings of Machine Learning Research}, pages 1920--1930. PMLR, 09--15 Jun 2019.

\bibitem[Ge et~al.(2015)Ge, Huang, Jin, and Yuan]{ge2015escaping}
Rong Ge, Furong Huang, Chi Jin, and Yang Yuan.
\newblock Escaping from saddle points---online stochastic gradient for tensor decomposition.
\newblock In Peter Grünwald, Elad Hazan, and Satyen Kale, editors, \emph{Proceedings of The 28th Conference on Learning Theory}, volume~40 of \emph{Proceedings of Machine Learning Research}, pages 797--842, Paris, France, 03--06 Jul 2015. PMLR.

\bibitem[Ghadimi et~al.(2020)Ghadimi, Ruszczynski, and Wang]{ghadimi2020single}
Saeed Ghadimi, Andrzej Ruszczynski, and Mengdi Wang.
\newblock A single timescale stochastic approximation method for nested stochastic optimization.
\newblock \emph{SIAM Journal on Optimization}, 30\penalty0 (1):\penalty0 960--979, 2020.

\bibitem[Haneveld and Van Der~Vlerk(2006)]{haneveld2006integrated}
Willem K~Klein Haneveld and Maarten~H Van Der~Vlerk.
\newblock Integrated chance constraints: reduced forms and an algorithm.
\newblock \emph{Computational Management Science}, 3:\penalty0 245--269, 2006.

\bibitem[Harvey et~al.(2010)Harvey, Liechty, Liechty, and M{\"u}ller]{harvey2010portfolio}
Campbell~R Harvey, John~C Liechty, Merrill~W Liechty, and Peter M{\"u}ller.
\newblock Portfolio selection with higher moments.
\newblock \emph{Quantitative Finance}, 10\penalty0 (5):\penalty0 469--485, 2010.

\bibitem[Lan(2020)]{lan2020first}
Guanghui Lan.
\newblock \emph{First-order and Stochastic Optimization Methods for Machine Learning}.
\newblock Springer Nature, 2020.

\bibitem[Lan and Zhou(2020)]{lan2020algorithms}
Guanghui Lan and Zhiqiang Zhou.
\newblock Algorithms for stochastic optimization with function or expectation constraints.
\newblock \emph{Computational Optimization and Applications}, 76\penalty0 (2):\penalty0 461--498, 2020.

\bibitem[Lan et~al.(2021)Lan, Romeijn, and Zhou]{lan2021conditional}
Guanghui Lan, Edwin Romeijn, and Zhiqiang Zhou.
\newblock Conditional gradient methods for convex optimization with general affine and nonlinear constraints.
\newblock \emph{SIAM Journal on Optimization}, 31\penalty0 (3):\penalty0 2307--2339, 2021.

\bibitem[Lemar{\'e}chal et~al.(1995)Lemar{\'e}chal, Nemirovskii, and Nesterov]{lemarechal1995new}
Claude Lemar{\'e}chal, Arkadii Nemirovskii, and Yurii Nesterov.
\newblock New variants of bundle methods.
\newblock \emph{Mathematical Programming}, 69\penalty0 (1):\penalty0 111--147, 1995.

\bibitem[Lin et~al.(2020)Lin, Nadarajah, Soheili, and Yang]{lin2020data}
Qihang Lin, Selvaprabu Nadarajah, Negar Soheili, and Tianbao Yang.
\newblock A data efficient and feasible level set method for stochastic convex optimization with expectation constraints.
\newblock \emph{Journal of Machine Learning Research}, 21\penalty0 (1), January 2020.
\newblock ISSN 1532-4435.

\bibitem[Liu et~al.(2022{\natexlab{a}})Liu, Wang, and So]{liu2022adaptive}
Huikang Liu, Xiaolu Wang, and Anthony Man-Cho So.
\newblock Adaptive coordinate sampling for stochastic primal--dual optimization.
\newblock \emph{International Transactions in Operational Research}, 29\penalty0 (1):\penalty0 24--47, 2022{\natexlab{a}}.

\bibitem[Liu et~al.(2022{\natexlab{b}})Liu, Cui, and Pang]{liu2022solving}
Junyi Liu, Ying Cui, and Jong-Shi Pang.
\newblock Solving nonsmooth and nonconvex compound stochastic programs with applications to risk measure minimization.
\newblock \emph{Mathematics of Operations Research}, 47\penalty0 (4):\penalty0 3051--3083, 2022{\natexlab{b}}.

\bibitem[Madavan and Bose(2021)]{madavan2021stochastic}
Avinash~N. Madavan and Subhonmesh Bose.
\newblock A stochastic primal-dual method for optimization with conditional value at risk constraints.
\newblock \emph{Journal of Optimization Theory and Applications}, 190\penalty0 (2):\penalty0 428--460, August 2021.

\bibitem[Nemirovski et~al.(2009)Nemirovski, Juditsky, Lan, and Shapiro]{nemirovski2009robust}
Arkadi Nemirovski, Anatoli Juditsky, Guanghui Lan, and Alexander Shapiro.
\newblock Robust stochastic approximation approach to stochastic programming.
\newblock \emph{SIAM Journal on Optimization}, 19\penalty0 (4):\penalty0 1574--1609, 2009.

\bibitem[Rockafellar et~al.(2000)Rockafellar, Uryasev, et~al.]{rockafellar2000optimization}
R~Tyrrell Rockafellar, Stanislav Uryasev, et~al.
\newblock Optimization of conditional value-at-risk.
\newblock \emph{Journal of Risk}, 2:\penalty0 21--42, 2000.

\bibitem[Ruszczy{\'n}ski and Shapiro(2006)]{ruszczynski2006optimization}
A.~Ruszczy{\'n}ski and A.~Shapiro.
\newblock Optimization of convex risk functions.
\newblock \emph{Mathematics of Operations Research}, 31\penalty0 (3):\penalty0 433--452, 2006.

\bibitem[Ruszczynski(2021)]{ruszczynski2021stochastic}
Andrzej Ruszczynski.
\newblock A stochastic subgradient method for nonsmooth nonconvex multilevel composition optimization.
\newblock \emph{SIAM Journal on Control and Optimization}, 59\penalty0 (3):\penalty0 2301--2320, 2021.

\bibitem[Slater(2014)]{slater2014lagrange}
Morton Slater.
\newblock Lagrange multipliers revisited.
\newblock In \emph{Traces and Emergence of Nonlinear Programming}, pages 293--306. Springer, 2014.

\bibitem[So et~al.(2009)So, Zhang, and Ye]{so2009stochastic}
Anthony Man-Cho So, Jiawei Zhang, and Yinyu Ye.
\newblock Stochastic combinatorial optimization with controllable risk aversion level.
\newblock \emph{Mathematics of Operations Research}, 34\penalty0 (3):\penalty0 522--537, 2009.

\bibitem[Tan et~al.(2018)Tan, Zhang, Ma, and Liu]{tan2018stochastic}
Conghui Tan, Tong Zhang, Shiqian Ma, and Ji~Liu.
\newblock Stochastic primal-dual method for empirical risk minimization with {O}(1) per-iteration complexity.
\newblock \emph{Advances in Neural Information Processing Systems}, 31, 2018.

\bibitem[Wang and Liu(2016)]{wang2016stochastic}
M.~Wang and J.~Liu.
\newblock A stochastic compositional gradient method using markov samples.
\newblock In \emph{Proceedings of the 2016 Winter Simulation Conference}, pages 702--713. IEEE Press, 2016.

\bibitem[Wang et~al.(2017{\natexlab{a}})Wang, Fang, and Liu]{wang2017stochastic}
M.~Wang, E.~X. Fang, and H.~Liu.
\newblock Stochastic compositional gradient descent: Algorithms for minimizing compositions of expected-value functions.
\newblock \emph{Mathematical Programming}, 161\penalty0 (1-2):\penalty0 419--449, 2017{\natexlab{a}}.

\bibitem[Wang et~al.(2017{\natexlab{b}})Wang, Liu, and Fang]{wang2017accelerating}
M.~Wang, J.~Liu, and E.~X. Fang.
\newblock Accelerating stochastic composition optimization.
\newblock \emph{The Journal of Machine Learning Research}, 18\penalty0 (1):\penalty0 3721--3743, 2017{\natexlab{b}}.

\bibitem[Yan et~al.(2019)Yan, Xu, Lin, Zhang, and Yang]{yan2019stochastic}
Yan Yan, Yi~Xu, Qihang Lin, Lijun Zhang, and Tianbao Yang.
\newblock Stochastic primal-dual algorithms with faster convergence than $ o (1/\sqrt{T} ) $ for problems without bilinear structure.
\newblock \emph{arXiv preprint arXiv:1904.10112}, 2019.

\bibitem[Yang et~al.(2019)Yang, Wang, and Fang]{yang2019multilevel}
S.~Yang, M.~Wang, and E.~X. Fang.
\newblock Multilevel stochastic gradient methods for nested composition optimization.
\newblock \emph{SIAM Journal on Optimization}, 29\penalty0 (1):\penalty0 616--659, 2019.

\bibitem[Yang et~al.(2024)Yang, Li, and Lan]{yang2024data}
Shuoguang Yang, Xudong Li, and Guanghui Lan.
\newblock Data-driven minimax optimization with expectation constraints.
\newblock \emph{Operations Research}, 2024.

\bibitem[Yu et~al.(2017)Yu, Neely, and Wei]{yu2017online}
Hao Yu, Michael Neely, and Xiaohan Wei.
\newblock Online convex optimization with stochastic constraints.
\newblock \emph{Advances in Neural Information Processing Systems}, 30, 2017.

\bibitem[Zafar et~al.(2019)Zafar, Valera, Gomez-Rodriguez, and Gummadi]{zafar2019fairness}
Muhammad~Bilal Zafar, Isabel Valera, Manuel Gomez-Rodriguez, and Krishna~P Gummadi.
\newblock Fairness constraints: A flexible approach for fair classification.
\newblock \emph{The Journal of Machine Learning Research}, 20\penalty0 (1):\penalty0 2737--2778, 2019.

\bibitem[Zhang and Xiao(2021)]{zhang2021multilevel}
Junyu Zhang and Lin Xiao.
\newblock Multilevel composite stochastic optimization via nested variance reduction.
\newblock \emph{SIAM Journal on Optimization}, 31\penalty0 (2):\penalty0 1131--1157, 2021.

\bibitem[Zhang and Lan(2024)]{zhang2024optimal}
Zhe Zhang and Guanghui Lan.
\newblock Optimal methods for convex nested stochastic composite optimization.
\newblock \emph{Mathematical Programming}, pages 1--48, 2024.

\bibitem[Zhang et~al.(2021)Zhang, Ahmed, and Lan]{zhang2019efficient}
Zhe Zhang, Shabbir Ahmed, and Guanghui Lan.
\newblock Efficient algorithms for distributionally robust stochastic optimization with discrete scenario support.
\newblock \emph{SIAM Journal on Optimization}, 31\penalty0 (3):\penalty0 1690--1721, 2021.

\end{thebibliography}
\bibliographystyle{plainnat}

\clearpage

\appendix
\part*{Appendix}

\section{Proof of Results in Section \ref{sec:single}}
\label{app:A}
Letting Assumptions \ref{assumption:01}, \ref{assumption:02}, and \ref{assumption:03} hold  and letting  $( x^*, \lambda,  \pi_{1:2}, v ) \in \cX \times \RR_+^m \times \Pi \times \cV $ be a general feasible point, we first introduce two technical lemmas to facilitate our analysis. 
\begin{lemma}[Lemma 3.8 of \citealt{lan2020first}]\label{lemma:three_point}
Assume function $g$ is $\mu$-strongly convex  with respect to some Bergman distance $V$, i.e., 
$g(y) - g(\bar y ) -g'(\bar y)^\top (y - \bar y ) \geq \mu V(y, \bar y)$. If $\hat y \in \argmin_{y \in Y} \{ \pi^\top y + g(y) + \tau V(y, \bar y)\}$, then 
\begin{equation*}
(\hat y - y)^\top \pi + g( \hat y) -g(y) \leq \tau V(y, \bar y ) -(\tau + \mu) V(y, \hat y) -\tau V(\hat y , \bar y).
\end{equation*}
\end{lemma}

\subsection{Proof of Theorem \ref{thm:01}}
\textit{Proof:}
Recall that $\cL(x,\pi_{1:2}, v,\lambda ) = \cL_{F} (x,\pi_{1:2}) + \lambda^\top \cL_{g} (x,v) $ for any $(x,\lambda,\pi_{1:2}, v ) \in \cX \times \RR_+^m \times \Pi \times \cV $, we consider $Q(z_t, z) $ and decompose it as 
\begin{equation}\label{eq:minmax_L}
\begin{split}
& Q(z_t, z)   =  \cL(x_t,\lambda,\pi_{1:2},v)-  \cL(x^*, \lambda_t,  \pi_{1:2,t} , v_t)   \\
& =  \cL_F(x_t,\pi_{1:2} ) - \cL_F( x^*,\pi_{1:2,t}) + \lambda^\top \cL_g(x_t, v ) - \lambda_t^\top \cL_g(x^*, v_t)  \\
& =  \cL_F(x_t, \pi_{1:2} ) - \cL_F( x^*,\pi_{1:2,t}) + \lambda^\top  \Big (   \cL_g(x_t, v ) - \cL_g(x_t, v_t )   \Big)     \\
& \quad + \lambda_t^\top  \Big (  \cL_g(x_t, v_t )   - \cL_g(x^*, v_t)  \Big ) + (\lambda  - \lambda_t)^\top \cL_g(x_t, v_t ) \\
& =  \cL_F(x_t, \pi_{1:2} )  - \cL_F(x_t, \pi_{1:2,t})  +    \cL_F(x_t, \pi_{1:2,t}) -    \cL_F(x^*, \pi_{1:2,t} )  \\
& \quad  + \lambda^\top  \Big (   \cL_g(x_t, v ) - \cL_g(x_t, v_t )   \Big)    
+ \lambda_t^\top \Big (  \cL_g(x_t, v_t )   - \cL_g(x^*, v_t)  \Big )  
 + (\lambda  - \lambda_t)^\top  \cL_g(x_t, v_t ) .
\end{split}
\end{equation}
We then provide bounds for the terms above in expectation. First, by  Lemma \ref{lemma:pi_ast_t} in Appendix Section \ref{app:A1}, we have  
\begin{equation}\label{eq:thm_single_01}
\begin{split}
& \sum_{t=1}^K \EE \Big [ \cL_F(x_t, \pi_{1:2})  - \cL_F(x_t, \pi_{1:2,t} )   \Big ]  \\
&\quad \leq     \frac{40 K C_{f_1}^2 C_{f_2}^2 }{\eta_0}  +   \sum_{t=1}^{K}  \frac{\eta_{t-1}  }{10} \EE [ \| x_{t}-x_{t-1}  \|^2 ] 
\quad +    \sqrt{ K } L_{f_1} \sigma_{f_2}^2  +   3 \sqrt{K} C_{f_1}\sigma_{f_2} .
\end{split}
\end{equation}
Secondly, by  Lemma \ref{lemma:F_lambda} in Appendix Section \ref{app:A2}, we obtain that 
\begin{equation}\label{eq:thm_single_02}
\begin{split}
  &\sum_{t=1}^K  \EE \Big [   \cL_F(x_t, \pi_{1:2,t}) -    \cL_F(x^*, \pi_{1:2,t} )   + \lambda_t^\top \Big (  \cL_g(x_t, v_t )   - \cL_g(x^*, v_t)  \Big )  \Big]   \\
& \quad \leq    \frac{\eta }{2} \| x_{0}- x^*\|^2 +   \sum_{t=1}^K \Big (   \frac{5C_{f_1}^2 C_{f_2}^2 }{  2\eta  }    -   \frac{3\eta }{10}  \EE [  \|x_{t-1} - x_{t} \|^2  
 ]  \Big )     \\
 & \qquad  + \sum_{t=1}^K \Big (  \frac{\alpha }{4} \EE[  \| \lambda_t - \lambda_{t-1}\|^2 ] +  \frac{D_X^2 C_g^2  }{\alpha } 
+   \frac{5 \EE[  \|  \lambda_{t-1}\|^2] C_g^2 }{2\eta }    \Big ).
\end{split} 
\end{equation}
Thirdly, by  Lemma \ref{lemma:L_G} in Appendix Section \ref{app:A3} and the condition that $\| \lambda\|\leq M_{\lambda}$, we obtain 
\begin{equation}\label{eq:thm_single_03}
\begin{split}
 \EE \Big [\sum_{t=1}^K \lambda^\top  \big (  \cL_g(x_t, v ) - \cL_g(x_t, v_t ) \big ) \Big ]
\leq \sum_{t=1}^K  \frac{10 C_g^2 M_{\lambda}^2 }{  \eta  } +  \frac{ \eta  \EE[ \| x_t - x_{t-1} \|^2]}{10}.
\end{split}
\end{equation}
Finally,  by  Lemma \ref{lemma:lambda} in Appendix Section \ref{app:A4}, we obtain that 
\begin{equation}\label{eq:thm_single_04}
\begin{split}
& \EE \Big [ \sum_{t=1}^K (\lambda  - \lambda_{t})^\top \cL_g(x_{t}, v_{t} )\Big ]  +  \frac{\alpha }{2}\EE[ \| \lambda_K - \lambda\|^2]\\
 &\quad \leq   \sqrt{K} M_{\lambda}  \sigma_g+\sum_{t=1}^K \Big ( \frac{5  C_g^2 \EE[  \|    \lambda  - \lambda_{t} \|^2 ] }{2\eta }  +  \frac{\eta }{10} \EE[\| x_{t-1} - x_t\|^2]  -   \frac{\alpha }{4} \EE[ \| \lambda_{t-1} - \lambda_{t} \|^2 ]  \Big ) \\
 & \qquad +   \frac{\alpha }{2} \EE [  \| \lambda_{0}- \lambda  \|^2   ]
+ \frac{K \sigma_{g}^2}{ \alpha }.
 \end{split}
\end{equation}

Summing \eqref{eq:minmax_L} over $t=1,2,\cdots, K$, taking expectations, and combining 
\eqref{eq:thm_single_01}, \eqref{eq:thm_single_02}, \eqref{eq:thm_single_03}, and \eqref{eq:thm_single_04}, 
 we conclude that
\begin{equation*}
\begin{split}
& \sum_{t=1}^K \EE \Big ( \cL(x_t,\lambda,  \pi_{1:2} , v) -  \cL(x^*, \lambda_t,  \pi_{1:2,t} , v_t)   \Big )  +  \frac{\alpha }{2}\EE[ \| \lambda_K - \lambda\|^2] 
\\
 & \leq      \frac{95 K C_{f_1}^2 C_{f_2}^2 }{2\eta } 
+  \sqrt{ K } L_{f_1} \sigma_{f_2}^2  +   3 \sqrt{K} C_{f_1}\sigma_{f_2}  + \frac{\eta }{2} \| x_{0}- x^*\|^2 
+     \frac{K}{\alpha } D_X^2 C_g^2 + \sqrt{K} M_{\lambda}  \sigma_g   \\ 
& \quad 
+   \frac{\alpha }{2}  \EE [ \| \lambda_{0}- \lambda  \|^2   ]
+ \frac{K \sigma_{g}^2}{ \alpha }
 +   \sum_{t=1}^K \Big ( \frac{5C_g^2  \EE[  \|  \lambda_{t-1}\|^2]  }{2\eta   }  +  \frac{10 C_g^2 M_{\lambda}^2 }{ \eta  } + \frac{5 C_g^2 \EE[  \| \lambda  - \lambda_{t}  \|^2 ] }{2\eta }  \Big ),
\end{split}
\end{equation*}
which completes the proof. 
\QED 

\begin{lemma}[Lemma 2 of \citealt{zhang2024optimal}] \label{lemma:sum_martingale}
Let $ \cF_t $ be a filtration and $\{ \delta_t \}_{t=1}^N$ be a martingale noise sequence such that $\delta_j \in \cF_{t+1}$ for $j=1,2,\cdots, t$,  $\EE[\delta_t | \cF_t] = 0$, and $\EE[ \| \delta_t\|^2 | \cF_t] \leq \sigma^2 $. For any random variable $\pi \in \Pi$ correlated with $\{ \delta_t \}_{t=1}^N$, suppose it is bounded such that $\| \pi \| \leq M_{\pi} $  uniformly, then
$$
\EE[\sum_{t=1}^N \hat \pi^\top \delta_t ] \leq \sqrt{N} M_{\pi}\sigma.
$$
\end{lemma}

\subsection{Lemma \ref{lemma:pi_ast_t} and Its Proof}\label{app:A1}

\noindent 
\begin{lemma}\label{lemma:pi_ast_t} Suppose Assumptions \ref{assumption:01}, \ref{assumption:02} and \ref{assumption:03} hold, and Algorithm \ref{alg:01}  generates
  $\{ (x_t,\lambda_t,\pi_{1:2,t}, v_t) \}_{t=1}^N$ by setting  $\tau_t = t/2$ and $\eta_t =  \eta $ for $t=0,1,\cdots, N$. Then for any integer $K \leq N$, we have
\begin{equation*}
\begin{split}
& \sum_{t=1}^K \EE \Big [ \cL_F(x_t, \pi_{1:2})  - \cL_F(x_t, \pi_{1:2,t} )   \Big ]  \\
& \quad \leq   \sum_{t=1}^{K}  \frac{\eta_{t-1}}{10} \EE [ \| x_{t}-x_{t-1}  \|^2 ]
+     \frac{40 K C_{f_1}^2 C_{f_2}^2 }{ \eta }  +      2\sqrt{ K } L_{f_1} \sigma_{f_2}^2  +   3 \sqrt{K} C_{f_1}\sigma_{f_2} .
\end{split}
\end{equation*}
\end{lemma}
\textit{Proof:} This result can be derived by combining Propositions 5 and 6 of \cite{zhang2024optimal}. Here we provide the proof sketch.
By decomposing $\cL_F(x_t, \pi_{1:2} )  - \cL_F(x_t, \pi_{1:2,t} )  $, we have 
\begin{align*}
& \cL_F(x_t, \pi_{1:2} )  - \cL_F(x_t, \pi_{1:2,t} )  
\\
& = \cL_{f_1}(x_t, \pi_1,\pi_2 )  -  \cL_{f_1}(x_t, \pi_1,\pi_{2,t} )   +  \cL_{f_1}(x_t, \pi_1,\pi_{2,t} )  - \cL_{f_1}(x_t, \pi_{1,t},\pi_{2,t} ) .
\end{align*}
We provide bound for each term after the above decomposition. First, we consider $\cL_{f_1}(x_t, \pi_1,\pi_2 )  -  \cL_{f_1}(x_t, \pi_1,\pi_{2,t} ) $ and obtain
\begin{align*} 
& \sum_{t=1}^K \EE \Big (  \cL_{f_1}(x_t, \pi_1,\pi_2 )  -  \cL_{f_1}(x_t, \pi_1,\pi_{2,t} )    \Big ) 
=  \sum_{t=1}^K  \EE[ \inner{ \cL_{f_2}(x_t, \pi_2) - \cL_{f_2}(x_t, \pi_{2,t}) }{ \pi_1}  ]
\\
& \quad =  \sum_{t=1}^K  \EE \Big [ \inner{ ( \pi_2 - \pi_{2,t})^\top x_{t } - f_2^*(\pi_2)  + f_2^*(\pi_{2,t}) } {\pi_1   }  \Big ]
\\
&\quad = \sum_{t=1}^K   \EE \Big [ \inner{ ( \pi_2 - \pi_{2,t})^\top x_{t -1} - f_2^*(\pi_2)  + f_2^*(\pi_{2,t}) } {\pi_1   }  \Big ] 
\\
& \qquad +  \sum_{t=1}^K  \EE \Big [ \inner{ ( \pi_2 - \pi_{2,t})^\top(x_t -  x_{t -1}) } {\pi_1   }  \Big ],
\end{align*}
By recalling \eqref{eq:pi_2t} that $\pi_{2, t} \in \argmax_{ \pi_2 \in \Pi_2} \{  \langle  \pi_2 ,  x_{t-1}  \rangle - f_2^*(  \pi_2)  \}$ and using Assumptions~\ref{assumption:03} that $\pi_1 \geq 0$ for non-affine $f_2$, we obtain for both affine ($\pi_{2}$ is a constant) and non-affine $f_2$, 
$$
\EE \Big [ \inner{ ( \pi_2 - \pi_{2,t})^\top x_{t -1} - f_2^*(\pi_2)  + f_2^*(\pi_{2,t}) } {\pi_1   }  \Big ]  \leq 0,
$$
which further implies that 
\begin{equation}\label{eq:Lf_1}
\begin{split}
& \sum_{t=1}^K \EE \Big (  \cL_{f_1}(x_t, \pi_1,\pi_2 )  -  \cL_{f_1}(x_t, \pi_1,\pi_{2,t} )    \Big ) 
 \leq  \sum_{t=1}^{K} \EE [   \inner{ (\pi_2 - \pi_{2, t})  \pi_1  }{x_{t} - x_{t-1}}] 
\\
&
\quad \leq \sum_{t=1}^{K}   \Big [  \frac{5 \EE [ \|  (\pi_2 - \pi_{2,t})   \pi_1 \|^2]  }{\eta }  +   \frac{\eta }{20} \EE [ \| x_{t}-x_{t-1}  \|^2 ]\Big ] 
 \\
 & \quad \leq \sum_{t=1}^{K}   \Big [  \frac{20  C_{f_1}^2 C_{f_2}^2  }{\eta }  +   \frac{\eta }{20} \EE [ \| x_{t}-x_{t-1}  \|^2 ]\Big ]  
.
\end{split}
\end{equation}
The next part is similar to Proposition 6 of \cite{zhang2024optimal}. Consider the following decomposition: 
\begin{align}
\begin{split}\label{eq:pfA3-1}
        &\EE[\cL_{f_1}(x_t, \pi_1, \pi_{2, t}) - \cL_{f_1}(x_t, \pi_{1, t}, \pi_{2, t})]   = \EE[\inner{\pi_1 - \pi_{1, t}}{f_2(x_t)} - \{f_1^*(\pi_1) - f_1^*(\pi_{1, t})\}]\\
    &\ \ = \EE[\inner{\pi_1 - \pi_{1, t}}{f_2(x_{t-1}, \xi_{2,t}^1)} - \{f_1^*(\pi_1) - f_1^*(\pi_{1, t})\}] \\
    & \ \quad + \EE[\inner{\pi_{1, t} - \pi_{1, t-1}}{f_2(x_{t-1})- f_2(x_{t-1}, \xi_{2,t}^1)}] + \EE[\inner{\pi_{1} - \pi_{1, t}}{f_2(x_t) - f_2(x_{t-1})}]\\
    & \ \quad + \EE[\inner{\pi_1}{f_2(x_{t-1}) - f_2(x_{t-1}, \xi_{2,t}^t)}] +\underbrace{\EE[\inner{\pi_{1, t-1}}{f_2(x_{t-1}) - f_2(x_{t-1}, \xi_{2,t}^1)}].}_{= 0}
\end{split}
\end{align}
Among these terms, the update rule of $\pi_{1, t}$, the strong convexity of $f_1^*$ with respect to $D_{f^*_1}$, and the choice of $\tau_t$ imply 
\begin{align*}
    &\EE[\sum_{t=1}^{K}\inner{\pi_1 - \pi_{1, t}}{f_2(x_{t-1}, \xi_{2,t}^1)} - \{f_1^*(\pi_1) - f_1^*(\pi_{1, t})\}] \\
    &\quad \leq \EE[\sum_{t=1}^K \tau_t  D_{f_1^*}(\pi_1; \pi_{1, t-1}) -  (\tau_t + 1)D_{f_1^*}(\pi_1; \pi_{1, t-1})- \tau_t D_{f_1^*}(\pi_{1, t}; \pi_{1, t-1})\}] \\
    &\quad \leq -\EE[\sum_{t=1}^K \tau_t D_{f_1^*}(\pi_{1, t}; \pi_{1, t-1})\}].
\end{align*}
The $1/L_{f_1}$-strong convexity of $D_{f_1^*}$, the choice of $\tau_t$ and the Cauchy-Schwartz inequality imply
\begin{align*}
 \EE &\sum_{t=1}^K \left[ (\pi_{1,t} - \pi_{1,t-1})(f_2(x_{t-1}, \xi_{2, t}^1) - f_2(x_{t-1}) ) - \tau_t \cD_{f_1^*}(\pi_{1,t}; \pi_{1,t-1})\right] \\
 &\leq \sum_{t=2}^{K} \frac{L_{f_1}\sigma_{f_2}^2}{\tau_t} + \sqrt{\EE [  \norm{\pi_{1,1} - \pi_{1,0}}^2] }\sqrt{\EE [\norm{f_2(x_{0}, \xi_{2, 0}^1) - f_2(x_{0})}^2]}\\
 &\leq 2 L_{f_1} \sigma^2_{f_2} \log K + 2 C_{f_1} \sigma_{f_2}\\
 &\leq 2\sqrt{K}  L_{f_1} \sigma^2_{f_2} + 2 \sqrt{K} C_{f_1} \sigma_{f_2}.
 \end{align*} 
 Lemma 2 of \cite{zhang2024optimal} implies that 
 $$\EE \left [\sum_{t=1}^K\inner{\pi_1}{f_2(x_{t-1}) - f_2(x_{t-1}, \xi_{2,t}^t)} \right ] \leq \sqrt{K} C_{f_1} \sigma_{f_2}.$$
 Similar to \eqref{eq:Lf_1}, we also have 
 \begin{align*}
     \EE[\sum_{t=1}^K \inner{\pi_{1} - \pi_{1, t}}{f_2(x_t) - f_2(x_{t-1})}] \leq \sum_{t=1}^{K}   \Big [  \frac{20  C_{f_1}^2 C_{f_2}^2  }{\eta }  +   \frac{\eta }{20} \EE [ \| x_{t}-x_{t-1}  \|^2 ]\Big ]. 
 \end{align*}
  Combining them with the decomposition in \eqref{eq:pfA3-1}, we obtain
\begin{align*}
& \sum_{t=1}^K  \EE  \Big (  \cL_{f_1}(x_t, \pi_1,\pi_{2,t} )  - \cL_{f_1}(x_t, \pi_{1,t},\pi_{2,t} ) \Big ) 
\\
& \quad \leq  \sum_{t=1}^{K}   \Big [  \frac{20 C_{f_1}^2 C_{f_2}^2 }{\eta }  +   \frac{\eta }{20} \EE [ \| x_{t}-x_{t-1}  \|^2 ]\Big ] 
+    2\sqrt{ K } L_{f_1} \sigma_{f_2}^2  +   3 \sqrt{K} C_{f_1}\sigma_{f_2}.
\end{align*}
The desired result then follows from adding to preceding inequality to \eqref{eq:Lf_1}.  
\QED 

\subsection{Lemma \ref{lemma:F_lambda} and Its Proof}\label{app:A2}
\begin{lemma}\label{lemma:F_lambda}
Suppose Assumptions \ref{assumption:01}, \ref{assumption:02}, and \ref{assumption:03} hold, and Algorithm \ref{alg:01} generates $\{ (x_t,\lambda_t,\pi_{1:2,t}, v_t) \}_{t=1}^N$ by setting $\tau_t = t/2$, $\alpha_{t} = \alpha $, and $\eta_t =  \eta $ for all $0 \leq t \leq N$. We have 
\begin{equation*}
\begin{split}
  &\sum_{t=1}^K \EE \Big [  \cL_F(x_t, \pi_{1:2,t} ) -    \cL_F(x^*, \pi_{1:2,t}) + \lambda_t^\top  \big (   \cL_g(x_t, v_t )   - \cL_g(x^*, v_t)  \big  )     \Big]   \\
&\quad \leq   \frac{\eta }{2} \EE [ \| x_{t-1}- x^*\|^2 ] - \sum_{t=1}^K   \frac{3\eta }{10}  \EE \Big [  \|x_{t-1} - x_{t} \|^2  
 \Big ]      +   \sum_{t=1}^K   \frac{5 C_{f_1}^2 C_{f_2}^2 }{2  \eta   }    \\
& \qquad 
+  \sum_{t=1}^K  \frac{\alpha }{4} \EE[  \| \lambda_t - \lambda_{t-1}\|^2 ] + \sum_{t=1}^K  \frac{1}{\alpha } D_X^2 C_g^2  
+   \sum_{t=1}^K   \frac{5 \EE[  \|  \lambda_{t-1}\|^2] C_g^2 }{2\eta  } ] .
\end{split}
\end{equation*}
\end{lemma}
\textit{Proof.}
Recall Algorithm~\ref{alg:01}, for $t=1,2,\cdots, K$, we note that $\pi_{2,t} \pi_{1,t} , v_t \lambda_t \in \RR^{d_x}$ and obtain  
\begin{equation*}
\begin{split}
  & \cL_F(x_t, \pi_{1:2,t} ) -    \cL_F(x^*, \pi_{1:2,t} ) +  \lambda_t^\top  \big (   \cL_g(x_t, v_t )   - \cL_g(x^*, v_t)  \big  )   \\
  &\quad = \inner{  \pi_{2,t}  \pi_{1,t}  }{ x_t - x^*  }+  \inner{  v_t  \lambda_t  }{   x_t - x^*} \\
  & \quad= \inner{ \pi_{2,t}^0    \pi_{1,t}^0   }{ x_t - x^*  }+  \inner{ v_t^0  \lambda_{t-1} }{  x_t - x^*}    \\
  &\qquad +\inner{  \pi_{2,t} \pi_{1,t}  -  \pi_{2,t}^0     \pi_{1,t}^0    }{ x_t - x^* }  +  \inner{  v_t  \lambda_t - v_t^0 \lambda_{t-1} }{  x_t - x^*},
\end{split}
\end{equation*}
and 
$$
x_t = \argmin_{x\in \cX} \Big \{ \inner { \pi_{2,t}^0  \pi_{1,t}^0  +  v_t^0   \lambda_{t-1} }{  x  }+   \frac{\eta_{t-1}}{2} \| x_{t-1}- x\|^2 \Big  \} .
$$
By applying Lemma \ref{lemma:three_point} to the above update rule, since $x^* \in \cX$, we have 
\begin{equation*}
\begin{split}
&   \inner{  \pi_{2,t}^0   \pi_{1,t}^0 +   v_t^0  \lambda_{t-1} }{ x_t - x^* }    \leq  \frac{\eta }{2} \| x_{t-1}- x^*\|^2 - \frac{\eta }{2} \|x_{t-1} - x_{t} \|^2 - \frac{\eta }{2} \|x_{t} -x^* \|^2 ,
\end{split}
\end{equation*}
which further yields 
\begin{equation}\label{eq:F_lambda}
\begin{split}
  & \cL_F(x_t, \pi_{1:2,t}  ) -    \cL_F(x^*, \pi_{1:2,t} ) +  \lambda_t^\top  \big (    \cL_g(x_t, v_t )   - \cL_g(x^*, v_t) \big )   \\
& \quad\leq
 \frac{\eta }{2} \| x_{t-1} - x^*\|^2 - \frac{\eta }{2} \|x_{t-1} - x_{t} \|^2 - \frac{\eta }{2} \|x_{t} -x^* \|  \\
  & \qquad + \underbrace{  \inner{ \pi_{2,t} \pi_{1,t}  -\pi_{2,t}^0   \pi_{1,t}^0 }{ x_t - x^* }  }_{\Delta_0^t} + \underbrace{ \inner{v_t    \lambda_t  - v_t^0 \lambda_{t-1}  }{ x_t - x^* } }_{\Lambda_0^t}.
\end{split}
\end{equation}
Then, we decompose $\Delta_0^t$ by 
\begin{equation*}
\begin{split}
\Delta_0^t = \big (  \pi_{2,t} \pi_{1,t}  -\pi_{2,t}^0  \pi_{1,t}^0  \big) (x_{t-1} - x^*)  + \big ( \pi_{2,t} \pi_{1,t}   - \pi_{2,t}^0  \pi_{1,t}^0   \big) (x_{t} - x_{t-1}). 
\end{split}
\end{equation*}
We note that  given $x_{t-1},y_t,$
$$
\EE \big [\pi_{2,t}   \pi_{1,t}  -  \pi_{2,t}^0 \pi_{1,t}^0  \big  | x_{t-1},y_t  \big ] = 0,
$$
and 
\begin{equation*}
\begin{split}
& \EE \big [  \| \pi_{2,t}   \pi_{1,t} -  \pi_{2,t}^0  \pi_{1,t}^0 \|^2 |  x_{t-1},y_t   \big ] \\
& \quad\leq  \EE \big [  \| \pi_{1,t} - \pi_{1,t}^0   \|^2  \big |   x_{t-1},y_t    \big ]  \EE \big [  \| \pi_{2,t}^0   \|^2  \big |x_{t-1} \big ] + \EE \big [  \| (\pi_{2,t}^0 - \pi_{2,t})\pi_{1,t} \|^2 |    x_{t-1},y_t   \big ]  \\
&  \quad\leq 2 C_{f_1}^2 C_{f_2}^2.
\end{split}
\end{equation*}
By the independence between $(x_{t-1} - x^*)$ and the mean-zero term $ (\pi_{2,t} \pi_{1,t}  - \pi_{2,t}^0  \pi_{1,t}^0 ) $, 
we have $\EE \Big [\big ( \pi_{2,t}  \pi_{1,t} -  \pi_{2,t}^0  \pi_{1,t}^0 \big)^\top (x_{t-1} - x^*)  \Big] = 0$. Moreover, to handle the  possible correlation between $(x_{t} - x_{t-1})$ and $(\pi_{2,t} \pi_{1,t}   -  \pi_{2,t}^0 \pi_{1,t}^0 )$, we utilize the fact that $ab \leq \frac{5a^2}{2} + \frac{b^2}{10}$ for all $a,b\in\RR$, and obtain 
\begin{equation}\label{eq:Delta_0}
\begin{split}
& \EE[\sum_{t=1}^K \Delta_0^t ] = \EE \Big [ \sum_{t=1}^K \inner{  \pi_{2,t} \pi_{1,t}  - \pi_{2,t}^0    \pi_{1,t}^0  }{  x_{t} - x_{t-1}  } \Big ]
\\
& \quad \leq   \sum_{t=1}^K  \Big ( \frac{5 C_{f_1}^2 C_{f_2}^2 }{2 \eta } + \frac{ \eta  \EE[ \| x_t - x_{t-1}\|^2]}{10 } \Big ).
 \end{split}
\end{equation}
We then decompose $\Lambda_0^t$ as
\begin{equation*}
\begin{split}
\Lambda_0^t = \underbrace{ \inner{ v_t  (\lambda_t - \lambda_{t-1} ) }{  x_t - x^* }  }_{\Lambda_{1}^t }  +  \underbrace{ \inner{ ( v_t - v_t^0 ) \lambda_{t-1}  }{ x_t - x^* }  }_{\Lambda_{2}^t }.
\end{split}
\end{equation*}
By the fact that $2ab\leq a^2 + b^2$ and   Assumptions~\ref{assumption:01} and  \ref{assumption:02} that $\| x_t - x^*\|\leq D_X$ and $\| v_t\|\leq C_g$, for any $\alpha >0 $, we have 
\begin{equation} \label{eq:Lambda_1}
\begin{split}
\EE [ \Lambda_{1}^t] & \leq  \frac{\alpha }{4} \EE[  \| \lambda_t - \lambda_{t-1}\|^2 ] + \frac{1}{ \alpha } \EE [ \| v_t^\top (x_t - x^*)  \|^2] 
\\
&
\leq   \frac{\alpha }{4} \EE[  \| \lambda_t - \lambda_{t-1}\|^2 ] + \frac{1}{\alpha } D_X^2 C_g^2. 
\end{split}
\end{equation}
Meanwhile, since $\EE[ v_t - v_t^0] = 0$ and $\EE[ \| v_t - v_t^0  \|^2 ] \leq C_g^2 $, we obtain 
\begin{equation}\label{eq:Lambda_2}
\begin{split}
\EE[\sum_{t=1}^K \Lambda_2^t] 
& =  \EE \Big [\sum_{t=1}^K  \inner{ (v_t - v_t^0 ) \lambda_{t-1} }{ x_{t-1} - x^* }   \Big ] + \EE \Big [\sum_{t=1}^K \inner{  (v_t - v_t^0 ) \lambda_{t-1} }{ x_{t} - x_{t-1} }  \Big  ] \\
&
=  \sum_{t=1}^K \EE \Big [   \inner{ (v_t - v_t^0 )  \lambda_{t-1} }{  x_{t} - x_{t-1} } \Big ] \\
&
\leq  \sum_{t=1}^K \EE \Big [ \frac{5\|   (v_t - v_t^0 ) \lambda_{t-1}  \|^2}{2\eta } \Big ] + \sum_{t=1}^K \frac{\eta \EE[  \|x_{t} - x_{t-1} \|^2] }{10}
\\
&
\leq \sum_{t=1}^K  \frac{5C_g^2  \EE[  \|  \lambda_{t-1}\|^2]  }{2\eta } ] + \sum_{t=1}^K \frac{\eta \EE[  \|x_{t} - x_{t-1} \|^2] }{10}
,
\end{split}
\end{equation}
where the second equality holds by the independence between $(v_t - v_t^0)\lambda_{t-1} $ and $(x_{t-1} - x^*)$, and the second last inequality holds due to the fact that $ab \leq \frac{5a^2}{2} + \frac{b^2}{10}$. Taking expectations on both sides of \eqref{eq:F_lambda}, summing over $t=1,\cdots, N$, and 
substituting \eqref{eq:Delta_0}, \eqref{eq:Lambda_1}, and \eqref{eq:Lambda_2} into it, we obtain
\begin{equation*}
\begin{split}
  &\sum_{t=1}^K\EE \Big [  \cL_F(x_t, \pi_{1:2,t} ) -    \cL_F(x^*, \pi_{1:2,t} ) + \lambda_t^\top  \Big (  \cL_g(x_t, v_t )   - \cL_g(x^*, v_t)  \Big )  \Big]   \\
& \quad \leq
\sum_{t=1}^K  \EE \Big [  \frac{\eta  }{2} \| x_t - x^*\|^2 - \frac{\eta  }{2} \|x_t - x_{t+1} \|^2 - \frac{\eta  }{2} \|x_{t+1} -x^* \|  + \Delta_0^t + \Lambda_1^t  + \Lambda_2^t  \Big ] \\
&\quad \leq    \frac{\eta }{2} \EE [ \| x_{t-1}- x^*\|^2 ] - \sum_{t=1}^K   \frac{3\eta }{10}  \EE \Big [  \|x_{t-1} - x_{t} \|^2  
 \Big ]      +   \sum_{t=1}^K   \frac{5 C_{f_1}^2 C_{f_2}^2 }{ 2 \eta_{t-1}  }    \\
& \qquad 
+  \sum_{t=1}^K  \frac{\alpha }{4} \EE[  \| \lambda_t - \lambda_{t-1}\|^2 ] + \sum_{t=1}^K  \frac{D_X^2 C_g^2 }{\alpha }  
+   \sum_{t=1}^K   \frac{5 \EE[  \|  \lambda_{t-1}\|^2] C_g^2 }{2\eta } ]   ,
\end{split} 
\end{equation*}
where the last inequality holds by the fact that $ \eta_t =  \eta  $ for $t=1,2,\cdots ,N$. This completes the proof.
\QED 

\subsection{Lemma \ref{lemma:L_G} and Its Proof}\label{app:A3}
\begin{lemma}\label{lemma:L_G}
Suppose Assumptions \ref{assumption:01}, \ref{assumption:02} and  \ref{assumption:03} hold, and Algorithm \ref{alg:01} generates  $ \{ (x_t,\lambda_t,\pi_t, v_t) \}_{t=1}^N$ by setting $\tau_t = t/2$, $\alpha_{t} = \alpha $, and $\eta_t =  \eta $ for all $0 \leq t \leq N$. Let $\lambda \in \RR_+^m$ be a nonnegative random variable whose norm is bounded such that $\| \lambda \| \leq M_{\lambda}$ uniformly. Then for any $K \leq N$, we have 
\begin{equation*}
\begin{split}
 \EE \Big [\sum_{t=1}^K \lambda^\top \big (  \cL_g(x_t, v ) - \cL_g(x_t, v_t ) \big ) \Big ] \leq  \frac{10M_{\lambda}^2 C_g^2   }{  \eta   } +  \frac{  \eta  \EE [ \| x_t - x_{t-1} \|^2 ] }{10 }.
\end{split}
\end{equation*}
\end{lemma}
\textit{Proof.}
Recall that the update rule is that 
$$
v_t \in \argmax_{v\in \cV }\Big \{ \inner{ v  }{ x_{t-1} } -g^*(v) \Big \}.
$$
By the fact that $\lambda \geq 0$, 
we have 
\begin{equation*}
\begin{split}
& \lambda^\top  \Big (\cL_g(x_t, v ) - \cL_g(x_t, v_t ) \Big ) 
\\
& =     \lambda^\top  \Big ( (v - v_t)^\top x_{t-1}  - (g^*(v) - g^*(v_t)) + (v - v_t)^\top (x_t - x_{t-1} ) \Big )  \\
& \leq \lambda^\top  \Big ( (v - v_t)^\top (x_t - x_{t-1} ) \Big ) .
\end{split}
\end{equation*}
By taking expectations on both sides of the above inequality, 
we have 
\begin{equation*}
\begin{split}
& \EE \left  [ \lambda^\top  \big (\cL_g(x_t, v ) - \cL_g(x_t, v_t ) \big )   \right ] 
 \leq   \EE\Big [  \frac{5\|  ( v - v_t) \lambda \|^2}{ 2 \eta  } +  \frac{ \eta \| x_t - x_{t-1} \|^2}{10} \Big ]
\\
& \quad\leq \frac{5 M_{\lambda}^2 \EE [\|  v - v_t\|^2]}{ 2 \eta_{t-1}  } +  \frac{ \eta  \EE [ \| x_t - x_{t-1} \|^2 ] }{10}
 \leq    \frac{10 C_g^2 M_{\lambda}^2 }{ \eta   } +  \frac{ \eta \EE[ \| x_t - x_{t-1} \|^2] }{10},
\end{split}
\end{equation*}
where the first inequality uses the fact that $ab \leq \frac{5a^2}{2} + \frac{b^2}{10}$, the second inequality holds by the condition that $\| \lambda\| \leq M_{\lambda}$, and the third inequality holds by Assumption \ref{assumption:02} that $\| v -v_t\|^2 \leq 2\| v\|^2 + 2\| v_t\|^2 \leq 4C_g^2$. 
Summing the above inequality over $t=1,\cdots, K$, we obtain 
\begin{equation*}
\begin{split}
 \EE \Big [\sum_{t=1}^K \lambda^\top \big (  \cL_g(x_t, v ) - \cL_g(x_t, v_t ) \big ) \Big ] \leq  \frac{10 C_g^2 M_{\lambda}^2 }{ \eta } +  \frac{  \eta  \EE [ \| x_t - x_{t-1} \|^2 ] }{10 },
\end{split}
\end{equation*}
which completes the proof. 
\QED 

\subsection{Lemma \ref{lemma:lambda} and Its Proof}\label{app:A4}
\begin{lemma}\label{lemma:lambda}
Suppose Assumptions \ref{assumption:01}, \ref{assumption:02}, and \ref{assumption:03} hold, and Algorithm \ref{alg:01} generates  $\{ (x_t,\lambda_t,\pi_{1:2,t}, v_t) \}_{t=1}^N$ by setting $\tau_t = t/2$, $\alpha_{t} = \alpha $, and $\eta_t =  \eta $ for all $0 \leq t \leq N$.  Let $\lambda \in \RR_+^m$ be a bounded  nonnegative random variable such that $ \| \lambda \| \leq M_{\lambda}$.
For any integer $K \leq N$, we have 
\begin{equation*}
\begin{split}
& \EE \Big [ \sum_{t=1}^K (\lambda  - \lambda_{t})^\top \cL_g(x_{t}, v_{t} )\Big ]  +  \frac{\alpha_{K-1}}{2}\EE[ \| \lambda_K - \lambda \|^2]\\
 &\quad \leq   \sqrt{K} M_{\lambda}  \sigma_g
+\sum_{t=1}^K \Big ( \frac{5 C_g^2  \EE[  \|    \lambda  - \lambda_{t} \|^2 ] }{2\eta }  +  \frac{\eta }{10} \EE[\| x_{t-1} - x_t\|^2]  -   \frac{\alpha }{4} \EE[ \| \lambda_{t-1} - \lambda_{t} \|^2 ]  \Big ) \\
& \qquad +   \frac{\alpha }{2} \EE [  \| \lambda_{0}- \lambda  \|^2   ]
+ \frac{K \sigma_{g}^2}{ \alpha }.
 \end{split}
\end{equation*}
\end{lemma}
\textit{Proof:}
By decomposing $(\lambda  - \lambda_{t})^\top \cL_g(x_{t}, v_{t} )$, we have, 
\begin{equation}\label{eq:lambda_1}
\begin{split}
&  (\lambda  - \lambda_{t}) ^\top \cL_g(x_{t}, v_{t} ) \\
& \quad =   (\lambda  - \lambda_{t})^\top \Big (  \cL_g(x_{t}, v_{t}) - \cL_g(x_{t-1}, v_{t}) \Big )  + (\lambda  - \lambda_{t})^\top  [  \cL_g(x_{t-1}, v_{t} )  - g(x_{t-1},\zeta_{t-1}^1 )  ] 
\\
& \qquad + (\lambda  - \lambda_{t}) g(x_{t-1},\zeta_{t-1}^1 ) \\
& =  \underbrace{ \inner{   v_{t} (\lambda  - \lambda_{t}) } { x_{t-1} - x_t  } }_{T_{1,t}} + \underbrace{(\lambda  - \lambda_{t})^\top [  g(x_{t-1})  - g(x_{t-1},\zeta_{t-1}^1 )  ] }_{T_{2,t}} 
\\
& \qquad + \underbrace{(\lambda  - \lambda_{t})^\top  g(x_{t-1},\zeta_{t-1}^1 )}_{T_{3,t}}.
\end{split}
\end{equation}
We provide bounds for $T_{1,t}, T_{2,t}$, and $T_{3,t}$ in our analysis. 

First,  for $T_{1,t}$, by using the fact that $\inner{a}{b} \leq \frac{5\| a\|^2}{2} + \frac{\| b\|^2}{10}$ and Assumption~\ref{assumption:02} that $ \| v_t \|^2 \leq C_g^2$, we obtain 
\begin{equation}\label{eq:lambda_t1}
\begin{split}
\EE [ \inner{   v_{t} (\lambda  - \lambda_{t}) } { x_{t-1} - x_t  }  ] & \leq  \frac{5\EE[  \| v_{t} (\lambda  - \lambda_{t})  \|^2 ] }{2\eta }+ \frac{\eta  \EE [ \| x_{t-1} - x_t\|^2 ] }{10}
\\
& \leq  \frac{5 C_g^2 \EE[  \| \lambda  - \lambda_{t}  \|^2 ] }{2\eta }+ \frac{\eta  \EE [ \| x_{t-1} - x_t\|^2 ] }{10}.
\end{split}
\end{equation}
Second, consider $T_{2,t}$, by  the independence between $ \lambda_{t-1}$ and the fact that  $\EE [ g(x_{t-1})  - g(x_{t-1},\zeta_{t-1}^1 ) ] = 0$, we have 
$$\EE\big [ \lambda_{t-1}^\top ( g(x_{t-1})    - g(x_{t-1},\zeta_{t-1}^1 ))  \big] =0,$$
which further implies that 
\begin{equation}\label{eq:single_t2_00}
\begin{split}
 \EE[T_{2,t}] 
& =  \EE \Big [ (\lambda  - \lambda_{t-1})^\top [ g(x_{t-1})   - g(x_{t-1},\zeta_{t-1}^1 ) ] \Big]  
\\
& \quad + \EE \Big [ (\lambda_{t-1}  - \lambda_{t})^\top [ g(x_{t-1})   - g(x_{t-1},\zeta_{t-1}^1 ) ] \Big]  \\
& =  \EE \Big [ \lambda^\top [ g(x_{t-1})  - g(x_{t-1},\zeta_{t-1}^1 ) ] \Big]  + \EE \Big [ (\lambda_{t-1}  - \lambda_{t})^\top [ g(x_{t-1})   - g(x_{t-1},\zeta_{t-1}^1 )  ] \Big] . 
\end{split}
\end{equation}
Recall that  $\lambda $ is random but bounded such that $\| \lambda\| \leq M_{\lambda}$,
by  Lemma~\ref{lemma:sum_martingale} and Assumption \ref{assumption:02} that $\EE[\|  g(x_{t-1})   - g(x_{t-1},\zeta_{t-1}^1 ) \|^2] \leq \sigma_g^2$, we obtain 
\begin{equation}\label{eq:single_t2_01}
\EE \Big [ \sum_{t=1}^K \lambda^\top  \Big ( g(x_{t-1})    - g(x_{t-1},\zeta_{t-1}^1 )  \Big ) \Big ] \leq \sqrt{K} M_{\lambda}  \sigma_g.
\end{equation}
Meanwhile, by the fact that $\inner { a} { b}\leq \frac{\alpha \| a\|^2}{4}+ \frac{\| b\|^2}{\alpha}$ for $\alpha >0$, we bound the second term within the right side of~\eqref{eq:single_t2_00} by 
\begin{equation}\label{eq:single_t2_02}
\begin{split}
 & \EE \Big [ (\lambda_{t-1}  - \lambda_{t})^\top [  g(x_{t-1})   - g(x_{t-1},\zeta_{t-1}^1 )  ] \Big]  
 \\
 & 
\leq  \frac{ \alpha  \EE[ \| \lambda_{t-1}  - \lambda_{t}\|^2 ] }{4} + \frac{\EE[ \| g(x_{t-1})   - g(x_{t-1},\zeta_{t-1}^1 )\|^2 ] }{\alpha }
 \leq \frac{ \alpha  \EE[ \| \lambda_{t-1}  - \lambda_{t}\|^2 ] }{4} + \frac{\sigma_g^2}{\alpha }.
\end{split}
\end{equation}
Summing \eqref{eq:single_t2_00} over $t=1,\cdots, K$ and applying \eqref{eq:single_t2_01} and \eqref{eq:single_t2_02}, we obtain 
\begin{equation}\label{eq:lambda_t2}
\sum_{t=1}^K \EE[T_{2,t} ] \leq \sqrt{K} M_{\lambda}  \sigma_g + \sum_{t=1}^K    \Big ( \frac{ \alpha  \EE [ \| \lambda_{t-1}  - \lambda_{t}\|^2 ] }{4} + \frac{\sigma_{g}^2}{\alpha }  \Big ). 
\end{equation}
Third, consider $T_{3,t}$, by recalling \eqref{eq:lambda_ssd_update} that 
\begin{equation*}
\begin{split}
\lambda_{t} \in  \argmax_{ \lambda \in \RR_+^{m}} \Big \{ \inner{  g(x_{t-1},\zeta_{t-1}^1 )   }{ \lambda }- \frac{\alpha }{2} \| \lambda_t - \lambda \|^2 \Big \},
\end{split}
\end{equation*}
and by Lemma \ref{lemma:three_point}, we have 
\begin{equation}\label{eq:lambda_t3}
\begin{split}
T_{3,t} & =  ( \lambda - \lambda_{t} )^\top g(x_{t-1},\zeta_{t-1}^1 )
=  - (\lambda_{t} - \lambda  )^\top g(x_{t-1},\zeta_{t-1}^1 )  \\
 & \leq  \frac{\alpha }{2} \| \lambda_{t-1}- \lambda  \|^2  - \frac{\alpha }{2} \| \lambda_{t-1} - \lambda_{t} \|^2 -  \frac{\alpha }{2} \| \lambda_{t} - \lambda  \|^2.
\end{split}
\end{equation} 
Finally, summing \eqref{eq:lambda_1} over $t=1,\cdots, N$, taking expectations on both sides, and plugging in \eqref{eq:lambda_t1}, \eqref{eq:lambda_t2}, and \eqref{eq:lambda_t3}, we conclude that
\begin{equation*}
\begin{split}
& \EE \Big [ \sum_{t=1}^K (\lambda  - \lambda_{t})^\top \cL_g(x_{t}, v_{t} )\Big ]  +  \frac{\alpha }{2}\EE[ \| \lambda_K - \lambda \|^2]\\
 &\quad \leq   \sqrt{K} M_{\lambda}  \sigma_g +\sum_{t=1}^K \Big ( \frac{5 C_g^2 \EE[  \|     \lambda  - \lambda_{t}  \|^2 ] }{2\eta }  +  \frac{\eta }{10} \EE[\| x_{t-1} - x_t\|^2]  -   \frac{\alpha }{4} \EE[ \| \lambda_{t-1} - \lambda_{t} \|^2 ]  \Big ) \\
 & \qquad +   \frac{\alpha }{2} \EE [  \| \lambda_{0}- \lambda  \|^2   ]
+ \frac{K \sigma_{g}^2}{ \alpha },
 \end{split}
\end{equation*}
which completes the proof. 
\QED 

\subsection{Proof of Proposition~\ref{prop:bounded_lambda_01}}
By applying Theorem \ref{thm:01} with  $\lambda = \lambda^*$, $\pi_{1:2} = \pi_{1:2}^*$, and $v = v^*$ defined in \eqref{eq:saddle_point}, we set $M_{\lambda^*} = \| \lambda^*\|$ and obtain 
\begin{equation*}
\begin{split}
& \sum_{t=1}^K \EE \Big ( \cL(x_t,\lambda^*,  \pi_{1:2}^* , v^* ) -  \cL(x^*, \lambda_t,  \pi_{1:2,t} , v_t)   \Big )  +  \frac{\alpha }{2}\EE[ \| \lambda_K - \lambda^*\|^2] 
\\
 & \quad \leq     \frac{95K C_{f_1}^2 C_{f_2}^2 }{2\eta } 
+  \sqrt{ K } L_{f_1} \sigma_{f_2}^2  +   3 \sqrt{K} C_{f_1}\sigma_{f_2}  + \frac{\eta }{2} \| x_{0}- x^*\|^2  \\
& \qquad +     \frac{K D_X^2 C_g^2 }{\alpha }
 +  \frac{10K C_g^2 \| \lambda^* \|^2 } {  \eta }  
 +   \sum_{t=1}^K \frac{5  C_g^2   \EE[  \|  \lambda_{t-1}\|^2]  }{2\eta  }  + \sqrt{K} \|\lambda^*\| \sigma_g  
 \\
 & \qquad +   \frac{\alpha }{2}   \| \lambda_{0} - \lambda^*  \|^2   
+ \frac{K \sigma_{g}^2}{ \alpha }
+\sum_{t=1}^K \frac{5  C_g^2   \EE[  \| \lambda^*  - \lambda_{t} \|^2]   }{2\eta }  .
\end{split}
\end{equation*}
By using the fact that $\| \lambda_t  \|^2 \leq 2 \| \lambda_t  - \lambda^* \|^2 + 2 \|  \lambda^* \|^2  $, we further express the above inequality as 
\begin{equation*}
\begin{split}
& \sum_{t=1}^K \EE \Big ( \cL(x_t,\lambda^*,  \pi_{1:2}^* , v^* ) -  \cL(x^*, \lambda_t,  \pi_{1:2,t} , v_t)   \Big )  +  \frac{\alpha }{2}\EE[ \| \lambda_K - \lambda^*\|^2] 
\\
 & \quad\leq     \frac{95 K C_{f_1}^2 C_{f_2}^2 }{2\eta } 
+  \sqrt{ K } L_{f_1} \sigma_{f_2}^2  +   3 \sqrt{K} C_{f_1}\sigma_{f_2}  + \frac{\eta }{2} \| x_{0}- x^*\|^2  \\
& \qquad +     \frac{KD_X^2 C_g^2}{\alpha }    + \sum_{t=1}^K \frac{15 C_g^2   \| \lambda^* \|^2}{  \eta  }     
 +   \sum_{t=1}^K \frac{5 C_g^2\EE[  \|  \lambda_{t-1} -  \lambda^*\|^2] }{\eta   }  
 + \sqrt{K} \|\lambda^*\| \sigma_g 
 \\
 & \qquad +   \frac{\alpha }{2}   \| \lambda_{0} - \lambda^*  \|^2   
+ \frac{K \sigma_{g}^2}{ \alpha }
+\sum_{t=1}^K \frac{5 C_g^2 \EE[  \| \lambda^*  - \lambda_{t} \|^2  ]   }{2\eta }.
\end{split}
\end{equation*}
By using the min-max relationship \eqref{eq:saddle_point} that 
$$\cL(x_t,\lambda^*,  \pi_{1:2}^* , v^*) -  \cL(x^*, \lambda_t,  \pi_{1:2,t} , v_{1,t})   \geq 0,$$ dividing  both sides of the above inequality by $\frac{\alpha }{2}$, and setting $\alpha  = 2\sqrt{N} $ and $\eta =  \frac{15 C_g^2  \sqrt{N}}{2}$, we obtain 
\begin{equation*}
\begin{split}
\EE [   \| \lambda_{K} - \lambda^*  \|^2 ] & \leq   \frac{2}{\alpha}\sum_{t=1}^K \EE \Big ( \cL(x_t,\lambda^*,  \pi_{1:2}^* , v^* ) -  \cL(x^*, \lambda_t,  \pi_{1:2,t} , v_t)   \Big )  +  \EE[ \| \lambda_K - \lambda^*\|^2]
\\
& \leq R_K + \frac{  \EE[  \| \lambda^*  - \lambda_{K}\|^2 ] }{3N} +\sum_{t=0}^{K-1} \frac{  \EE[  \| \lambda^*  - \lambda_{t}\|^2 ] }{N} , 
\end{split}
\end{equation*}
where 
\begin{equation*}
\begin{split}
R_K & =  \frac{K}{N} \Big (\frac{  7 C_{f_1}^2 C_{f_2}^2  }{C_g^2}+   \frac{D_X^2 C_g^2}{2}  + 2  \|  \lambda^* \|^2  +  \frac{ \sigma_{g}^2}{2} \Big ) + \frac{15 C_g^2 }{4} \| x_{0}- x^*\|^2 +  \| \lambda_{0} - \lambda^*  \|^2  
\\
&  \quad  +  \frac{\sqrt{K}}{\sqrt{N}}  \Big (  
 L_{f_1} \sigma_{f_2}^2  +   3 C_{f_1}\sigma_{f_2}  +   \| \lambda^*\|\sigma_g \Big ) . 
\end{split}
\end{equation*}
Next, letting $R$ be the constant defined in \eqref{eq:R_1}, we note that $R_K \leq R$, and 
the above inequality further implies that 
\begin{equation*}
\begin{split}
\Big (1-\frac{1 }{3N} \Big ) \EE [   \| \lambda_{K} - \lambda^*  \|^2 ] \leq R 
+\sum_{t=0}^{K-1} \frac{\EE[  \| \lambda^*  - \lambda_{t}\|^2 ] }{N} .
\end{split}
\end{equation*}
Dividing  both sides of the above inequality by $(1-\frac{1}{3N})$, using the fact that $(1-\frac{1 }{3N}) \geq 1/2$ for $N \geq 2$, and applying Lemma \ref{lemma:recursive_bound}, we conclude that
 \begin{equation*}
\begin{split}
 \EE [   \| \lambda_{K} - \lambda^*  \|^2 ] & \leq \Big  (1-\frac{1 }{3N} \Big )^{-1} \Big ( R 
+\sum_{t=0}^{K-1} \frac{ \EE[  \| \lambda^*  - \lambda_{t}\|^2 ] }{N}  \Big ) 
\\
& \leq 2 R + 2  \sum_{t=0}^{K-1} \frac{\EE[  \| \lambda^*  - \lambda_{t}\|^2 ] }{N} 
\leq 2R \Big (1 + \frac{2}{N} \Big )^K \leq 2R e^{2},
\end{split}
\end{equation*}
 for all $1\leq K \leq N$. 
This completes the proof. 
\QED

\subsection{Proof of Theorem~\ref{thm:rate}}
\textit{Proof:} In our analysis, we set $\lambda_0 = 0$ for ease of presentation. 
Recall that $\bar x_N = \frac{1}{N} \sum_{t=1}^N x_t$, for any fixed $(\lambda, \pi_{1:2}, v)$, by  the convexity of $\cL(x,\lambda,  \pi_{1:2} , v) $ with respect to $x$, we have 
\begin{equation}\label{eq:f_xN_convexity}
 \frac{1}{N} \sum_{t=1}^N \cL(x_t,\lambda,  \pi_{1:2} , v) \geq \cL(\bar x_N, \lambda,  \pi_{1:2} , v).
 \end{equation}
Meanwhile, we denote by 
\begin{equation}\label{eq:hat_pi_single}
\hat \pi_2 \in  \partial  f_2(\bar x_N) , \, \hat \pi_1=  \nabla f_1\big ( f_2( \bar x_N)\big ), \mbox{ and } \, \hat v \in \partial g(\bar x_N),
\end{equation}
the dual variables associated with $\bar x_N$. By  the definition of composite Lagrangian \eqref{def:Lagrangian_single}, we have 
\begin{equation}\label{eq:f_xN}
F(\bar x_N) = \cL_{F}(\bar x_N, \hat \pi_{1:2}) = \cL(\bar x_N, 0, \hat \pi_{1:2}, \hat v  ).
\end{equation}

First, we derive the convergence rate of the objective optimality gap $F(\bar x_N) - F(x^*) $. Let $\bar \lambda_N = \frac{1}{N} \sum_{t=1 }^N \lambda_t$, $\bar \pi_{1:2,N} = \frac{1}{N}\sum_{t=1}^N \pi_{1:2,t}$, and $\bar v_N = \frac{1}{N} \sum_{t=1}^N  v_t $. By setting $\lambda = 0 $, $\pi_{1:2} = \hat \pi_{1:2}$, and $v = \hat v$ within \eqref{eq:hat_pi_single}, and using \eqref{eq:saddle_point} that $\cL(x^*,\lambda^*,\pi_{1:2}^*,v^*)  \geq \cL(x^*,\lambda,\pi_{1:2},v)  $
for any $(\lambda,\pi_{1:2},v) \in \RR_+^m \times \Pi \times \cV$,  we have 
\begin{equation}\label{eq:prop1_01}
\begin{split}
 & \frac{1}{N} \sum_{t=1}^N \Big ( \cL(x_t, 0, \hat \pi_{1:2}, \hat v  ) - \cL(x^*, \lambda_t,\pi_{1:2,t}, v_t  ) \Big ) 
\\
 & \quad \quad  \geq  \cL(\bar x_N, 0, \hat \pi_{1:2}, \hat v  ) - \cL(x^*, \lambda^*,\pi_{1:2}^*, v^*  ) 
= F(\bar x_N) - F(x^*),
\end{split}
\end{equation}
where the last equality holds by \eqref{eq:f_xN}. 
Therefore, by  \eqref{eq:prop1_01}, setting the feasible point as $(x^*, \lambda,  \pi_{1:2},  v ) = (x^*, 0, \hat \pi_{1:2}, \hat v )$, and applying Theorem \ref{thm:01} with the fact $\|\lambda\| = 0$ when $\lambda = 0$,
we have 
\begin{equation*}
\begin{split}
& \EE [ F(\bar x_N) ]- F(x^*) \leq   \frac{1}{N} \sum_{t=1}^N \EE \Big ( \cL(x_t, 0, \hat \pi_{1:2}, \hat v  ) - \cL(x^*, \lambda_t,\pi_{1:2,t}, v_t  ) \Big )  \\
& \quad \ \leq   \frac{1}{ \sqrt{N}} \Big ( \frac{ 7 C_{f_1}^2 C_{f_2}^2  }{C_g^2}
+  L_{f_1} \sigma_{f_2}^2  +   3C_{f_1}\sigma_{f_2}  + \frac{15 C_g^2 }{4} \| x_{0}- x^*\|^2 
   \Big )
  \\
 & \qquad  +  \frac{1}{N} \sum_{t=1}^N \Big ( \frac{ \EE[  \|  \lambda_{t-1}\|^2]  }{3 \sqrt{N} }  +  \frac{\EE[  \| \lambda_{t}\|^2 ] }{3\sqrt{N}}  \Big ) 
\\
& \quad \ \leq   \frac{1}{ \sqrt{N}} \Big ( \frac{ 7 C_{f_1}^2 C_{f_2}^2  }{C_g^2 }
+ L_{f_1} \sigma_{f_2}^2  +   3C_{f_1}\sigma_{f_2}  + \frac{15 C_g^2 }{4} \| x_{0}- x^*\|^2 
\Big ) +  \frac{(4\| \lambda^*\|^2  + 8Re^2 )  }{3 \sqrt{N}} ,
\end{split}
\end{equation*}
where the last inequality comes from the bound of $\| \lambda_t\|^2$ provided by Proposition \ref{prop:bounded_lambda_01} that 
\begin{equation}\label{eq:lambda_t_radius_01}
\EE[\| \lambda_t\|^2] \leq 2 \| \lambda^*\|^2 + 2 \EE[ \| \lambda_t - \lambda^*\|^2] \leq   2\| \lambda^*\|^2  + 4Re^2, \text{ for all }t=1,2,\cdots, N.
\end{equation}
This establishes the convergence rate for the objective optimality gap $F(\bar x_N) - F(x^*) $. 

Second, we consider the feasibility residual  $\| g(\bar x_N)_{+} \|_2$.  
Recall that $\hat \pi_{1:2}$ and $\hat v $ are the dual variables associated with $\bar x_N$ such that  $\hat \pi_{1:2} \in \argmax_{\pi_{1,2} \in \Pi} \cL_F(\bar x_N,\pi_{1:2})$ and  $\hat v \in \argmax_{v \in \cV} \cL_g(\bar x_N,v)$, we have 
$$
F(\bar x_N) = \cL_F (\bar x_N, \hat \pi_{1:2}) \geq \cL_F(\bar x_N, \pi_{1:2}^*) 
\text{ and }
g(\bar x_N) = \cL_g(\bar x_N, \hat v) \geq \cL_g(\bar x_N, v^*).
$$
Since $\lambda^* \geq 0$, by  the above inequalities and the min-max relationship \eqref{eq:saddle_point}, we can see
\begin{equation*}
\begin{split}
 & \cL(\bar x_N, \lambda^*, \hat \pi_{1:2}, \hat v)  =  \cL_F(\bar x_N, \hat \pi) + \inner{ \lambda^* }{  \cL_g(\bar x_N, \hat v)  } \\
&  \quad\geq \cL_F(\bar x_N, \pi_{1:2}^*) + \inner{ \lambda^* }{  \cL_g(\bar x_N,  v^*) } =  \cL(\bar x_N, \lambda^*, \pi_{1:2}^*, v^*)  \geq \cL(x^*, \lambda^*,\pi_{1:2}^*, v^*),
\end{split}
\end{equation*}
which further implies 
$$
F(\bar x_N) + \inner{  \lambda^* }{  g(\bar x_N)  } -F(x^*) \geq 0.
$$
Meanwhile, due to the facts that $\lambda^* \geq 0$ and $g(\bar  x_N) \leq g(\bar  x_N)_{+}$, we have  $ \langle \lambda^*, g(\bar  x_N) \rangle \leq  \langle \lambda^*, g(\bar x_N)_{+} \rangle$, hence, 
\begin{equation}\label{eq:lambda_pos}
F(\bar x_N) + \| \lambda^* \|  \| g(\bar x_N)_{+} \|   -F(x^*) \geq  F(\bar x_N) + \langle \lambda^*, g(\bar  x_N) \rangle -F(x^*) \geq 0.
\end{equation}
Let $ \tilde \lambda = (\| \lambda^*\|_2  +1 ) g(\bar x_N)_{+} /  \| g(\bar x_N)_{+} \| $. Consider another feasible point $ (x^*, \tilde \lambda, \hat \pi_{1:2}, \hat v ) $, by  \eqref{eq:f_xN_convexity} and the facts that $g(\bar x_N)^\top g(\bar x_N)_+ = \| g(\bar x_N)_+\|^2 $ and $ \cL(x^*, \lambda_t, \pi_{1:2,t}, v_t)
\leq  \cL(x^*, \lambda^*, \pi_{1:2}^*, v^*)$, we have $\langle g(\bar x_N),  \tilde \lambda  \rangle  = (\| \lambda^*\|_2  +1 ) \| g(\bar x_N)_{+} \|  $, which further yields that 
\begin{equation*}
\begin{split}
& \frac{1}{N} \sum_{t=1}^N \Big ( \cL(x_t,  \tilde \lambda, \hat \pi_{1:2}, \hat v) - \cL(x^*, \lambda_t, \pi_{1:2,t}, v_t)  \Big ) \geq  \cL(\bar x_N,  \tilde \lambda, \hat \pi_{1:2}, \hat v) - \cL(x^*, \lambda^*, \pi_{1:2}^*, v^*) \\
& \quad = F(\bar x_N) + \langle \tilde \lambda,  g(\bar x_N)  \rangle -F(x^*) = F(\bar x_N) +    (\| \lambda^*\|_2  +1 ) \| g(\bar x_N)_{+} \|   -F(x^*).
\end{split}
\end{equation*}
By rearranging the terms in the  inequality above and applying \eqref{eq:lambda_pos}, we obtain 
\begin{equation*}\label{eq:prop1_02}
\begin{split}
 \| g(\bar x_N)_{+} \| 
 & \leq   \frac{1}{N} \sum_{t=1}^N \Big ( \cL(x_t,  \tilde \lambda, \hat \pi_{1:2}, \hat v) - (  \cL(x^*, \lambda_t, \pi_{1:2,t}, v_t)  \Big ) 
 \\
 & \quad -  \Big   ( F(\bar x_N) + \| \lambda^* \|  \| g(\bar x_N)_{+} \|  - F(x^*)  \Big ) \\
& \leq    \frac{1}{N} \sum_{t=1}^N \Big ( \cL(x_t,  \tilde \lambda, \hat \pi_{1:2}, \hat v) - \cL(x^*, \lambda_t, \pi_{1:2,t}, v_t)  \Big )   .
\end{split}
\end{equation*}
Next, we note that $ \| \tilde \lambda \| =    \| \lambda^*\|  +1$.
By  the  inequality above, considering the feasible point $(x^*, \tilde \lambda, \hat \pi_{1:2}, \hat v )$, and applying Theorem \ref{thm:01} with $\lambda_0  = 0$, $M_{\tilde \lambda} =  \| \tilde \lambda\| = \| \lambda^*\| + 1$, $\alpha = 2 \sqrt{N}$,  and $\eta = \frac{15 C_g^2  \sqrt{N}}{2}$,  we have 
\begin{equation*}
\begin{split}
&  \EE[ \| g(\bar x_N)_{+} \|  ] \\
 & \leq    \frac{95  C_{f_1}^2 C_{f_2}^2 }{ 2\eta } 
+  \frac{ L_{f_1} \sigma_{f_2}^2 }{\sqrt{ N }  } +   \frac{ 3  C_{f_1}\sigma_{f_2}}{\sqrt{N}}  + \frac{\eta }{2N} \| x_{0}- x^*\|^2 
+     \frac{D_X^2 C_g^2}{\alpha }   + \frac{ \| \tilde \lambda \|  \sigma_g }{\sqrt{N}}   \\ 
& \quad 
+   \frac{1}{N} \Big[ \frac{\alpha }{2}    \|  \tilde \lambda  \|^2   
+ \frac{N \sigma_{g}^2}{ \alpha }
 +   \sum_{t=1}^N \Big ( \frac{5 C_g^2 \EE[  \|  \lambda_{t-1}\|^2]  }{2\eta   }  +  \frac{10 C_g^2 \EE [  \| \tilde  \lambda\|^2]}{ \eta  } + \frac{5 C_g^2 \EE[  \| \tilde  \lambda  - \lambda_{t} \|^2 ] }{2\eta }  \Big ) \Big ] 
 \\
  &\leq  \frac{1}{\sqrt{N} }  \Big (  \frac{ 7 C_{f_1}^2 C_{f_2}^2  }{C_g^2}
+  L_{f_1} \sigma_{f_2}^2  +   3  C_{f_1}\sigma_{f_2}  + \frac{15 }{4} \| x_{0}- x^*\|^2 
+     \frac{D_X^2 C_g^2}{2}    + \| \tilde \lambda\|   \sigma_g  +  \| \tilde \lambda\| ^2 
 \Big ) \\ 
& \quad + \frac{ \sigma_{g}^2}{ 2\sqrt{N}}  
+   \frac{1}{N} \sum_{t=1}^N \Big ( \frac{  \EE[  \|  \lambda_{t-1}\|^2]  }{3\sqrt{N}  }  +  \frac{4   \| \tilde \lambda\| ^2  }{ 3 \sqrt{N} } + \frac{2 (  \| \tilde \lambda\| ^2   + \EE[ \| \lambda_t\|^2] ) }{3\sqrt{N}}  \Big )
 \\
   &  \leq  \frac{1}{\sqrt{N} }  \Big ( \frac{ 7 C_{f_1}^2 C_{f_2}^2 }{C_g^2}
+  L_{f_1} \sigma_{f_2}^2  +   3  C_{f_1}\sigma_{f_2}  + \frac{15  C_g^2 }{4} \| x_{0}- x^*\|^2 
+     \frac{D_X^2 C_g^2}{2}  
+ \frac{ \sigma_{g}^2}{ 2}  \Big ) \\ 
& \quad 
 +    \frac{1}{\sqrt{N} }  \Big (   (2\| \lambda^*\|^2  + 4Re^2 )   +   3 \| \tilde \lambda \|^2   +  \| \tilde \lambda\|   \sigma_g   \Big ),
\end{split}
\end{equation*}
where the second inequality holds because $\| \tilde \lambda - \lambda_t\|^2 \leq 2\| \tilde \lambda\|^2 + 2 \| \lambda_t\|^2$, and the last inequality holds by the bound of $\| \lambda_t\|^2$ in \eqref{eq:lambda_t_radius_01}.
The desired result can be acquired by substituting $\| \tilde \lambda\| =  \| \lambda^*\|_2  +1$ and $ \| \tilde \lambda\|^2 = ( \| \lambda^*\|_2  +1)^2$ into the above inequality. This completes the proof. 
\QED 

\section{Proof of Results in Section \ref{sec:composition}}\label{app:B}
\subsection{Proof of Lemma~\ref{lemma:comp_lambda_G}}
First, we decompose $(\lambda  - \lambda_{t})^\top  \cL_G(x_t, v_{1:2,t} ) $ that 
\begin{equation}\label{eq:comp_t0}
\begin{split}
& (\lambda  - \lambda_{t})^\top  \cL_G(x_t, v_{1:2,t} ) \\
& \quad =  (\lambda  - \lambda_{t})^\top  \Big (  \cL_G(x_t, v_{1:2,t} )- \cL_G(x_{t-1}, v_{1:2,  t} ) \Big )  + (\lambda  - \lambda_{t})^\top   \cL_G(x_{t-1}, v_{1:2,t}  )  
\\
& \quad =   \underbrace{ \inner{   v_{2,t} v_{1,t}(\lambda  - \lambda_{t}) }{ x_{t} - x_{t-1} }  }_{T_{1,t}  } + \underbrace{(\lambda  - \lambda_{t})^\top [  \cL_G(x_{t-1}, v_{1:2,t}  )  - H_t ] }_{T_{2,t}}+ \underbrace{(\lambda  - \lambda_{t})^\top   H_t  }_{T_{3,t}  }.
\end{split}
\end{equation}
We now provide bounds for the three terms $T_{1,t}, T_{2,t} $, and $T_{3,t} $. 
First, for $T_{1,t}$, by  the fact that $\inner{a}{b} \leq \frac{5\| a\|^2}{2\eta} + \frac{\eta \| b \|^2}{10}$ for any vectors $a,b$, we have 
\begin{equation}\label{eq:comp_t1}
\begin{split}
\sum_{t=1}^K \EE [ T_{1,t}] & =   \sum_{t=1}^K \EE [  (x_{t} - x_{t-1})^\top  v_{2, t} v_{1, t}(\lambda  - \lambda_{t})  ] \\
 & \leq   \sum_{t=1}^K \frac{5\EE[  \|  v_{2,t}   v_{1,t}   (\lambda  - \lambda_{t}) \|^2 ] }{2\eta }+  \sum_{t=1}^K \frac{\eta \EE [ \| x_{t} - x_{t-1}\|^2 ] }{10}
 \\
  & \leq   \sum_{t=1}^K \frac{5 C_{g_1}^2 C_{g_2}^2 \EE[  \|  \lambda  - \lambda_{t} \|^2 ] }{2\eta }+  \sum_{t=1}^K \frac{\eta \EE [ \| x_{t} - x_{t-1}\|^2 ] }{10},
\end{split}
\end{equation}
where the last inequality holds by Assumption~\ref{assumption:compositional_g} (d)-(e) that $\| v_{1,t}\|^2 \leq C_{g_1}^2$ and $\| v_{2,t}\|^2 \leq C_{g_2}^2$. 
Second, we consider $T_{2,t} $ and denote by $ \tilde \Delta_{2,t} = \cL_G(x_{t-1}, v_{1:2 ,t } )  - H_t $. Given $x_{t-1}$, we observe from Algorithm~\ref{alg:02} that 
$H_t $ is conditionally independent of $\lambda_{t-1}$. 
Together with Lemma \ref{lemma:H} (a) that $H_t$ is also an unbiased estimator for $ \cL_G(x_{t-1}, v_{1:2 , t}  ) $, we have $\EE[\lambda_{t-1}^\top \tilde \Delta_{2,t}] =0$, which further implies that 
\begin{equation*}
\begin{split}
\EE[T_{2,t}  ] & = \EE[  \lambda^\top \tilde  \Delta_{2,t}   ]  -  \EE[   \tilde  \Delta_{2,t}^\top \lambda_{t-1}]  +  \EE[  (\lambda_{t-1} - \lambda_t)^\top  \tilde   \Delta_{2,t} ]
\\
& =  \EE[  \lambda^\top  \tilde  \Delta_{2,t}   ]   +  \EE[  (\lambda_{t-1} - \lambda_t)^\top  \tilde   \Delta_{2,t} ] .
\end{split}
\end{equation*}
By setting $M_{\lambda}  = \| \lambda \|$
and recalling Lemma \ref{lemma:H} (b) that $\text{Var}(\tilde \Delta_{2,t}) \leq \sigma_H^2$,  we further have 
\begin{equation}\label{eq:comp_t2}
\begin{split}
\sum_{t=1}^K \EE[T_{2,t}  ] & = \sum_{t=1}^K   \EE \Big [ \lambda^\top  \tilde  \Delta_{2,t}    \Big]  + \sum_{t=1}^K  \EE \Big [  (\lambda_{t-1} - \lambda_t)^\top  \tilde  \Delta_{2,t}  \Big]   \\
& \leq   \sqrt{K} M_{\lambda}  \sigma_H  + \sum_{t=1}^K  \EE \Big [  (\lambda_{t-1} - \lambda_t)^\top \tilde  \Delta_{2,t}   \Big]  \\
& \leq  \sqrt{K} M_{\lambda}  \sigma_H + \sum_{t=1}^K  \Big ( \frac{ \alpha   \EE [ \| \lambda_{t-1}  - \lambda_{t}\|^2]}{4} + \frac{\sigma_{H}^2}{\alpha }  \Big ),
\end{split}
\end{equation}
where the first inequality holds by Lemma \ref{lemma:sum_martingale} in Appendix Section \ref{app:A}, and the last inequality holds by the fact $\inner{a}{b}\leq \frac{\| a \|^2}{4} + \| b \|^2$ for any vectors $a,b$. 

Finally, for $T_{3,t}$, by applying the three-point Lemma \ref{lemma:three_point} to the update rule that 
\begin{equation*}
\begin{split}
\lambda_{t} \in  \argmax_{ \lambda \in \RR_+^m } \Big \{ \lambda^\top H_t - \frac{\alpha }{2} \| \lambda_{t-1} - \lambda \|^2 \Big \},
\end{split}
\end{equation*}
we have that for any $\lambda \in \RR_+^m$, 
\begin{equation*}
\begin{split}
 ( \lambda - \lambda_{t} )^\top H_t
=  - (\lambda_{t} - \lambda  )^\top H_t   \leq  \frac{\alpha }{2} \| \lambda_{t-1}- \lambda  \|^2  - \frac{\alpha }{2} \| \lambda_{t-1} - \lambda_{t} \|^2 -  \frac{\alpha }{2} \| \lambda_{t} - \lambda  \|^2,
\end{split}
\end{equation*} 
which implies 
\begin{equation}\label{eq:comp_t3}
\EE[T_{3,t}  ]  \leq  \sum_{t=1}^K \Big ( \frac{\alpha }{2} \| \lambda_{t-1}- \lambda  \|^2  - \frac{\alpha }{2} \| \lambda_{t-1} - \lambda_{t} \|^2 -  \frac{\alpha }{2} \| \lambda_{t} - \lambda  \|^2 \Big ).
\end{equation} 
Summing  \eqref{eq:comp_t0} over $t=1,\cdots, K$, taking expectations, and plugging   \eqref{eq:comp_t1}, \eqref{eq:comp_t2}, and \eqref{eq:comp_t3} in, we obtain that for all $\lambda \in \RR_+^m$,
\begin{equation*}
\begin{split}
&\EE \Big [   \sum_{t=1}^{K} (\lambda  - \lambda_{t})^\top  \cL_G(x_t, v_{1,t},v_{2,t} ) \Big ] \\
 & \leq    \frac{\alpha }{2} \EE [ \| \lambda_{0}- \lambda  \|^2 ] -  \frac{\alpha }{2} \EE [ \| \lambda_{K} - \lambda  \|^2 ]  - \sum_{t=1}^K \frac{\alpha }{2} \EE[ \| \lambda_{t-1} - \lambda_{t} \|^2 ]  + \sqrt{K} M_{\lambda}\sigma_H    \\
 & \quad   + \sum_{t=1}^K   \Big (   \frac{ \alpha  \EE [\| \lambda_{t-1}  - \lambda_{t}\|^2 ] }{4} +  \frac{\sigma_{H}^2}{\alpha } +   \frac{5C_{g_1}^2 C_{g_2}^2 \EE[  \|  \lambda  - \lambda_{t} \|^2 ]  }{2\eta } +  \frac{\eta \EE [ \| x_{t} - x_{t-1}\|^2 ] }{10} \Big ) 
  \\
 &  \leq    \frac{\alpha }{2} \| \lambda_{0}- \lambda  \|^2   -   \frac{\alpha }{2} \EE [    \| \lambda_{K} - \lambda  \|^2 ]  +  \sqrt{K} M_{\lambda}  \sigma_H \\
 & \quad  + \sum_{t=1}^K\Big (  \frac{\sigma_{H}^2}{\alpha }  - \frac{ \alpha \EE[ \| \lambda_{t-1}  - \lambda_{t}\|^2]}{4}   + \frac{5C_{g_1}^2 C_{g_2}^2 \EE[  \|  \lambda  - \lambda_{t} \|^2 ]  }{2\eta } + \frac{\eta \EE [ \| x_{t} - x_{t-1}\|^2 ] }{10} \Big )  , 
 \end{split}
\end{equation*} 
which completes the proof.
\QED

\subsection{Proof of Theorem~\ref{thm:composition_01}} 
We provide bound on each term after decomposing \eqref{eq:decom_comp}. 
Consider any integer $K$ such that $K \leq N$. 
First, by  Lemma \ref{lemma:pi_ast_t_compositional}  in Appendix Section \ref{app:B}, we obtain 
\begin{equation}\label{eq:comp_lambda_0}
\begin{split}
  & \sum_{t=1}^K \EE  \left [ \lambda^\top  \Big ( \cL_G(x_t, v_{1:2} ) - \cL_G(x_t, v_{1:2 , t} ) \Big ) \right ] \\
& \quad\leq    \sum_{t=1}^{K}  \frac{\eta }{10} \EE \Big [ \| x_{t}-x_{t-1}  \|^2 \Big ]
+     \frac{40 K C_{g_1}^2 C_{g_2}^2  M_{\lambda}^2 }{ \eta }  +      2\sqrt{ K } L_{g_1} \sigma_{g_2}^2 M_{\lambda} +   3 \sqrt{K} C_{g_1}\sigma_{g_2} M_{\lambda} .
\end{split}
\end{equation}
Second, by  Lemma \ref{lemma:comp_lambda_F} in Appendix Section \ref{app:B}, we have
\begin{equation}\label{eq:comp_lambda_1}
\begin{split}
&\sum_{t=1}^K \EE \Big [  \cL_F(x_t, \pi_{1:2,t}   ) -    \cL_F(x^*, \pi_{1:2,t}  ) + \lambda_t^\top \Big (  \cL_G(x_t, v_{1:2,t} )   - \cL_G(x^*, v_{1:2,t} )  \Big ) \Big]   \\
& \quad\leq    \frac{\eta }{2} \| x_{0}- x^*\|^2  -  \sum_{t=1}^K  \frac{3\eta }{10}    \EE \Big [  \|x_{t-1} - x_{t} \|^2  
 \Big ]      +  \sum_{t=1}^K \frac{5 C_{f_1}^2 C_{f_2}^2 }{2\eta   }   \\
& 
\qquad +  \sum_{t=1}^K  \frac{\alpha }{4} \EE \Big [  \| \lambda_t - \lambda_{t-1}\|^2  \Big ] + \sum_{t=1}^K  \frac{ D_X^2 C_{g_1}^2 C_{g_2}^2}{\alpha }
+ \sum_{t=1}^K \frac{5C_{g_1}^2 C_{g_2}^2 \EE [ \| \lambda_{t-1}\|^2]  }{2\eta }.
\end{split} 
\end{equation}
Meanwhile, Lemma \ref{lemma:pi_ast_t}  in Appendix Section \ref{app:A1} implies that 
\begin{equation}\label{eq:comp_lambda_2}
\begin{split}
& \sum_{t=1}^K \EE \Big [ \cL_F(x_t, \pi_{1:2})  - \cL_F(x_t, \pi_{1:2,t} )   \Big ]  \\
& \quad\leq   \sum_{t=1}^{K}  \frac{\eta }{10} \EE \Big  [ \| x_{t}-x_{t-1}  \|^2 \Big ]
+     \frac{40 K C_{f_1}^2 C_{f_2}^2 }{ \eta }  +      2\sqrt{ K } L_{f_1} \sigma_{f_2}^2  +   3 \sqrt{K} C_{f_1}\sigma_{f_2} .
\end{split}
\end{equation}
Finally, by  Lemma \ref{lemma:comp_lambda_G} and the assumption that $\| \lambda\| \leq M_{\lambda}$,  we have 
\begin{equation}\label{eq:lambda_LG}
\begin{split}
&\EE \Big [   \sum_{t=1}^{K} (\lambda  - \lambda_{t})^\top  \cL_G(x_t, v_{1:2 , t} ) \Big ]  +  \frac{\alpha }{2} \EE \Big [    \| \lambda_{K} - \lambda  \|^2 \Big ]   \\
 & \quad \leq  \sum_{t=1}^K\Big (  \frac{\sigma_{H}^2}{\alpha }  - \frac{ \alpha \EE[ \| \lambda_{t-1}  - \lambda_{t}\|^2] }{4}   + \frac{5C_{g_1}^2 C_{g_2}^2 \EE[  \| \lambda  - \lambda_{t} \|^2 ] }{2\eta } + \frac{\eta \EE [ \| x_{t} - x_{t-1}\|^2 ] }{10} \Big )  \\
 & \qquad  + \frac{\alpha }{2}\EE \Big [  \| \lambda_{0}- \lambda  \|^2 \Big ]    +  \sqrt{K} M_{\lambda} \sigma_H    .
 \end{split}
\end{equation} 
By substituting \eqref{eq:comp_lambda_0}, \eqref{eq:comp_lambda_1}, \eqref{eq:comp_lambda_2}, and \eqref{eq:lambda_LG} into \eqref{eq:decom_comp}, and omitting some algebras, 
 we obtain 
\begin{equation*}
\begin{split}
& \sum_{t=1}^K \EE \Big ( \cL(x_t,\lambda,  \pi_{1:2} , v_{1:2}) -  \cL(x^*, \lambda_t,  \pi_{1:2,t}, v_{1:2,t})   \Big )  +  \frac{\alpha }{2} \EE [    \| \lambda_{K} - \lambda  \|^2 ] 
\\
& \quad \leq  
\frac{40 K C_{g_1}^2 C_{g_2}^2  M_{\lambda}^2 }{ \eta }  +      2\sqrt{ K } L_{g_1} \sigma_{g_2}^2 M_{\lambda} +   3 \sqrt{K} C_{g_1}\sigma_{g_2} M_{\lambda} + \frac{\eta }{2} \| x_{0}- x^*\|^2  
\\
& 
\qquad    +   \frac{95 K C_{f_1}^2 C_{f_2}^2 }{2\eta   }    +  \frac{K}{\alpha } ( D_X^2 C_{g_1}^2 C_{g_2}^2 + \sigma_H^2 )
+ \frac{\alpha }{2}\EE[  \| \lambda_{0}- \lambda  \|^2]    +  \sqrt{K} M_{\lambda}  \sigma_H  
\\
& 
\qquad 
 +      2\sqrt{ K } L_{f_1} \sigma_{f_2}^2  +   3 \sqrt{K} C_{f_1}\sigma_{f_2} 
+  \sum_{t=1}^K  \frac{5  C_{g_1}^2 C_{g_2}^2 }{2\eta }  \Big ( \EE[  \| \lambda  - \lambda_{t} \|^2] + \EE [ \| \lambda_{t-1}\|^2]  \Big ).
\end{split}
\end{equation*}
This completes the proof. 
\QED

\begin{lemma}\label{lemma:pi_ast_t_compositional} 
Suppose Assumptions \ref{assumption:01}, \ref{assumption:03}, and  \ref{assumption:compositional_g} hold, and Algorithm \ref{alg:02} generates 
  $ \{ (x_t,\lambda_t,\pi_{1:2,t}, v_{1:2,t}) \}_{t=1}^N$ by setting $\tau_t = \rho_t = t/2$, $\eta_t =  \eta $, and $\alpha_t = \alpha $ for $t=1,2,\cdots, N$. Letting $\lambda \in \RR_+^m$ be a bounded random variable satisfying $\| \lambda \| \leq M_{\lambda}$ uniformly, then for any $K \leq N$, we have
\begin{equation*}
\begin{split}
& \sum_{t=1}^K \EE \Big [\lambda^\top \Big ( \cL_G(x_t, v_{1:2})  - \cL_G(x_t, v_{1:2,t} )   \Big ) \Big ]  \\
&\quad \leq   \sum_{t=1}^{K}  \frac{\eta }{10} \EE [ \| x_{t}-x_{t-1}  \|^2 ]
+     \frac{60 K C_{g_1}^2 C_{g_2}^2 M_{\lambda}^2 }{ \eta }  +      2\sqrt{K}   L_{g_1} M_{\lambda} \sigma^2_{g_2}  + 4 \sqrt{K} M_\lambda C_{g_1} \sigma_{g_2}.
\end{split}
\end{equation*}
\end{lemma}
\textit{Proof:}
This result is similar to Lemma \ref{lemma:pi_ast_t} with an additional nonnegative (random) variable $\lambda$. 
We start with the following decomposition that
\begin{equation*}
\begin{split}
& \lambda^\top (\cL_G(x_t, v_{1:2} )  - \cL_G(x_t, v_{1:2,t} ) )  
\\
& \quad = \lambda^\top (\cL_{g_1}(x_t, v_1,v_2 )  -  \cL_{g_1}(x_t, v_1,v_{2,t} ) )  + \lambda^\top ( \cL_{g_1}(x_t, v_1,v_{2,t} )  - \cL_{g_1}(x_t, v_{1,t},v_{2,t} ) ).
\end{split}
\end{equation*}
 First, we recall \eqref{def:v2} that $
v_{2, t} \in \argmax_{ v_2 \in \cV_2} \{ v_2^\top  x_{t-1}- g_2^*( v_2)  \}$. By  Assumption~\ref{assumption:compositional_g} (g) that for each compositional constraint $G^{(j)} = g_1^{(j)} \circ g_2^{(j)}(x)$, the outer-level function $g_1^{(j)}$ is monotone nondecreasing for non-affine inner-level function $g_2^{(j)}$, we obtain 
$$
\EE \Big [ \inner{ ( v_2 - v_{2,t})^\top x_{t -1} - g_2^*(v_2)  + g_2^*(v_{2,t}) } {v_1   }  \Big ]  \leq 0.
$$
Because $\lambda \geq 0$, by following \eqref{eq:Lf_1} in Lemma \ref{lemma:pi_ast_t}, we have 
\begin{align*} 
& \sum_{t=1}^K \EE \Big [ \lambda^\top \Big (  \cL_{g_1}(x_t, v_1,v_2 )  -  \cL_{g_1}(x_t, v_1,v_{2,t} )    \Big )  \Big ] 
\leq  \sum_{t=1}^{K} \EE [ \langle ( v_2 - v_{2, t})  v_1 \lambda , x_{t} - x_{t-1} \rangle ] \\
&
\quad \leq \sum_{t=1}^{K}   \Big [  \frac{10 \EE [ \|   (v_2 - v_{2,t})  v_1 \lambda \|^2]  }{\eta }  +   \frac{\eta }{20} \EE [ \| x_{t}-x_{t-1}  \|^2 ]\Big ] 
 \\
 & \quad \leq \sum_{t=1}^{K}   \Big [  \frac{40  C_{g_1}^2 C_{g_2}^2  M_{\lambda}^2}{\eta }  +   \frac{\eta }{20} \EE [ \| x_{t}-x_{t-1}  \|^2 ]\Big ],
\end{align*} 
where the last inequality uses the fact that $\norm{\lambda} \leq M_{\lambda}$. Second, by  Proposition 15 in \cite{zhang2024optimal}, the $L_{g_1}$ smoothness of $g_1$ implies 
$\lambda^\top D_{g_1^*}(v_1; \bar{v}_1)\geq \tfrac{\|\lambda^\top (v_1 - \bar{v}_1)\|^2}{ 2\norm{\lambda}  L_{g_1}  } \geq \tfrac{\|\lambda^\top (v_1 - \bar{v}_1)\|^2}{ 2M_\lambda L_{g_1}  }$, 
which further yields
\begin{align*}
&\sum_{t=1}^K \EE  \left[ \lambda^\top (v_{1,t} - v_{1,t-1})^\top (g_2(x_{t-1}, \zeta_{2,t-1}^1) - g_2(x_{t-1})) - \rho_t \lambda^\top D_{g_1^*}(v_{1,t}; v_{1,t-1})\right] \\
 &\quad\leq \sum_{t=2}^{K} \frac{M_\lambda L_{g_1}\sigma_{g_2}^2  }{\rho_t} + \sqrt{\EE [  \|\lambda\|^2 \norm{v_{1,1} - v_{1,0}}^2] }\sqrt{\EE [\norm{g_2(x_{0}, \zeta_{2,0}^1) - g_2(x_{0})}^2]}\\
 &\quad\leq 2 L_{g_1} M_{\lambda} \sigma^2_{g_2}  \log K + 2 M_\lambda C_{g_1} \sigma_{g_2}\\
 &\quad\leq 2\sqrt{K}   L_{g_1} M_{\lambda} \sigma^2_{g_2}  + 2 \sqrt{K} M_\lambda C_{g_1} \sigma_{g_2}.
 \end{align*} 
 By substituting the above bound into Proposition 6 of \cite{zhang2024optimal}, we further obtain that 
\begin{align*}
& \sum_{t=1}^K \EE  \Big [ \lambda^\top \Big (  \cL_{g_1}(x_t, v_1,v_{2,t} )  - \cL_{g_1}(x_t, v_{1,t},v_{2,t} ) \Big )  \Big ] \\
&\quad \leq  \sum_{t=1}^{K}   \Big [  \frac{20 C_{g_1}^2 C_{g_2}^2  M_{\lambda}^2 }{\eta }  +   \frac{\eta }{20} \EE [ \| x_{t}-x_{t-1}  \|^2 ]\Big ] 
+   2\sqrt{K}   L_{g_1} M_{\lambda} \sigma^2_{g_2}  + 4 \sqrt{K} M_\lambda C_{g_1} \sigma_{g_2}.
\end{align*}
The desired result then follows by combining the preceding inequalities.
\QED 

\begin{lemma}\label{lemma:comp_lambda_F}
Suppose Assumptions \ref{assumption:01}, \ref{assumption:03}, and \ref{assumption:compositional_g} hold, and  Algorithm \ref{alg:02} generates $\{ (x_t,\lambda_t,\pi_t, v_t) \}_{t=1}^N $ by setting $\rho_t = \tau_t = \frac{t}{2}$, $\eta_{t} = \eta$, and $\alpha_t = \alpha$ for $t \leq N$. Then for any integer $K \leq N$, we have 
\begin{equation*}
\begin{split}
  &\sum_{t=1}^K \EE \Big [  \cL_F(x_t,  \pi_{1:2,t}   ) -    \cL_F(x^*, \pi_{1:2,t}  ) + \lambda_t^\top \Big (  \cL_G(x_t, v_{1:2,t} )   - \cL_G(x^*, v_{1:2,t} )  \Big ) \Big]   \\
&\quad \leq    \frac{\eta }{2} \| x_{0}- x^*\|^2  -  \sum_{t=1}^K  \frac{3\eta_{t-1}}{10}    \EE \Big [  \|x_{t-1} - x_{t} \|^2  
 \Big ]      +  \sum_{t=1}^K \frac{5 C_{f_1}^2 C_{f_2}^2 }{2\eta  }   \\
& 
\qquad +  \sum_{t=1}^K  \frac{\alpha }{4} \EE[  \| \lambda_t - \lambda_{t-1}\|^2 ] + \sum_{t=1}^K  \frac{ D_X^2 C_{g_1}^2 C_{g_2}^2
}{\alpha }+ \sum_{t=1}^K \frac{5C_{g_1}^2 C_{g_2}^2  \EE [ \| \lambda_{t-1}\|^2] }{2\eta }.
\end{split} 
\end{equation*}
\end{lemma}
\textit{Proof:}
Recall the update rule for $x_t$ \eqref{eq:update_x_comp_ssd}  that 
$$
x_t = \argmin_{x\in \cX} \Big \{ \big \langle  \pi_{2,t}^0   \pi_{1,t}^0 + v_{2,t}^0   v_{1,t}^0  \lambda_{t-1} , x  \big \rangle +   \frac{\eta }{2} \| x_{t-1}- x\|^2 \Big  \} .
$$
By following analogous analysis in Lemma \ref{lemma:F_lambda}, and applying the three-point Lemma \ref{lemma:three_point} to the above update rule, we have 
\begin{equation}\label{eq:comp_lambda_F_00}
\begin{split}
  & \cL_F(x_t, \pi_{1:2,t}  ) -    \cL_F(x^*, \pi_{1:2,t} ) + \lambda_t^\top \Big (  \cL_G(x_t, v_{1:2,t} )   - \cL_G(x^*, v_{1:2,t})  \Big )    \\ 
 &  \quad =     (x_t - x^*)^\top  \pi_{2,t}  \pi_{1,t}  +  (x_t - x^*)^\top   v_{2,t} v_{1,t}   \lambda_t\\
& \quad \leq
  \frac{\eta }{2} \| x_{t-1} - x^*\|^2 - \frac{\eta }{2} \|x_{t-1} - x_{t} \|^2 - \frac{\eta }{2} \|x_{t} -x^* \|  \\
  & \qquad \ + \langle \pi_{2,t}  \pi_{1,t}  -  \pi_{2,t}^0  \pi_{1,t}^0 , x_t - x^* \rangle  + \tilde \Lambda_{1,t}+ \tilde \Lambda_{2,t}, 
\end{split}
\end{equation}
  where 
  $$
  \tilde \Lambda_{1,t}  = (x_t - x^*)^\top   v_{2,t}  v_{1,t} ( \lambda_t -\lambda_{t-1} ) \mbox{  and  } \ \tilde \Lambda_{2,t} = (x_t - x^*)^\top  \big  (  v_{2,t}  v_{1,t}- v_{2,t}^0   v_{1,t}^0   \big  )\lambda_{t-1} . 
  $$
We note that a bound for $\inner{ \pi_{2,t}   \pi_{1,t}- \pi_{2,t}^0  \pi_{1,t}^0 }{ x_t - x^*}$  has been provided in \eqref{eq:Delta_0}. Specifically, we have 
\begin{equation}\label{eq:comp_tilde_delta_0}
\begin{split}
\EE \Big [ \sum_{t=1}^K \inner{  \pi_{2,t} \pi_{1,t}  - \pi_{2,t}^0    \pi_{1,t}^0  }{  x_{t} - x^* } \Big ] & = \EE \Big [ \sum_{t=1}^K \inner{  \pi_{2,t} \pi_{1,t}  - \pi_{2,t}^0    \pi_{1,t}^0  }{  x_{t} - x_{t-1}  } \Big ]
\\
&  \leq   \sum_{t=1}^K  \Big ( \frac{5 C_{f_1}^2 C_{f_2}^2 }{2 \eta } + \frac{ \eta  \EE[ \| x_t - x_{t-1}\|^2]}{10 } \Big ).
 \end{split}
\end{equation}
 To bound $\tilde \Lambda_{1,t}$, by following the analysis  of
 \eqref{eq:Lambda_1}, we have 
 \begin{equation}\label{eq:comp_tilde_delta_1}
\begin{split}
\EE [ \tilde \Lambda_{1,t}] & \leq  \frac{\alpha }{4} \EE[  \| \lambda_t - \lambda_{t-1}\|^2 ] + \frac{1}{ \alpha } \EE [ \|(x_t - x^*)^\top   v_{2,t} v_{1,t}  \|^2] \\
& \leq   \frac{\alpha }{4} \EE[  \| \lambda_t - \lambda_{t-1}\|^2 ] + \frac{D_X^2 C_{g_1}^2 C_{g_2}^2}{\alpha } .
\end{split}
\end{equation}
Meanwhile, we observe from Algorithm~\ref{alg:02} that 
$\EE[ x_{t-1}^\top (v_{2,t}^0v_{1,t}^0 - v_{2,t} v_{1,t}) \lambda_{t-1}] = 0$. 
By following the analysis of \eqref{eq:Lambda_2}, we have 
 \begin{equation}\label{eq:comp_tilde_delta_2}
\begin{split}
\EE[\sum_{t=1}^K \tilde \Lambda_{2,t}] 
= & \sum_{t=1}^K \EE[ (x_{t} - x_{t-1})^\top  ( v_{2,t} v_{1,t} - v_{2,t}^0 v_{1,t}^0  ) \lambda_{t-1} ] \\
\leq & \sum_{t=1}^K \EE[ \frac{5\| (v_{2,t}^0v_{1,t}^0 - v_{2,t} v_{1,t}) \lambda_{t-1}  \|^2}{2\eta  } ] + \sum_{t=1}^K  \frac{\eta  \EE[  \|x_{t} - x_{t-1} \|^2] }{10}
\\
\leq & \sum_{t=1}^K \frac{5C_{g_1}^2 C_{g_2}^2  \EE [ \| \lambda_{t-1}\|^2] }{2\eta  }   + \sum_{t=1}^K  \frac{\eta  \EE[   \|x_{t} - x_{t-1} \|^2] }{10},
\end{split}
\end{equation}
where the last inequality uses the fact that $\lambda_{t-1}$ and $v_{2,t}^0v_{1,t}^0 - v_{2,t} v_{1,t}  $ are independent. 
Summing \eqref{eq:comp_lambda_F_00} over $t=1,2,\cdots,K$, and combining \eqref{eq:comp_tilde_delta_0}, \eqref{eq:comp_tilde_delta_1}, and \eqref{eq:comp_tilde_delta_2}, we obtain
\begin{equation*}
\begin{split}
  &\sum_{t=1}^K \EE \Big [  \cL_F(x_t, \pi_{1,t},\pi_{2,t} ) -    \cL_F(x^*, \pi_{1,t},\pi_{2,t} ) + \lambda_t^\top \Big (  \cL_G(x_t, v_t )   - \cL_G(x^*, v_t)  \Big ) \Big]   \\
& \quad \leq    \frac{\eta }{2} \| x_{0}- x^*\|^2 
+  \sum_{t=1}^K \Big (  \frac{\eta }{10} +  \frac{\eta }{10}  - \frac{\eta }{2} \Big )  \EE \Big [  \|x_{t-1} - x_{t} \|^2  
 \Big ]       +  \sum_{t=1}^K \frac{5 C_{f_1}^2 C_{f_2}^2 }{2\eta  }    \\
& 
\qquad +  \sum_{t=1}^K  \frac{\alpha }{4} \EE[  \| \lambda_t - \lambda_{t-1}\|^2 ] + \sum_{t=1}^K  \frac{ D_X^2 C_{g_1}^2 C_{g_2}^2 }{\alpha }+ \sum_{t=1}^K  \frac{5  C_{g_1}^2 C_{g_2}^2 \EE[ \| \lambda_{t-1} \|^2 ]  }{2\eta }    \\
& \quad=    \frac{\eta }{2} \| x_{0}- x^*\|^2  -  \sum_{t=1}^K  \frac{3\eta }{10}    \EE \Big [  \|x_{t-1} - x_{t} \|^2  
 \Big ]      +  \sum_{t=1}^K \frac{5 C_{f_1}^2 C_{f_2}^2 }{2\eta  }   \\
& 
\qquad +  \sum_{t=1}^K  \frac{\alpha }{4} \EE[  \| \lambda_t - \lambda_{t-1}\|^2 ] + \sum_{t=1}^K  \frac{D_X^2 C_{g_1}^2 C_{g_2}^2
}{\alpha } + \sum_{t=1}^K \frac{5 C_{g_1}^2 C_{g_2}^2 \EE [ \| \lambda_{t-1}\|^2] }{2\eta }.
\end{split} 
\end{equation*}
This completes the proof. 
\QED

\subsection{Proof of Proposition \ref{prop:dual_comp}}
For the feasible point $(x^*,\lambda, \pi_{1:2}, v_{1:2})$, 
we set $\pi_{1:2} = \pi_{1:2}^*$, $v_{1:2} = v_{1:2}^*$, and  $\lambda = \lambda^*$ as defined in \eqref{def:saddle_2_ssd}. By applying  Theorem \ref{thm:composition_01} with $M_{\lambda^*} = \| \lambda^*\|$, we obtain
\begin{equation*}
\begin{split}
& \sum_{t=1}^K \EE \Big ( \cL(x_t,\lambda^*,  \pi_{1:2}^*  , v_{1:2}^*  ) -  \cL(x^*, \lambda_t,  \pi_{1:2,t} , v_{1:2,t})   \Big )  +  \frac{\alpha }{2} \EE [    \| \lambda_{K} - \lambda^*  \|^2 ] 
\\
& \leq  \frac{40 K C_{g_1}^2 C_{g_2}^2 \|\lambda^* \|   }{ \eta }  +      2\sqrt{ K } L_{g_1} \sigma_{g_2}^2 \|\lambda^* \| +   3 \sqrt{K} C_{g_1}\sigma_{g_2} \|\lambda^* \| 
+ \frac{\eta }{2} \| x_{0}- x^*\|^2   \\
& 
\quad    +   \frac{95 K C_{f_1}^2 C_{f_2}^2 }{2\eta   }    +  \frac{K}{\alpha } ( D_X^2 C_{g_1}^2 C_{g_2}^2 + \sigma_H^2 )
+ \frac{\alpha }{2} \| \lambda_{0}- \lambda^*  \|^2    +  \sqrt{K} \|\lambda^* \|    \sigma_H   
\\
& 
\quad 
 +      2\sqrt{ K } L_{f_1} \sigma_{f_2}^2  +   3 \sqrt{K} C_{f_1}\sigma_{f_2} 
+  \sum_{t=1}^K  \frac{5  C_{g_1}^2 C_{g_2}^2 }{2\eta }  \Big ( \EE[  \| \lambda^*  - \lambda_{t}\|^2] + \EE [ \| \lambda_{t-1}\|^2]  \Big ).
\end{split}
\end{equation*}
We note that $\cL(x_t,\lambda^*,  \pi_{1:2}^*  , v_{1:2}^*  ) -  \cL(x^*, \lambda_t,  \pi_{1:2,t} , v_{1:2,t}  )   \geq 0$ by the definition of a saddle point~\eqref{def:saddle_2_ssd}. Thus, we have 
\begin{equation*}
\begin{split}
&  \frac{\alpha }{2} \EE [   \| \lambda_{K} - \lambda^*  \|^2 ] \\
&\leq  \frac{40 K C_{g_1}^2 C_{g_2}^2 \|\lambda^* \|   }{ \eta }  
+ \frac{\eta }{2} \| x_{0}- x^*\|^2     +   \frac{95 K C_{f_1}^2 C_{f_2}^2 }{2\eta   }    +  \frac{K}{\alpha } ( D_X^2 C_{g_1}^2 C_{g_2}^2 + \sigma_H^2 )
\\
& \quad  + \frac{\alpha }{2} \| \lambda_{0}- \lambda^*  \|^2  
 +  \sum_{t=1}^K  \frac{5 C_{g_1}^2 C_{g_2}^2  }{2\eta }  \Big ( \EE[  \| \lambda^*  - \lambda_{t}\|^2] + \EE [ \| \lambda_{t-1}\|^2]  \Big ) 
+ \sqrt{K} Q, 
\end{split}
\end{equation*}
where 
$Q =    2L_{g_1} \sigma_{g_2}^2 \|\lambda^* \| +   3C_{g_1}\sigma_{g_2} \|\lambda^* \|  +  \|\lambda^* \|   \sigma_H      +        2L_{f_1} \sigma_{f_2}^2  +   3 C_{f_1}\sigma_{f_2}  $ is defined in \eqref{eq:R_comp}.
By using the fact that $\| \lambda_t  \|^2 \leq 2 \| \lambda_t  - \lambda^* \|^2 + 2 \|  \lambda^* \|^2  $, setting $\alpha = 2 \sqrt{N}$ and $\eta = \frac{15 C_{g_1}^2 C_{g_2}^2 \sqrt{N}  }{2}$ , and dividing both sides of the above inequality by  $\sqrt{N}$, we obtain 
\begin{equation*}
\begin{split}
 &\EE [  \| \lambda_{K} - \lambda^*  \|^2 ]
\\
& \leq \frac{16 K   \|\lambda^* \|   }{ 3N}  
+ \frac{15 C_{g_1}^2 C_{g_2}^2 }{4} \| x_{0}- x^*\|^2     +   \frac{19 K C_{f_1}^2 C_{f_2}^2 }{3C_{g_1}^2 C_{g_2}^2 N }    +  \frac{K}{2N} ( D_X^2 C_{g_1}^2 C_{g_2}^2 + \sigma_H^2 )
\\
& \quad  + \| \lambda_{0}- \lambda^*  \|^2  
 +  \frac{ 1 }{3N}  \sum_{t=1}^K \Big ( \EE[  \| \lambda^*  - \lambda_{t}\|^2] +2 \EE [ \| \lambda_{t-1} - \lambda^*\|^2]  + 2\| \lambda^* \|^2  \Big )  + \frac{\sqrt{K}Q}{\sqrt{N} }
 \\
 & = 
R_K 
 + \frac{1}{3N}  \EE[  \| \lambda^*  - \lambda_{K}\|^2]  +  \frac{1}{N}  \sum_{t=1}^K \EE [ \| \lambda_{t-1} - \lambda^*\|^2],
\end{split}
\end{equation*}
where \begin{align*}
R_K & = \frac{K}{N} \Big (   6  \|\lambda^* \|      + \frac{ 7C_{f_1}^2 C_{f_2}^2 }{C_{g_1}^2 C_{g_2}^2}  +  \frac{ D_X^2 C_{g_1}^2 C_{g_2}^2 + \sigma_H^2  }{2}  \Big )   +  \frac{15 C_{g_1}^2 C_{g_2}^2 }{4}  \| x_{0}- x^*\|^2     + \| \lambda_{0}- \lambda^*  \|^2      + \frac{\sqrt{K}Q}{\sqrt{N} }.
\end{align*}
By noting that $R_{K} \leq R$ defined in \eqref{eq:R_comp}, the above inequality further implies that 
\begin{equation*}
\begin{split}
 \Big (1-\frac{1}{3N} \Big  )\EE [  \| \lambda_{K} - \lambda^*  \|^2 ] \leq R + \sum_{t=1}^K \frac{\EE [ \| \lambda_{t-1} - \lambda^*\|^2 ] }{N}.
\end{split}
\end{equation*}
Assuming 
$N \geq 2$, for all $K \leq N$, we have 
\begin{equation*}
\begin{split}
\EE [  \| \lambda_{K} - \lambda^*  \|^2 ] \leq 2R + \sum_{t=0}^{K-1} \frac{2}{N} \cdot\EE [ \| \lambda_{t} - \lambda^*\|^2].
\end{split}
\end{equation*}
Finally, by applying Lemma \ref{lemma:recursive_bound}, we recursively bound $\EE[ \| \lambda_{K} - \lambda^*  \|^2 ]$ by 
$$
\EE [  \| \lambda_{K} - \lambda^*  \|^2 ]  \leq 2R e^{2}
$$
for all $K \leq N$, which completes the proof. 
\QED 

\subsection{Proof of Theorem~\ref{thm:comp_rate}} We set $\lambda_0 = 0$ throughout our analysis. We first consider the objective optimality gap $F(\bar x_N) - F(x^*)$. We denote by 
\begin{equation*}
\hat \pi_2 \in  \partial  f_2(\bar x_N) , \, \hat \pi_1=  \nabla f_1\big ( f_2( \bar x_N)\big ), \, \hat v_2 \in \partial g_2(\bar x_N), \mbox{ and } \, \hat v_1 = \nabla g_1( g_2( \bar x_N) ). 
\end{equation*}
Setting the feasible point as $(x^*, 0, \hat \pi_{1:2}, \hat v_{1:2})$ and following Theorem~\ref{thm:rate} with $M_{\lambda} = 0$ for $\lambda = 0$, $\alpha_{t} = 2 \sqrt{N}$, and $\eta_t = \frac{15 C_{g_1}^2 C_{g_2}^2 \sqrt{N}}{2}$, we obtain that 
\begin{equation*}
\begin{split}
& \EE [ F(\bar x_N) ]- F(x^*) \leq   \frac{1}{N} \sum_{t=1}^N \EE \Big ( \cL(x_t, 0 , \hat \pi_{1:2} , \hat v_{1:2}) -  \cL(x^*, \lambda_t,  \pi_{1:2, t} , v_{1:2 , t} )  \Big )  \\
& \quad  \leq     \frac{1}{\sqrt{N}}\Big ( 
 \frac{15 C_{g_1}^2 C_{g_2}^2}{4} \| x_{0}- x^*\|^2     +   \frac{ 7  C_{f_1}^2 C_{f_2}^2  }{C_{g_1}^2 C_{g_2^2}}  +  \frac{ D_X^2 C_{g_1}^2 C_{g_2}^2 + \sigma_H^2}{2 }   +      2 L_{f_1} \sigma_{f_2}^2  +   3 C_{f_1}\sigma_{f_2}  \Big ) \\
& 
\qquad 
+  \frac{1}{N}\sum_{t=1}^N  \frac{ 1  }{3\sqrt{N}}\Big ( \EE[  \|  \lambda_{t}\|^2 ]+ \EE [ \| \lambda_{t-1}\|^2]  \Big ).
\end{split}
\end{equation*}
We then derive the convergence rate of the optimality gap by the bound of $\| \lambda_t\|^2$ in Proposition~\ref{prop:dual_comp} that 
\begin{equation*}
\EE[\| \lambda_t\|^2] \leq 2 \| \lambda^*\|^2 + 2 \EE[ \| \lambda_t - \lambda^*\|^2] \leq   2\| \lambda^*\|^2  + 2R e^{2}, \text{ for all }t=1,2,\cdots, N.
\end{equation*} 
Second, we consider  the feasibility residual $G(\bar x_N)_+$. By adopting the feasible  point $(x^*, \tilde \lambda, \hat \pi_{1:2}, \hat v_{1:2})$ where $\tilde \lambda = (\lambda^*+1)\frac{G(\bar x_N)_+}{ \| G(\bar x_N)_+ \|_2}$ and following Theorem~\ref{thm:rate} with $M_{\tilde \lambda } = \| \tilde \lambda\|$, $\alpha_{t} = 2 \sqrt{N}$, and $\eta_t = \frac{15 C_{g_1}^2 C_{g_2}^2 \sqrt{N}}{2}$, we obtain 
\begin{equation*}
\begin{split}
&   \EE[ \| G(\bar x_N)_{+} \|_2 ] 
 \leq    \frac{1}{N} \sum_{t=1}^N  \EE \Big ( \cL(x_t,  \tilde \lambda, \hat \pi_{1:2}, \hat v_{1:2}) - \cL(x^*, \lambda_t, \pi_{1:2,t}, v_{1:2,t})  \Big )  
\\
&\quad \leq \frac{1}{\sqrt{N}} \Big ( \frac{ 16    \| \tilde \lambda\|^2 }{3}  +      2L_{g_1} \sigma_{g_2}^2   \| \tilde \lambda\| +   3 \sqrt{K} C_{g_1}\sigma_{g_2}   \| \tilde \lambda\| +  \frac{15 C_{g_1}^2 C_{g_2}^2}{4} \| x_{0}- x^*\|^2      \Big )
\\
& 
\qquad + \frac{1}{\sqrt{N}} \Big (    \frac{ 7  C_{f_1}^2 C_{f_2}^2  }{C_{g_1}^2 C_{g_2^2}}    +     \frac{D_X^2 C_{g_1}^2 C_{g_2}^2 + \sigma_H^2}{2 } 
+ \| \tilde \lambda\|^2    +  \sigma_H   \| \tilde \lambda\|    +      2L_{f_1} \sigma_{f_2}^2  +   3  C_{f_1}\sigma_{f_2}   \Big )
\\
& 
\qquad 
+  \frac{1}{N}\sum_{t=1}^N  \frac{1 }{3\sqrt{N} }  \Big ( \EE[  \| \tilde \lambda  - \lambda_{t}\|^2 ]+ \EE [ \| \lambda_{t-1}\|^2]  \Big ).
 \end{split}
\end{equation*}
We obtain the convergence rate of feasibility residual $\| G(x)_+ \|_2$ by using the facts that 
\begin{equation*}
\| \lambda_t\|^2 \leq 2 \| \lambda^*\|^2 + 2 \| \lambda_t - \lambda^*\|^2 \text{ and } \| \tilde \lambda  - \lambda_{t}\|^2  \leq 3 \| \tilde \lambda  \|^2 + 3\| \lambda^* \|^2  + 3 \| \lambda_t - \lambda^*\|^2,
\end{equation*} 
and applying $\EE[\| \lambda_t - \lambda^*\|^2] \leq 2R e^{2}$ provided by Proposition~\ref{prop:dual_comp}. This completes the proof.  \QED

\end{document}